\documentclass{article}

\usepackage[T1]{fontenc}
\usepackage[french,english]{babel} 
\usepackage[utf8]{inputenc} 
\usepackage{amsmath,amssymb,amsthm}
\usepackage[colorlinks=true]{hyperref}
\usepackage{pdfsync}
\usepackage[dvips]{graphicx}
\usepackage{appendix}
\usepackage{algorithm2e}
\usepackage{float}

\def\R{\mathbb{R}}
\def\N{\mathbb{N}}

\def\T{\mathcal{T}}

\newtheorem{thm}{Theorem}
\newtheorem{lem}{Lemma}

\newtheorem{cor}{Corollary}
\theoremstyle{definition}
\theoremstyle{definition}
\theoremstyle{definition}\newtheorem{defi}{Definition}

\renewcommand{\geq}{\geqslant}
\renewcommand{\leq}{\leqslant}

\def\U{\mathcal U}
\def\G{\mathcal G}
\def\P{\mathcal P}

\title{Geometric Optimal Control and Applications to Aerospace}

\author{Jiamin Zhu\footnote{Sorbonne Universit\'es, UPMC Univ Paris 06, CNRS UMR 7598, Laboratoire Jacques-Louis Lions, F-75005, Paris, France (\texttt{zhu@ann.jussieu.fr}).}
\and
Emmanuel Tr\'elat\footnote{Sorbonne Universit\'es, UPMC Univ Paris 06, CNRS UMR 7598, Laboratoire Jacques-Louis Lions, F-75005, Paris, France (\texttt{emmanuel.trelat@upmc.fr}).}
\and
Max Cerf\footnote{Airbus Defence and Space, Flight Control Unit, 66 route de Verneuil, BP 3002, 78133 Les Mureaux Cedex, France (\texttt{max.cerf@astrium.eads.net}).}
}

\date{}

\begin{document}

\maketitle

\begin{abstract}
This survey article deals with applications of optimal control to aerospace problems with a focus on modern geometric optimal control tools and numerical continuation techniques.
Geometric optimal control is a theory combining optimal control with various concepts of differential geometry. The ultimate objective is to derive optimal synthesis results for general classes of control systems.
Continuation or homotopy methods consist in solving a series of parameterized problems, starting from a simple one to end up by continuous deformation with the initial problem. They help overcoming the difficult initialization issues of the shooting method.
The combination of geometric control and homotopy methods improves the traditional techniques of optimal control theory.

A nonacademic example of optimal attitude-trajectory control of (classical and airborne) launch vehicles, treated in details, illustrates how geometric optimal control can be used to analyze finely the structure of the extremals. This theoretical analysis helps building an efficient numerical solution procedure combining shooting methods and numerical continuation.
Chattering is also analyzed and it is shown how to deal with this issue in practice.
\end{abstract}

\vspace{0.5cm}


\newpage
\tableofcontents
\newpage

\section{Introduction} \label{sec_Intro}
This article makes a survey of the main issues in optimal control theory, with a specific focus on numerical solution methods and applications to aerospace problems.
The purpose is to show how to address optimal control problems using modern techniques of geometric optimal control and how to build solution algorithms based on continuation techniques.
The geometric optimal control (stated in the early 1980s and having widely demonstrated its advantages over the classical theory of the 1960s) and the continuation techniques (which are not new, but have been somewhat neglected until recently in optimal control) are powerful approaches for aerospace applications.

As motivation, an overview of optimal control problems raised by aerospace missions is first presented.
These problems are classified in three categories depending on the departure and the arrival point. The interested reader will thus have a general view on how space transportation missions translate into optimal control problems.

A detailed example is then presented to illustrate the application of geometric optimal control techniques and numerical continuation methods on a practical problem.
This example deals with a minimum time maneuver of a coupled attitude-trajectory dynamic system.
Due to the system high nonlinearity and the existence of a chattering phenomenon (see Sections \ref{chatphenom} and \ref{fullexample} for details), the standard techniques of optimal control do not provide adequate solutions to this problem.
Through this example, we will show step by step how to build efficient numerical procedures with the help of theoretical results obtained by applying geometric optimal control techniques.

Before this example, we will recall briefly the main techniques of optimal control theory, including the Pontryagin Maximum Principle, the first-order and higher order optimality conditions, the associated numerical methods, and the numerical continuation principles.
Most mathematical notions presented here are known by many readers, and can be skipped at the first reading.

\medskip

In Section \ref{sec_appliaero}, several optimal control problems stemming from various aerospace missions are systematically introduced.
In Section \ref{sec_GOCR}, we provide a brief survey of geometric optimal control, including the use of Lie and Poisson brackets with first and higher order optimality conditions. In Section \ref{sec_NMOC}, we recall classical numerical methods for optimal control problems, namely indirect and direct methods.
In Section \ref{sec_CM}, we recall the concept of continuation methods, which help overcoming the initialization issue for indirect methods.
In Section \ref{fullexample}, we detail a full nonacademic example in aerospace, in order to illustrate how to solve optimal control problems with the help of geometric optimal control theory and the continuation methods.
Finally in Section \ref{sec_ATO}, we shortly give other applications of geometric optimal control and of continuation for space trajectory optimization problems.

\section{Applications to Aerospace Problems}
\label{sec_appliaero}
Transport in space gives rise to a large range of problems that can be addressed by optimal control and mathematical programming techniques. Three kinds of problems can be distinguished depending on the departure and the arrival point: ascent from the Earth ground to an orbit, reentry from an orbit to the Earth ground (or to another body of the solar system), transfer from an orbit to another one. A space mission is generally composed of successive ascent, transfer and reentry phases, whose features are presented in the following paragraphs.

\begin{figure}[h]
\centering
	\includegraphics[scale=0.73]{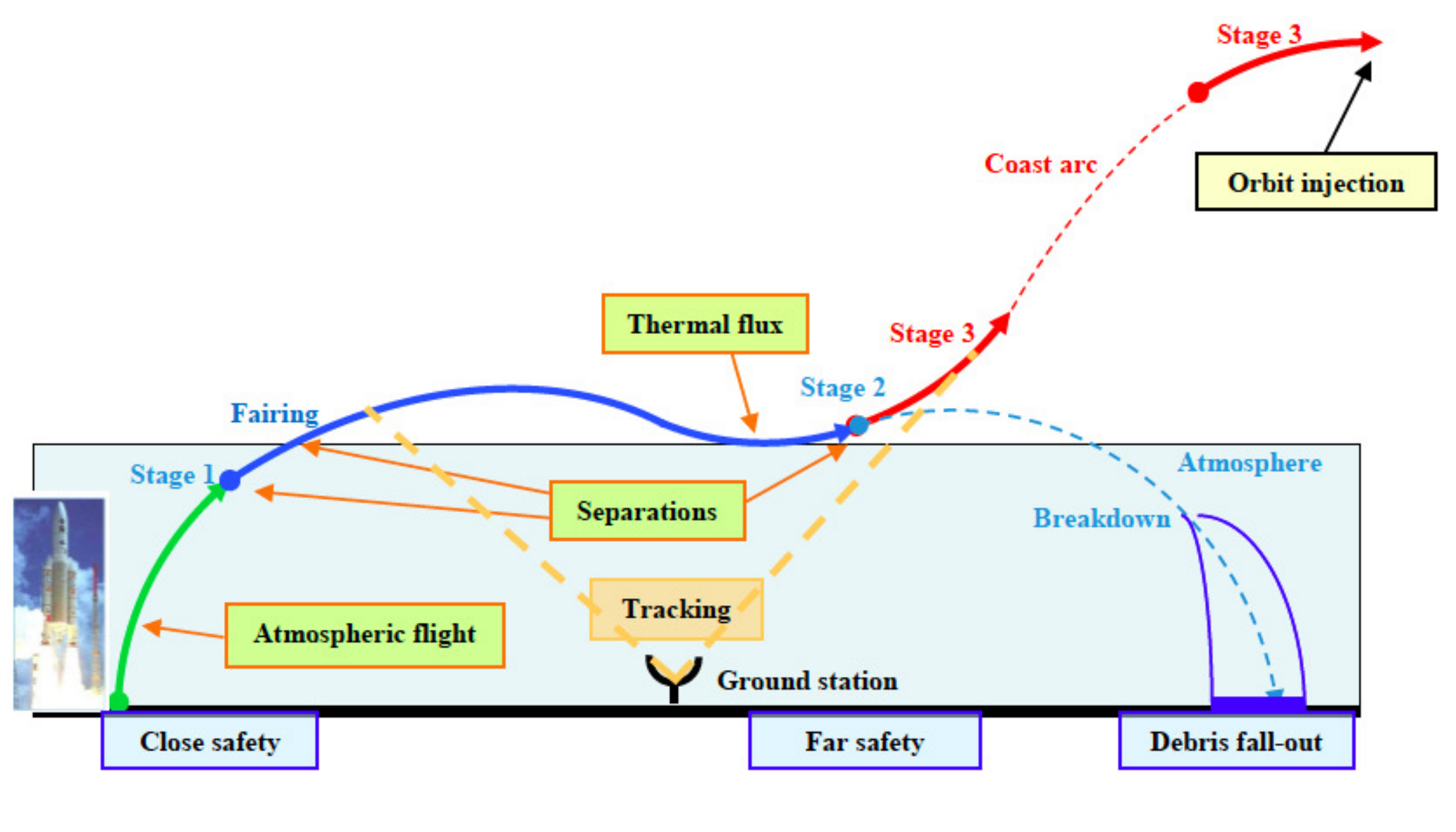}
	\caption{Ascent trajectory.}
	\label{ascenttraj}
\end{figure}

Ascent missions necessitate huge propellant masses to reach the orbital velocity and deliver large payloads such as telecommunications satellites. Due to the large lift-off mass, only chemical rocket engines are able to deliver the required thrust level. Consumption minimization is the main concern for these missions whose time of flight is generally about half an hour. Heavy launchers lift off vertically from a fixed ground launch pad, whereas airborne launchers are released horizontally by an airplane, benefiting thus from a higher initial altitude and an initial subsonic velocity. The first part of the trajectory occurs in the Earth atmosphere at increasing speed. The large aerodynamics loads met during the atmospheric flight require flying at near zero angle of attack, so that the atmospheric leg is completely driven by the initial conditions. Due to the large masses of propellants carried on board, the whole flight must be track by ground radar stations and stringent safety constraints must be applied regarding the area flown over. Once in vacuum the vehicle attitude is no longer constrained and the thrust direction can be freely chosen. When the orbital velocity is reached the thrust level can be reduced and coast arcs may help sparing propellant to reach the targeted orbit.
Figure \ref{ascenttraj} gives an overview of the constraints applied to an ascent trajectory.

\begin{figure}[H]
\centering
	\includegraphics[scale=0.73]{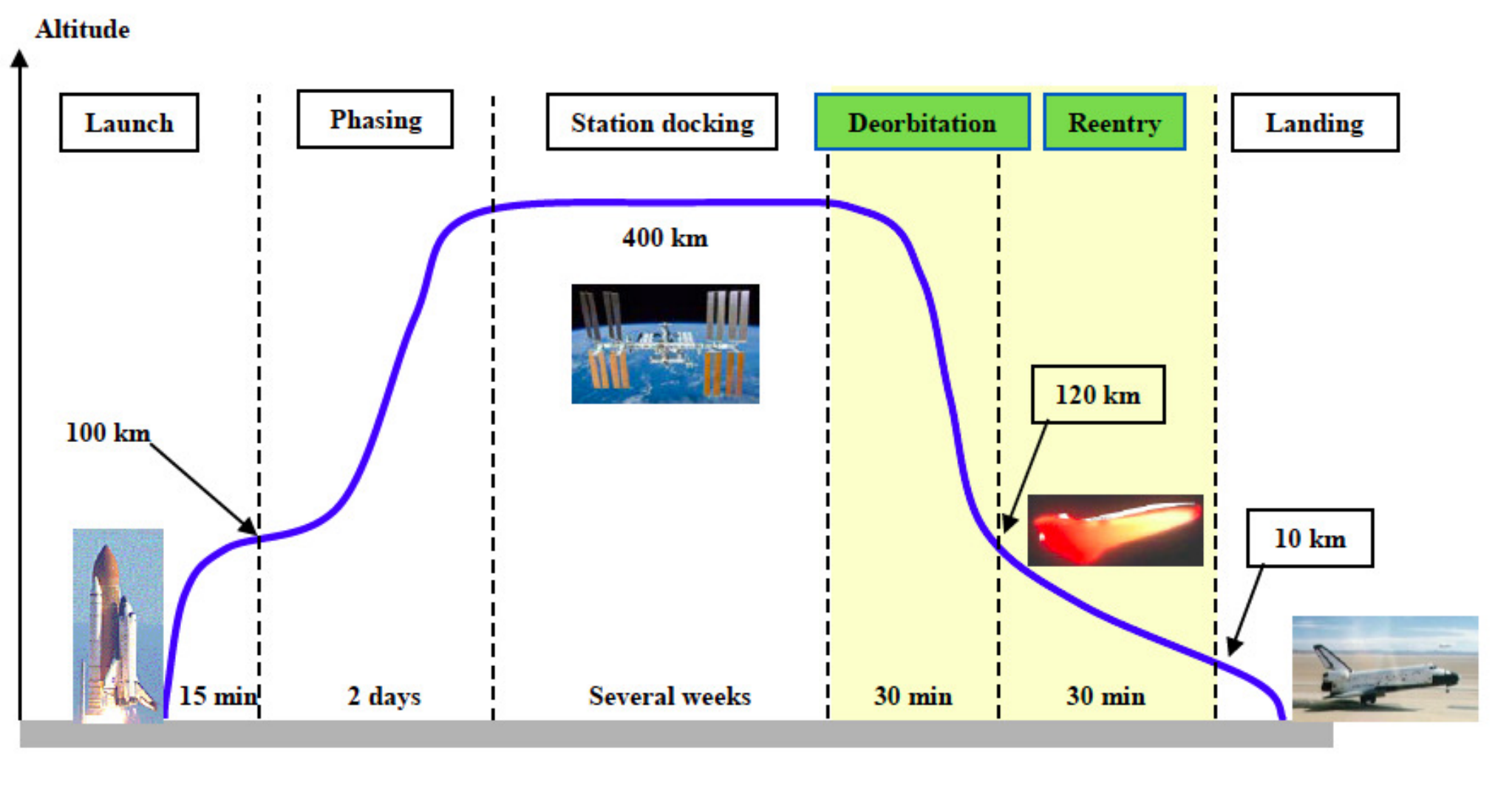}
	\caption{ISS servicing and shuttle reentry.}
	\label{ISSreentry}
\end{figure}

Reentry missions aim at retrieving either experiment results or space crews. The trajectory is split into a coast arc targeting accurate conditions at the atmospheric entry interface and a gliding atmospheric leg of about half an hour until the landing. The most stringent constraint comes from the convection flux that grows quickly when entering the dense atmosphere layers at hypersonic speeds. A near-horizontal flight is mandatory to achieve a progressive braking at limited thermal flux and load factor levels. The aerodynamic forces are controlled through the vehicle attitude. The angle of attack modulates the force magnitude and the loads applied to the vehicle. The bank angle orientates the lift left or right to follow an adequate descent rate and achieve the required downrange and cross-range until the targeted landing site. The landing may occur vertically in the sea or on the ground, or horizontally on a runway. Depending on the landing options the final braking is achieved by thrusting engines or by parachutes. If necessary the touchdown may also be damped by airbags or legs, for example for delivering scientific payloads on the Mars surface.
The reentry is always the final part of a space mission. The example of the Space Shuttle servicing the International Space Station is pictured on Figure \ref{ISSreentry}.

\begin{figure}[H]
\centering
	\includegraphics[scale=0.7]{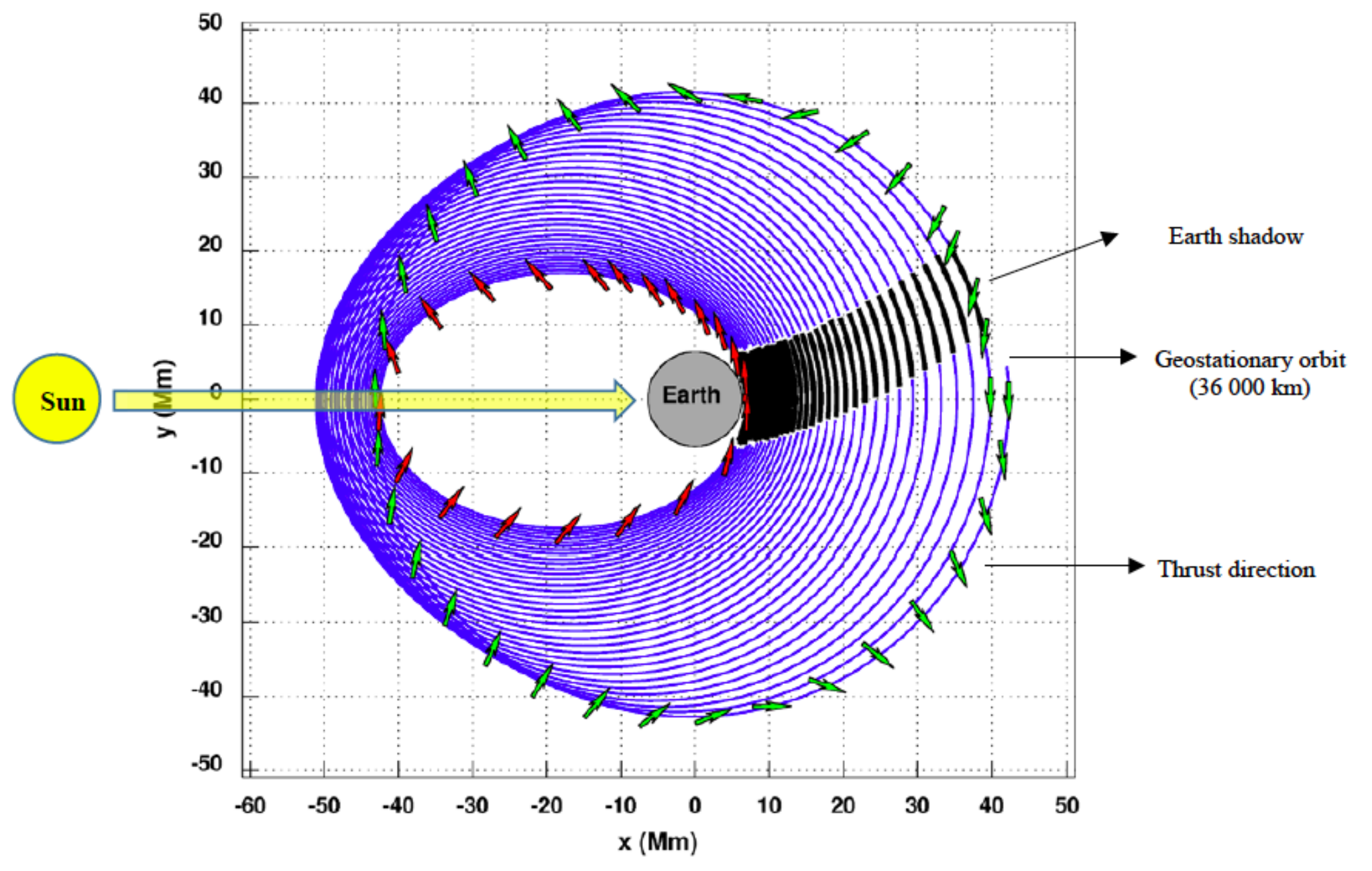}
	\caption{Low thrust orbit raising toward the geostationary orbit.}
	\label{Lowthrust}
\end{figure}

Orbital missions deal with orbit changes around the Earth and also with interplanetary travels. A major difference with ascent and reentry trajectories is the much larger duration, which ranges from days to months or even years to reach the farthest planets of the solar system. The motion is essentially due to the gravity field of the nearest body and possibly of a second one. The vehicle operational life is limited by its onboard propellant so that all propelled maneuvers must be achieved as economically as possible. Depending on the engine thrust level the maneuvers are modeled either as impulsive velocity changes (impulsive modelling) or as short duration boosts (high thrust modelling) or as long duration boosts (low thrust modelling). Low thrust engines are particularly attractive due to their high specific impulse, but they require a high electrical power that cannot be delivered by onboard batteries. The energy is provided by large solar panels and the engine must be cut-off when the vehicle enters the Earth shadow. Low thrust orbit raising of telecommunication satellites toward the geostationary orbit at 36000 km lead thus to quite complex optimal control problems as pictured on Figure \ref{Lowthrust}.

\begin{figure}[H]
\centering
	\includegraphics[scale=0.7]{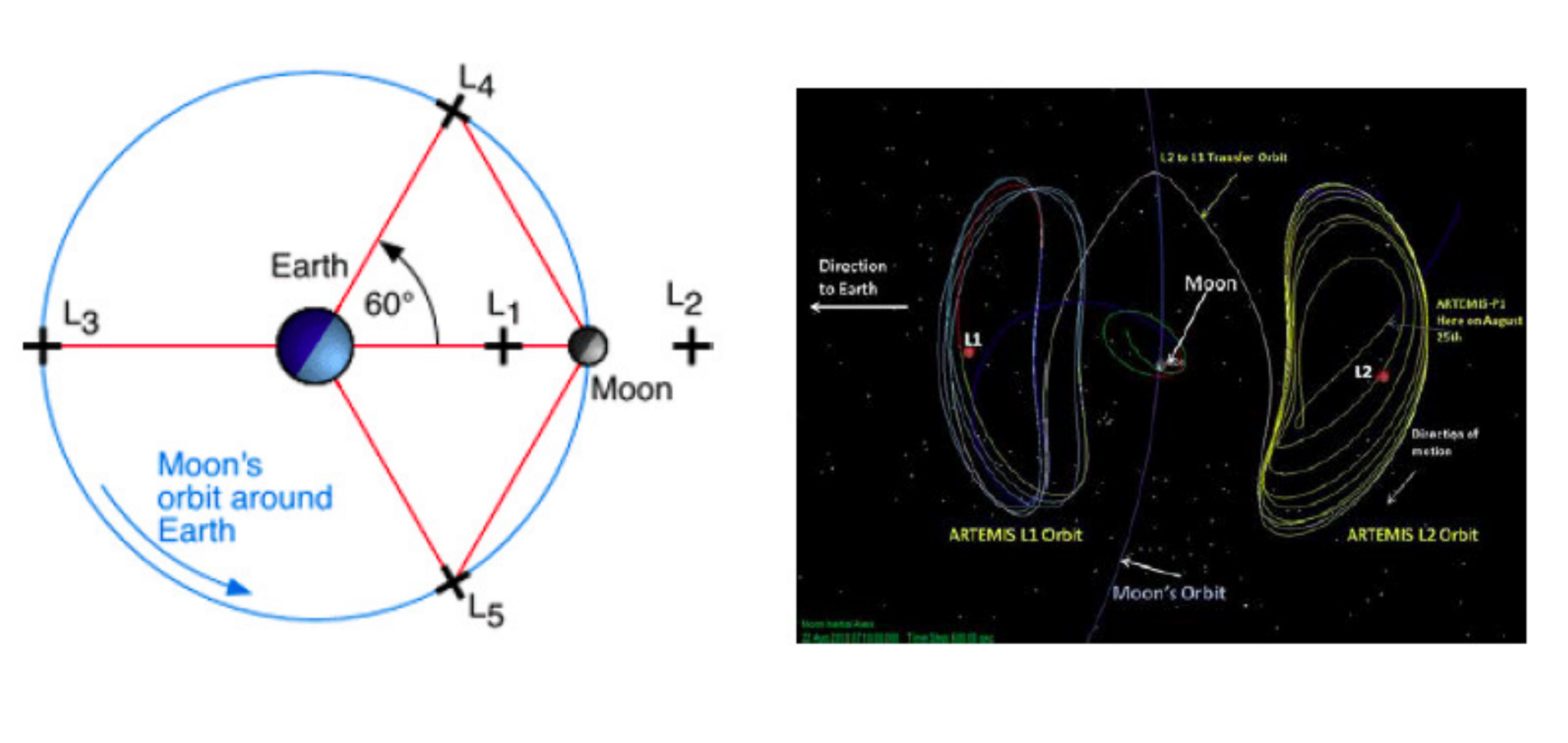}
	\caption{Earth-Moon system Lagrange points and orbits. (Sources : lagrangianpoints.com / space.stackexchange.com)}
	\label{EMorbits}
\end{figure}

Other orbital transfer problems are the removal of space debris or the rendezvous for orbit servicing. Interplanetary missions raise other difficulties due to the gravity of several attracting bodies. 
For missions towards the Lagrange points (see Figure \ref{EMorbits}) the detailed analysis of manifolds in the three body problem can provide very inexpensive transfer solutions. 

For farther solar system travels successive fly-bys around selected planets allow ``free'' velocity gains. The resulting combinatorial problem with optional intermediate deep space maneuvers is challenging.

\begin{figure}[H]
\centering
	\includegraphics[scale=0.72]{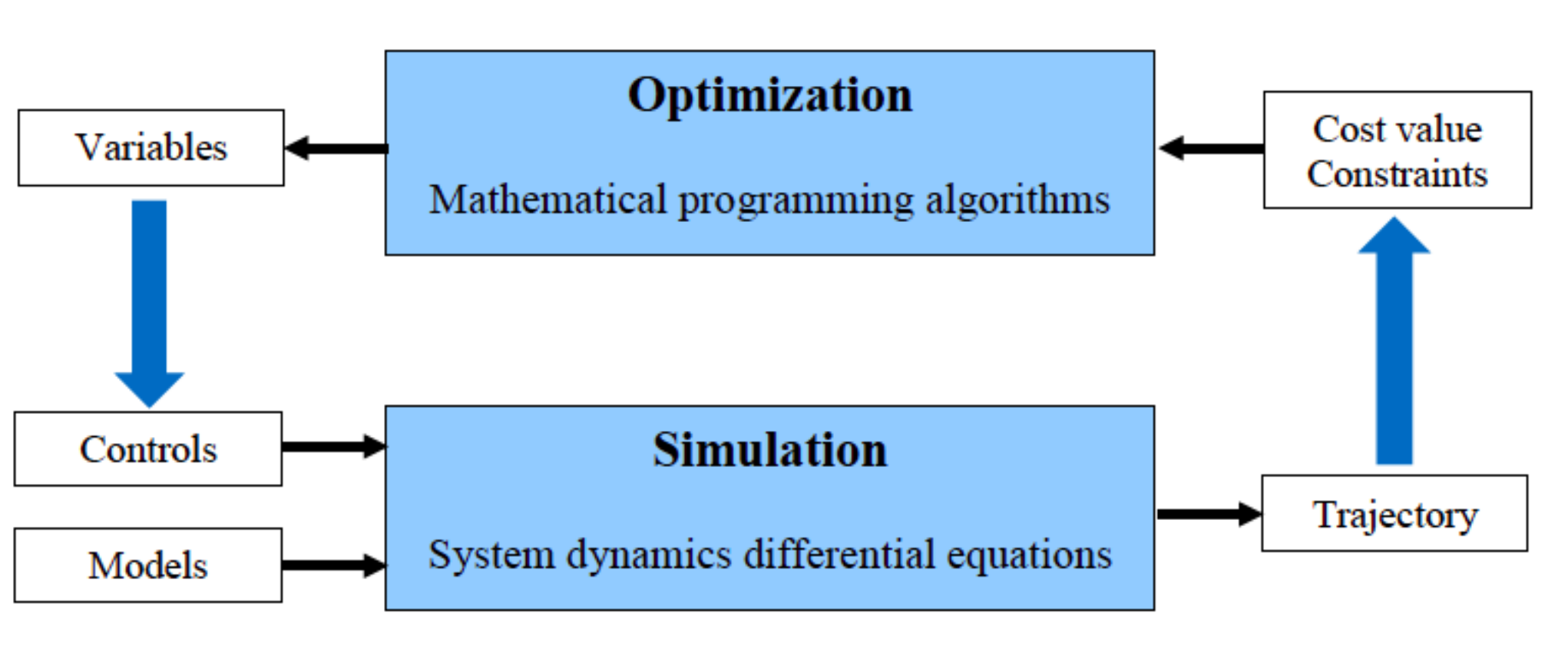}
	\caption{Simulation and optimization tasks.}
	\label{simuoptitasks}
\end{figure}

The above non exhaustive list gives a preview of various space transportation problems. In all cases the mission analysis comprises a simulation task and an optimization task (see Figure \ref{simuoptitasks}). Various formulations and methods are possible regarding these two tasks. Selecting an adequate approach is essential in order to build a satisfying numerical solution process.

The simulation task consists in integrating the dynamics differential equations derived from mechanics laws. The vehicle is generally modeled as a solid body. The motion combines the translation of the center of gravity defining the trajectory and the body rotation around its center of gravity defining the attitude. The main forces and torques originate from the gravity field (always present), from the propulsion system (when switched on) and possibly from the aerodynamics shape when the vehicle evolves in an atmosphere. In many cases a gravity model including the first zonal term due to the Earth flattening is sufficiently accurate at the mission analysis stage. The aerodynamics is generally modeled by the drag and lift components tabulated versus the Mach number and the angle of attack. The atmosphere parameters (density, pressure, temperature) can be represented by an exponential model or tabulated with respect to the altitude. A higher accuracy may be required on some specific occasions, for example to forecast the possible fall-out of dangerous space debris, to assess correctly low thrust orbital transfers or complex interplanetary space missions. In such cases the dynamical model must be enhanced to account for effects of smaller magnitudes. These enhancements include higher order terms of the gravitational field, accurate atmosphere models depending on the season and the geographic position, extended aerodynamic databases, third body attraction, etc, and also other effects such as the solar wind pressure or the magnetic induced forces.

\begin{figure}[H]
\centering
	\includegraphics[scale=0.74]{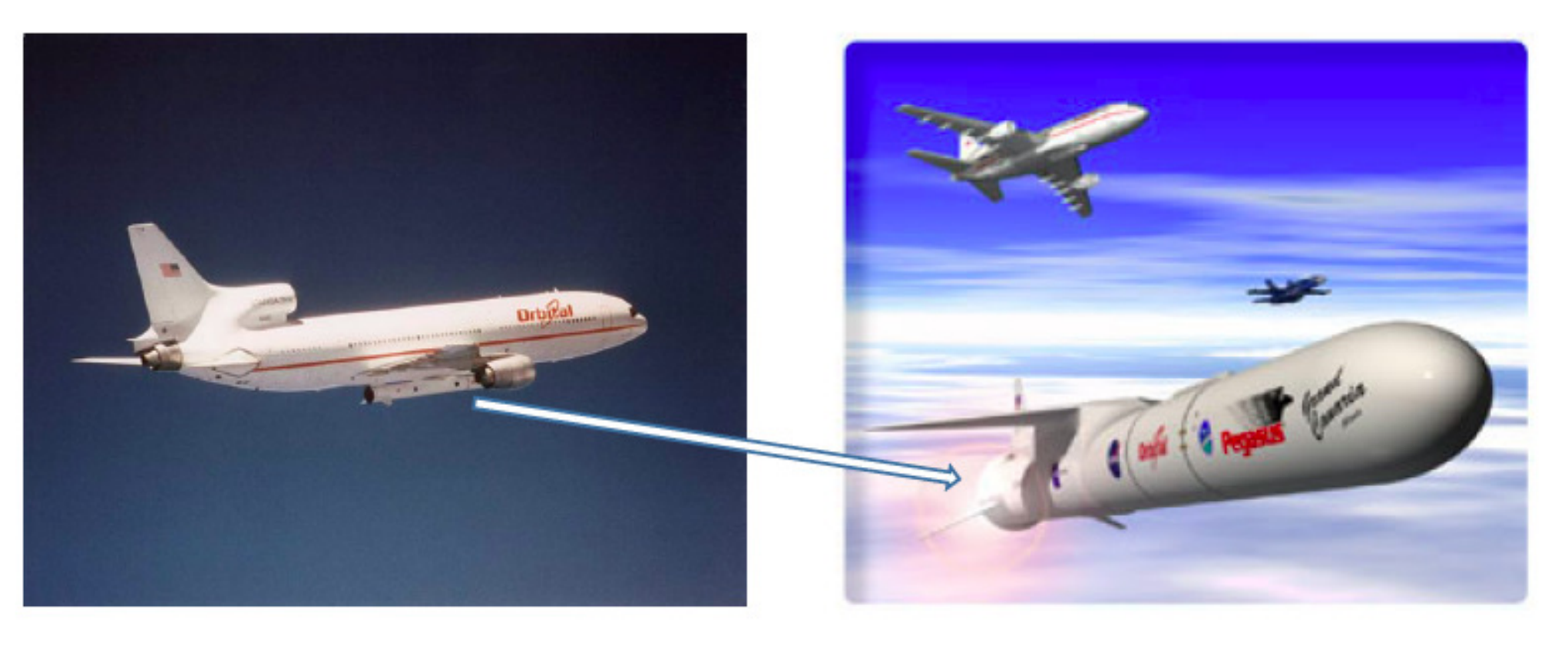}
	\caption{Pegasus airborne launcher before and after release. (Sources : air-and-space.com / spacewar.com)}
	\label{airbornelaunchers}
\end{figure}

Complex dynamical models yield more representative results at the expense of larger computation times. In view of trajectory optimization purposes the simulation models have to make compromises between accuracy and speed. A usual simplification consists in assuming that the translation and the rotation motions are independent. With this assumption the trajectory problem (also called the guidance problem) and the attitude problem (also called the control problem) can be addressed separately. This uncoupling of the guidance and the control problem is valid either when the torque commands have a negligible effect on the CoG motion or when the control time scale is much shorter than the guidance time scale. Most space vehicles fall into one of these two categories. The main exceptions are atmospheric maneuvering vehicles such as cruise or anti-ballistic missiles and airborne launchers (see Figure \ref{airbornelaunchers}).

Such vehicles have to perform large reorientation maneuvers requiring significant durations. These maneuvers have a sensible influence of the CoG motion and they must be accounted for a realistic trajectory optimization.

Another way to speed up the simulation consists in splitting the trajectory into successive sequences using different dynamical models and propagation methods. Ascent or reentry trajectories are thus split into propelled, coast and gliding legs, while interplanetary missions are modeled by patched conics. Each leg is computed with its specific coordinate system and numerical integrator. Usual state vector choices are Cartesian coordinates for ascent trajectories, orbital parameters for orbital transfers, spherical coordinate for reentry trajectories. The reference frame is usually Galilean for most applications excepted for the reentry assessment. In this case an Earth rotating frame is more suited to formulate the landing constraints. The propagation of the dynamics equations may be achieved either by semi-analytical or numerical integrators. Semi-analytical integrators require significant mathematical efforts prior to the implementation and they are specialized to a given modelling. For example averaging techniques are particularly useful for long time-scale problems, such as low thrust transfers or space debris evolution, in order to provide high speed simulations with good differentiability features. On the other hand numerical integrators can be applied very directly to any dynamical problem. An adequate compromise has then to be found between the time-step as large as possible and the error tolerance depending on the desired accuracy.

\begin{figure}[H]
\centering
	\includegraphics[scale=0.73]{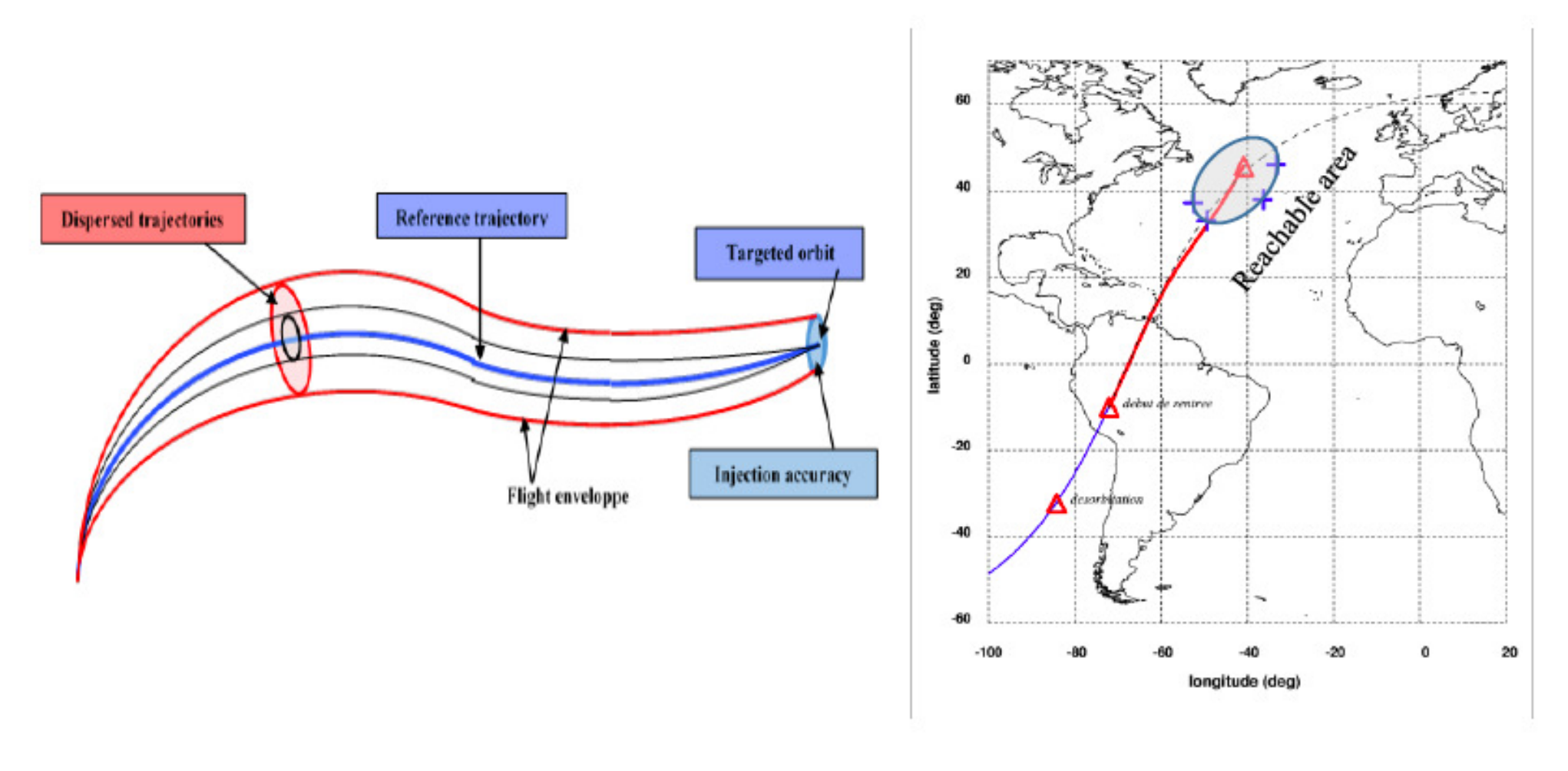}
	\caption{Dispersed flight trajectories (left) and reachable landing area (right)}
	\label{Dispertraj}
\end{figure}

The dynamics models consider first nominal features of the vehicle and of its environment in order to build a reference mission profile. Since the real flight conditions are never perfectly known, the analysis must also be extended with model uncertainties, first to assess sufficient margins when designing a future vehicle, then to ensure the required success probability and the flight safety when preparing an operational flight. The desired robustness may be obtained by additional propellant reserves for a launcher, or by reachable landing areas for a reentry glider (see Figure \ref{Dispertraj}).

The optimization task consists in finding the vehicle commands and optionally some design parameters in order to fulfill the mission constraints at the best cost. In most cases, the optimization deals only with the path followed by one vehicle. In more complicated cases, the optimization must account for moving targets or other vehicles that may be jettisoned parts of the main vehicle. Examples or such missions are debris removal, orbital rendezvous, interplanetary travel or reusable launchers with recovery of the stages after their separation.

\medskip

\begin{figure}[h]
\centering
	\includegraphics[scale=0.73]{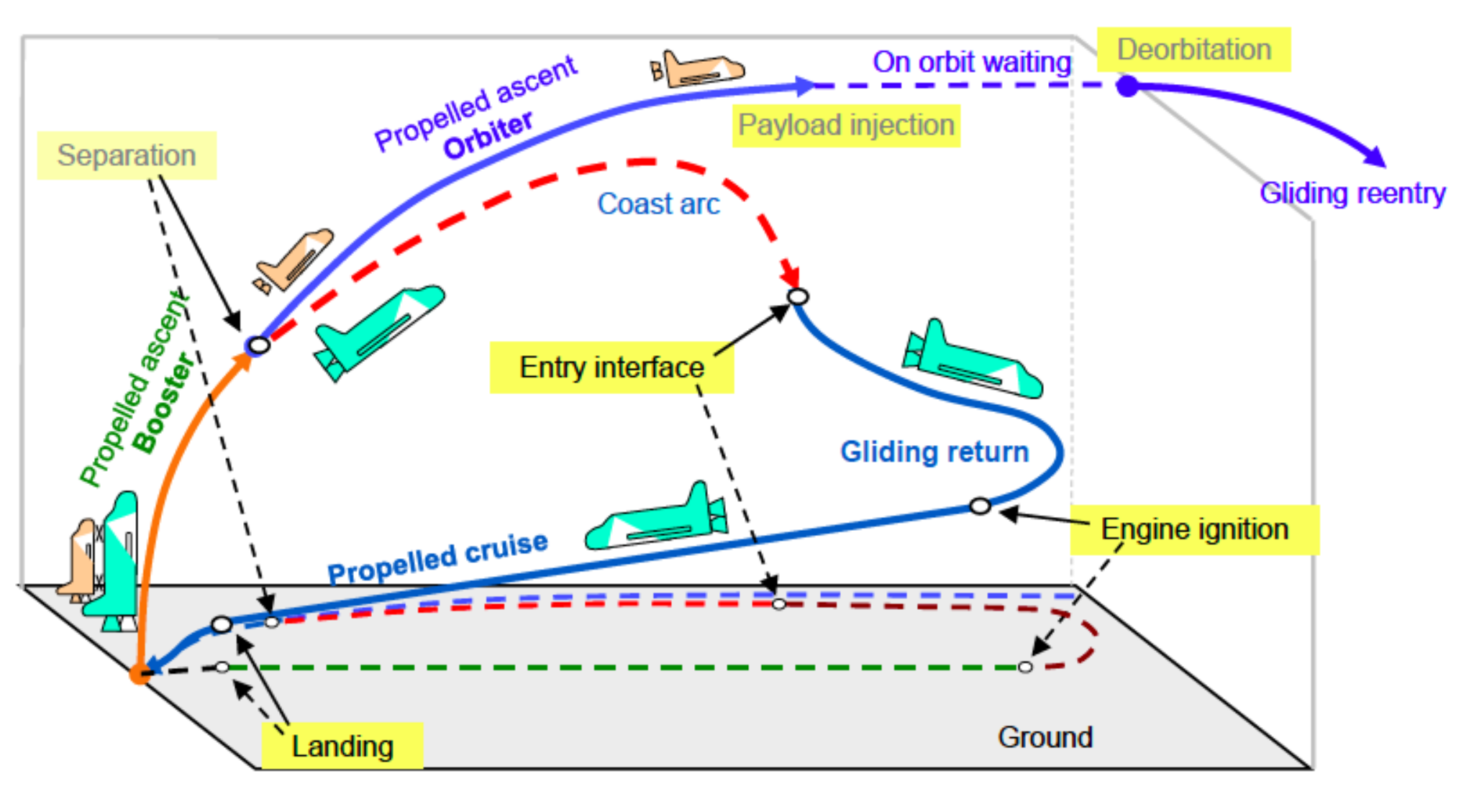}
	\caption{Reusable launch vehicle}
	\label{Rlv}
\end{figure}

A typical reusable launcher mission is pictured on Figure \ref{Rlv}. The goal is to reach the targeted orbit with the upper stage carrying the payload, while the lower and the upper stage must be recovered safely for the next launches. This problem necessitates a multi-branch modelling and a coordinated optimization method.

For preliminary design studies, the vehicle configuration is not defined. The optimization has to deal simultaneously with the vehicle design and the trajectory control.  Depending on the problem formulation the optimization variables may thus be functions, reals or integers.

In almost all cases an optimal control problem must be solved to find the vehicle command law along the trajectory. The command aims at changing the magnitude and the direction of the forces applied, namely the thrust and the aerodynamic force. The attitude time scale is often much shorter than the trajectory time scale so that the attitude control can be considered as nearly perfect, i.e., instantaneous or with a short response time. The rotation dynamics is thus not simulated and the command is directly the vehicle attitude. If the rotation and the translation motions are coupled, the 6 degrees of freedom must be simulated. The command are then the nozzle or the flap deflections depending on the vehicle control devices. The choice of the attitude angles depends on the mission dynamics. For a propelled launcher, the motion is controlled by the thrust force which is nearly aligned with the roll axis. This axis is orientated by inertial pitch and yaw angles. For a gliding reentry vehicle, the motion is controlled by the drag and lift forces. The angle of attack modulates the force magnitude while the bank angle only acts on the lift direction. For orbital maneuvering vehicles, the dynamics is generally formulated using the orbital parameters evolution, e.g., by Gauss equations, so that attitude angles in the local orbital frame are best suited.

\medskip

If the trajectory comprises multiple branches or successive flight sequences with dynamics changes and interior point constraints, discontinuities may occur in the optimal command law. This occurs typically at stage separations and engine ignitions or shutdowns. The commutation dates between the flight sequences themselves may be part of the optimized variables, as well as other finite dimension parameters, leading to a hybrid optimal control problem. A further complexity occurs with path constraints relating either to the vehicle design (e.g., dynamic pressure or thermal flux levels), or to the operations (e.g., tracking, safety, lightening). These constraints may be active along some parts of the trajectory, and the junction between constrained and unconstrained arcs may raise theoretical and numerical issues.

\medskip

The numerical procedures for optimal control problems are usually classified between direct and indirect methods. Direct methods discretize the optimal control problem in order to rewrite it as a nonlinear large scale optimization problem. The process is straightforward and it can be applied in a systematic manner to any optimal control problem. New variables or constraints may be added easily. But achieving an accurate solution requires a careful discretization and the convergence may be difficult due to the large number of variables. On the other hand indirect methods are based on the Pontryagin Maximum Principle which gives a set of necessary conditions for a local minimum. The problem is reduced to a nonlinear system that is generally solved by a shooting method using a Newton-like algorithm. The convergence is fast and accurate, but the method requires both an adequate starting point and a high integration accuracy. The sensitivity to the initial guess can be lowered by multiple shooting which breaks the trajectory into several legs linked by interface constraints, at the expense of a larger nonlinear system. The indirect method requires also prior theoretical work for problems with singular solutions or with state constraints. Handling these constraints by penalty method can avoid numerical issues, but yields less optimal solutions.

\medskip

\begin{figure}[H]
\centering
	\includegraphics[scale=0.74]{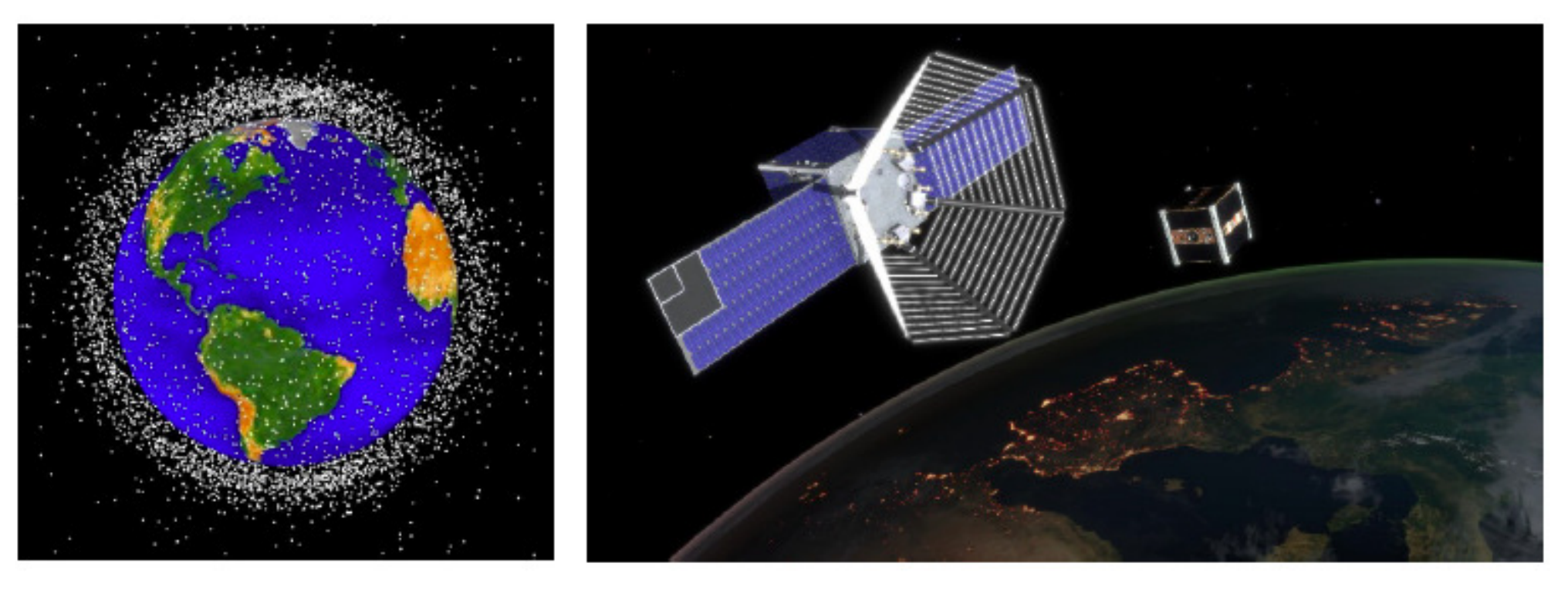}
	\caption{Space debris. (Sources : nasa.gov / leonarddavid.com : Credit: 2015 EPFL/Jamani Caillet)}
	\label{debris}
\end{figure}

In some cases the mission analysis may address discrete variables. Examples of such problems are the removal of space debris by a cleaner vehicle or interplanetary travels with multiple fly-bys. For a debris cleaning mission (see Figure \ref{debris}) the successive targets are moving independently of the vehicle, and the propellant required to go from one target to another depends on the rendezvous dates. The optimization aims at selecting the targets and the visiting order in order to minimize the required propellant. The path between two given targets is obtained by solving a time-dependent optimal control problem. The overall problem is thus a combinatorial variant of the well-known Traveling Salesman Problem, with successive embedded optimal control problems.

For an interplanetary mission successive fly-bys around planets are necessary to increase progressively the velocity in the solar system and reach far destinations. Additional propelled maneuvers are necessary either at the fly-by or in the deep space in order to achieve the desired path. An impulsive velocity modelling is considered for these maneuvers in a first stage. If a low thrust engine is used, the maneuver assessment must be refined by solving an embedded optimal control problem. The optimization problem mixes discrete variables (selected planets, number of revolutions between two successive fly-bys, number of propelled maneuvers) and continuous variables (fly-bys dates, maneuver dates, magnitudes and orientations).

\medskip

\begin{figure}[H]
\centering
	\includegraphics[scale=0.73]{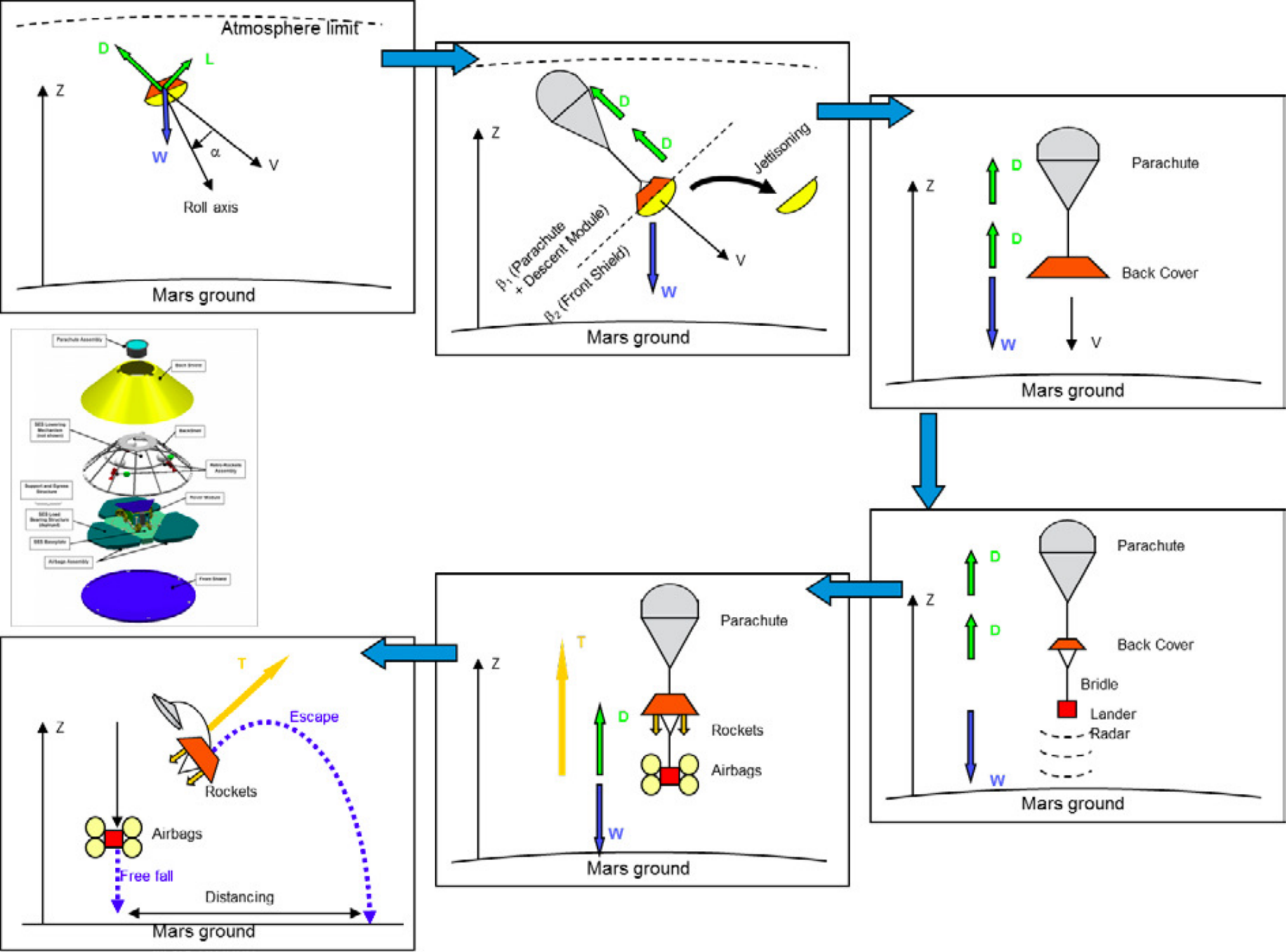}
	\caption{Entry, descent and landing system design}
	\label{mission}
\end{figure}

In preliminary design studies, the optimization problem addresses simultaneously the vehicle configuration and its command along the trajectory. The goal is usually to find the minimal gross weight vehicle able to achieve the specified mission. The configuration parameters are either continuous or discrete variables. For a propelled vehicle the main design parameters are the number of stages, the number of engines, the thrust level, the propellant type and the propellant masses. For a reentry vehicle the design is driven by the aerodynamic shape, the surface and by the auxiliary braking sub-systems if any. The gross mass minimization is essential for the feasibility of interplanetary missions. An example is given by a Mars lander composed of a heat shield, one or several parachutes, braking engines, airbags and legs. The sub-system designs drive the acceptable load levels and thus the state constraints applied to the entry trajectory. The successive sequence of the descent trajectory are depicted on Figure \ref{mission}. Large uncertainties have also to be accounted regarding the Mars environment in order to define a robust vehicle configuration.

Multidisciplinary optimization deals with such problems involving both the vehicle design and the mission scenario. The overall problem is too complex to be address directly, and a specific optimization procedure must be devised for each new case. A bi-level approach consists in separating the design and the trajectory optimization. The design problem is generally non differentiable or may present many local minima. It can be addressed in some cases by mixed optimization methods like branch and bound, or more generally by meta-heuristics like simulated annealing, genetic algorithms, particle swarm, etc. None is intrinsically better than another and a specific analysis is needed to formulate the optimization problem in a way suited to the selected method. These algorithms are based partly on a random exploration of the variable space. In order to be successful the exploration strategy has to be customized to the problem specificities. Thousands or millions of trials may be necessary to yield a candidate configuration, based on very simplified performance assessment (e.g., analytical solutions, impulsive velocities, response surface models etc.). The trajectory problem is then solved for this candidate solution in order to assess the real performance, and if necessary iterate on the configuration optimization with a corrected the performance model. Meta-heuristics may also be combined with multi-objective optimization approaches since several criteria have to be balanced at the design stage of a new space vehicle. The goal is to build a family of launchers using a common architecture of propelled stages with variants depending the targeted orbit and payload. By this way the development and manufacturing costs are minimized while the launcher configuration and the launch cost can be customized for each flight.

\section{Geometric Optimal Control}
\label{sec_GOCR}
Geometric optimal control (see, e.g., \cite{Agrachev,Schattler,Trelatsurvey}) combines classical optimal control and geometric methods in system theory, with the goal of achieving optimal synthesis results.
 More precisely, by combining the knowledge inferred from the Pontryagin Maximum Principle (PMP) with geometric considerations,
 such as the use of Lie brackets and Lie algebras,
 of differential geometry on manifolds,
 and of symplectic geometry and Hamiltonian systems,
the aim is to describe in a precise way the structure of optimal trajectories.
We refer the reader to \cite{Trelatsurvey,Schattler2012} for a list of references on geometric tools used in geometric optimal control.
The foundations of geometric control can be dated back to the Chow's theorem and to \cite{Brunovsky1978,Brunovsky1980}, where Brunovsky found that it was possible to derive regular synthesis results by using geometric considerations for a large class of control systems.
Apart from the main goal of achieving a complete \emph{optimal synthesis}, geometric control aims also at deriving higher-order optimality conditions in order to better characterize the set of candidate optimal trajectories.

\medskip

In this section, we formulate the optimal control problem on differentiable manifolds and recall some tools and results from geometric optimal control.
More precisely, the Lie derivative is used to define the order of the state constraints, the Lie and Poisson brackets are used to analyze the singular extremals and to derive higher order optimality conditions, and the optimality conditions (order one, two and higher) are used to analyze the chattering extremals (see Section \ref{chatphenom} for the chattering phenomenon). These results will be applied in Section \ref{fullexample} on a coupled attitude and trajectory optimization problem.

\subsection{Optimal Control Problem}
Let $M$ be a smooth manifold of dimension $n$, let $N$ be a smooth manifold of dimension $m$, let $M_0$ and $M_1$ be two subsets of $M$, and let $U$ be a subset of $N$.
We consider the general nonlinear optimal control problem ($\P_0$), of minimizing the cost functional
$$
C(t_f,u) = \int_0^{t_f} f^0(x(t),u(t))dt + g(t_f,x(t_f)),
$$
over all possible trajectories solutions of the control system
\begin{equation} \label{sys_general}
\dot{x}(t) = f(x(t),u(t)),
\end{equation}
and satisfying the terminal conditions
\begin{equation} \label{tercond_general}
x(0) \in M_0, \quad x(t_f) \in M_1 ,
\end{equation}
where the mappings $f: M \times N \rightarrow TM$, $f^0: M \times N \rightarrow \R$, and $g: \R \times M \rightarrow \R$ are smooth,
and where the controls are bounded and measurable functions defined on $[0,t_f(u)]$ of $\R^+$, taking values in $U$.
The final time $t_f$ may be fixed or not.
We denote $\U$ the set of \emph{admissible controls} such that the corresponding trajectories steer the system from an initial point of $M_0$ to a final point in $M_1$.

For each $x(0) \in M_0$ and $u \in \U$, we can integrate the system \eqref{sys_general} from $t=0$ to $t=t_f$, and assess the cost $C(t_f,u)$ corresponding to $x(t)=x(t;x_0,u(t))$ and $u(t)$ for $t=[0,t_f]$.
Solving the problem ($\P_0$) consists in finding a pair $(x(t),u(t))=(x(t;x_0,u(t)),u(t))$ minimizing the cost.
For convenience, we define the end-point mapping to describe the final point of the trajectory solution of the control system \eqref{sys_general}.

\begin{defi} \label{def_endpmap}
The \emph{end-point mapping} $E: M \times \R \times \U$ of the system is defined by 
$$E(x_0,t_f,u)=x(x_0,t_f,u),$$
where $t \mapsto x(x_0,t,u)$ is the trajectory solution of the control system \eqref{sys_general} associated to $u$ such that $x(x_0,0,u)=x_0$
\end{defi}

Assuming moreover that $\U$ is endowed with the standard $L^\infty$ topology, then the end-point mapping is $C^1$ on $\U$,
and in terms of the end-point mapping, the optimal control problem under consideration can be written as the infinite-dimensional minimization problem
$$
\min \left\{ C(t_f,u)\, | \, x_0 \in M_0, E(x_0,t_f,u)\in M_1, u \in L^\infty([0,t_f];U)\right\}.
$$
This formulation of the problem will be used when we introduce the Lagrange multipliers rule in Section \ref{subsubsec_langrange} in a simpler case when $M_0 = \left\{ x_0 \right\}$ and $M_1 = \left\{x_1\right\}$ and $U=\R^m$.

If the optimal control problem has a solution, we say that the corresponding control and trajectory are minimizing or optimal.
We refer to \cite{Cesari,Trelatbook} for existence results in optimal control.

\bigskip

Next, we introduce briefly the concept of Lie derivative, and of Lie and Poisson brackets (used in Section \ref{hooc} for higher order optimality conditions). These concepts will be applied in Section \ref{fullexample} to analyze the pull-up maneuver problem.

\subsection{Lie Derivative, Lie Bracket, and Poisson Bracket} \label{Lie}
Let $\Omega$ be an open and connected subset in $M$, and denote the space of all infinitely continuously differentiable functions on $\Omega$ by $C^\infty(\Omega)$.
Let $X \in C^\infty(\Omega)$ be a vector field. $X$ can be seen as defining a first-order differential operator from the space $C^\infty(\Omega)$ into $C^\infty(\Omega)$ by taking at every point $q \in \Omega$ the directional derivative of a function $\varphi \in C^\infty(\Omega)$ in the direction of the vector field $X(q)$, i.e.,
$$
X:C^\infty(\Omega) \to C^\infty(\Omega),\quad
\varphi \mapsto X\varphi,
$$
defined by
$$
(X. \varphi)(q)=\nabla \varphi(q) \cdot X(q).
$$

We call $(X. \varphi)(q)$ the \emph{Lie derivative} of the function $\varphi$ along the vector field $X$, and generally one denote the operator by $L_X$, i.e., 
$$L_X(\varphi)(q) = (X. \varphi)(q).$$
In general, the order of the state constraints in optimal control problems is defined through Lie derivatives as we will show on the example in Section \ref{example2}.

\begin{defi}
The \emph{Lie bracket} of two vector fields $X$ and $Y$ defined on a domain $\Omega$ is the operator defined by the commutator
$$
[X,Y]=X \circ Y - Y\circ X = XY-YX.
$$
\end{defi}
The Lie bracket actually defines a first-order differential operator. For any function $\varphi$ we have 
\begin{equation*}
\begin{split}
[X,Y](\varphi)
& = X(Y.\varphi)-Y(X.\varphi) \\
& = X(\nabla \varphi Y)-Y(\nabla \varphi X)\\
& = \nabla(\nabla \varphi Y)X - \nabla(\nabla \varphi X)Y \\
& = \nabla(\nabla \varphi) (Y,X)- \nabla(\nabla \varphi)(X,Y) + \nabla \varphi (DY \cdot X-DX \cdot Y) \\
& =  \nabla \varphi (DY \cdot X- DX \cdot Y),
\end{split}
\end{equation*}
where $\nabla(\nabla \varphi)(X,Y)$ denotes the action of the Hessian matrix of the function $\varphi$ on the vector fields $X$ and $Y$, and $DX$ and $DY$ denote the matrices of the partial derivatives of the vector fields $X$ and $Y$.
Therefore, if $X:\Omega \to M, z \mapsto X(z)$, and $Y:\Omega \to M, z \mapsto Y(z)$, are coordinates for these vector fields, then
$$
[X,Y](z)= DY(z)\cdot X(z)-DX(z) \cdot Y(z).
$$

\begin{lem}
Let $X$, $Y$, and $Z$ be three $C^\infty$ vector fields defined on $\Omega$, and let $\alpha$, $\beta$ be smooth functions on $\Omega$. The Lie bracket has the following properties:
\begin{itemize}
\item $[\cdot,\cdot]$ is a bilinear operator;
\item $[X,Y] = -[Y,X]$;
\item $[X+Y,Z]=[X,Z]+[Y,Z]$;
\item $[X,[Y,Z]] + [Y,[Z,X]]+[Z,[X,Y]] = 0$ (Jacobi identity);
\item $[\alpha X, \beta Y] = \alpha \beta [X,Y] +\alpha(L_X\beta)Y - \beta (L_Y\alpha)X$.
\end{itemize}
\end{lem}

These properties show that the vector fields (as differential operators) form a Lie algebra.
A \emph{Lie algebra} over $\R$ is a real vector space $\G$ together with a bilinear operator $[\cdot,\cdot]: \G \times \G \to \G$ such that for all $X,Y,Z \in \G$ we have $[X,Y] = -[Y,X]$ and $[X+Y,Z]=[X,Z]+[Y,Z]$.

\medskip
Going back to the problem ($\P_0$), we assume that $f(x,u) = f_0(x)+uf_1(x)$,
$f^0(x,u)=1$, and $g(t,x)=0$, and we define a $C^1$ function by
$$
h(x,p) = \langle p,Z(x) \rangle,
$$
where $p$ is the adoint vector and $Z$ is a vector field. The function $h$ is the Hamiltonian lift of the vector field $Z$. Accordingly, and with a slight abuse of notation, we denote by $h(t) = h(x(t),p(t))$ the value at time $t$ of $h$ along a given extremal.
The derivative of this function is
\begin{equation} \label{deriLiebra}
\begin{split}
\dot{h}(t) & = \langle \dot{p},Z(x) \rangle + \langle p,DZ(x) \dot{x} \rangle \\
& = - \langle p(Df_0(x)+uDf_1(x)),Z(x) \rangle +  \langle p,DZ(x) (f_0(x)+uf_1(x)) \rangle \\
& =  \langle p, DZ(x)f_0(x)-Df_0(x)Z(x) \rangle + u \langle p, DZ(x)f_1(x)-Df_1(x)Z(x) \rangle \\
& = \langle p, [f_0,Z](x) \rangle + u \langle p,[f_1,Z](x) \rangle.
\end{split}
\end{equation}

Let us recall also the concept of the \emph{Poisson bracket}. The Poisson bracket is related to the Hamiltonians. In the canonical coordinates $z=(x,p)$, given two $C^1$ functions $\alpha_1(x,p)$ and $\alpha_2(x,p)$, the Poisson bracket takes the form
$$
\left\{ \alpha_1,\alpha_2\right\} (x,p) = \frac{\partial \alpha_2}{\partial x} \frac{\partial \alpha_1}{\partial p} - \frac{\partial \alpha_1}{\partial x} \frac{\partial \alpha_2}{\partial p}.
$$

According to \eqref{deriLiebra}, taking
$$
\alpha_1(x(t),p(t))=H(x(t),p(t)), \quad \alpha_2(x(t),p(t))=h(x(t),p(t)),
$$
we have
\begin{equation*} 
\dot{h}(t)  = \left\{H,h\right\}(x(t),p(t)) = \left\{h_0,h\right\}(x(t),p(t)) + u \left\{h_1,h\right\} (x(t),p(t)),
\end{equation*}
where $h_0(t)=\langle p(t),f_0(x(t)) \rangle$ and $h_1(t)=\langle p(t),f_1(x(t)) \rangle$.

For convenience, we adopt the usual notations
$$
\mathrm{ad}\, f_0.f_1 = [f_0,f_1],\, \textrm{resp.} \, \mathrm{ad}\, h_0.h_1 = \left\{h_0,h_1\right\},
$$
and
$$
\mathrm{ad}^if_0.f_1 = [f_0,\mathrm{ad}^{i-1}f_0.f_1],\, \textrm{resp.} \, \mathrm{ad}^ih_0.h_1 = \left\{h_0,\mathrm{ad}^{i-1}h_0.h_1\right\}.
$$

\medskip

We will see in Section \ref{opticondisec} (and also in Section \ref{fullexample}) that the Lie brackets and the Poisson brackets are very useful for deriving higher order optimality conditions in simpler form and for calculating the singular controls.

\subsection{Optimality Conditions}
\label{opticondisec}
This section gives an overview of necessary optimality conditions.

For the first-order optimality conditions, we recall the Lagrange multipliers method for the optimal control problem without control constraints. Such constraints can be accounted in the Lagrangian with additional Lagrange multipliers \cite{Bryson}. This method leads to weaker results than the Pontryagin Maximum Principle which considers needle-like variations accounting directly for the control constraints.

In some cases, the first-order conditions do not provide adequate information of the optimal control, and the higher order optimality conditions are needed.
Therefore we recall the second and higher order necessary optimality conditions that must be met by any trajectory associated to an optimal control $u$.  These conditions are especially useful to analyze the singular solutions because the first-order optimality conditions do not provide any information in such cases.

\subsubsection{First-Order Optimality Conditions}

\paragraph{Lagrange multipliers rule.}
\label{subsubsec_langrange}
We consider the simplified problem ($\P_0$) with $M=\R^n$, $M_0=\left\{x_0\right\}$, $M_1=\left\{ x_1\right\}$, and $U=\R^m$.
According to the well known Lagrange multipliers rule (and assuming the $C^1$ regularity of the problem), if $x \in M$ is optimal then there exists a nontrivial couple $(\psi,\psi^0) \in \R^n \times \R $ such that
\begin{equation} \label{lagcond}
\psi.dE_{x_0,t_f}(u)+\psi^0 dC_{t_f}(u)=0,
\end{equation}
where $dE(\cdot)$ and $dC(\cdot)$ denote the Fr\'echet derivative of $E(\cdot)$ and $C(\cdot)$, respectively.
Defining the Lagrangian by 
$$L_{t_f} = \psi E_{x_0,t_f}(u)+\psi^0 dC_{t_f}(u),$$
this first-order necessary condition can be written in the form
$$\frac{\partial L_{t_f}}{\partial u} (u,\psi,\psi^0)=0.$$
If we define as usual the intrinsic second-order derivative $Q_{t_f}$ of the Lagrangian as the Hessian $\frac{\partial^2 L_{t_f}}{\partial^2 u}(u,\psi,\psi^0)$ restricted to the subspace $\ker \frac{\partial L_{t_f}}{\partial u}$, a second-order necessary condition for optimality is the nonpositivity of $Q_{t_f}$ (with $\psi^0 \leq 0$), and a second-order sufficient condition for local optimality is the negative definiteness of $Q_{t_f}$.

\medskip
These results are weaker to those obtained with the PMP. The Lagrange multiplier $(\psi,\psi^0)$ is in fact related to the adjoint vector introduced in the PMP.
More precisely,  the Lagrange multiplier is unique up to a multiplicative scalar if and only if the trajectory $x(\cdot)$ admits a unique extremal lift up to a multiplicative scalar, and
the adjoint vector $(p(\cdot),p^0)$ can be constructed such that $(\psi,\psi^0)=(p(t_f),p^0)$ up to some multiplicative scalar. This relation can be observed from the proof of the PMP.
The Lagrange multiplier $\psi^0=p^0$ is associated with the instantaneous cost. The case with $p^0$ null is said abnormal, which means that there are no neighboring trajectories having the same terminal point
(see, e.g., \cite{Agrachev1988,Trelatsurvey}).

\paragraph{Pontryagin Maximum Principle.}
The Pontryagin Maximum Principle (PMP, see \cite{PONTRYAGIN}) for the problem ($\P_0$) with control constraints and without state constraints is recalled in the following statement.

\begin{thm} \label{thm_pmp}
If the trajectory $x(\cdot)$, associated to the optimal control $u$ on $[0,t_f]$, is optimal, then it is the projection of an extremal $(x(\cdot),p(\cdot),p^0,u(\cdot))$ where $p^0 \leq 0$, and $p(\cdot):[0,t_f] \mapsto T^\ast_{x(t)} M$\footnote{Given any $x\in M$, $T^\ast_{x} M$ is the cotangent space to $M$ at $x$.} is an absolutely continuous mapping (called adjoint vector) with $(p(\cdot),p^0)\neq 0$, such that almost everywhere on $[0,t_f]$,
\begin{equation} \label{extremals}
	\dot{x}(t) = \frac{\partial H}{\partial p}(x(t),p(t),p^0,u(t)),\quad
	\dot{p}(t) = -\frac{\partial H}{\partial x}(x(t),p(t),p^0,u(t)) ,
\end{equation}	
where the Hamiltonian is defined by
\begin{equation*}
H(x,p,p^0,u) = \langle p, f(x,u)\rangle + p^0 f^0(x,u),
\end{equation*}
and there holds almost everywhere on $[0,t_f]$.
\begin{equation} \label{Hamiltonian}
H(x(t),p(t),p^0,u(t)) = \max_{v \in U} H(x(t),p(t),p^0,v),
\end{equation}
If moreover, the final time $t_f$ is not fixed, then
\begin{equation} \label{condtrans_hamil}
\max_{v \in U} H(x(t),p(t),p^0,v) = -p^0 \frac{\partial g}{\partial t}(t_f,x(t_f)).
\end{equation}
If $M_0$ and $M_1$ (or just one of them) are submanifolds of $M$ locally around $x(0) \in M_0$ and $x(t_f) \in M_1$, then the adjoint vector satisfies the transversality conditions at both endpoints (or just one of them)
\begin{equation} \label{condtrans_adjoint}
	p(0) \perp T_{x(0)} M_0 ,\quad
	p(t_f) -p^0 \frac{\partial g}{\partial x}(t_f,x(t_f))\perp T_{x(t_f)} M_1 ,
\end{equation}
where $T_x M_0$ (resp., $T_x M_1$) denote the tangent space to $M_0$ (resp., $M_1$) at the point $x$.
\end{thm}

The quadruple $(x(\cdot),p(\cdot),p^0,u(\cdot))$ is called the \emph{extremal lift} of $x(\cdot)$. An extremal is said to be \emph{normal} (resp., \emph{abnormal}) if $p^0 < 0$ (resp., $p^0 = 0$).
According to the convention chosen in the PMP,  we consider $p^0 \leq 0$. If we adopt the opposite convention $p^0 \geq 0$, then we have to replace the maximization condition \eqref{condtrans_hamil} with a minimization condition.
When there are no control constraints, abnormal extremals project exactly onto singular trajectories.

\medskip

The proof of the PMP is based on needle-like variations and uses a conic implicit function theorem (see, e.g., \cite{Agrachev,HTrelat,STrelat}).
Since these needle-like variants are of order one, the optimality conditions given by the PMP are necessary conditions of the first-order.
For singular controls, higher order control variations are needed to obtain optimality conditions. A singular control is defined precisely as follows.
\begin{defi} \label{def_sin_E}
Assume that $M_0=\left\{x_0\right\}$. A control $u$ defined on $[0,t_f]$ is said to be \emph{singular} if and only if the Fr\'echet differential $\frac{\partial E}{\partial u}(x_0,t_f,u)$ is not of full rank.
The trajectory $x(\cdot)$ associated with a singular control $u$ is called singular trajectory.
\end{defi}

In practice the condition $\frac{\partial^2 H}{\partial u^2}(x(\cdot),p(\cdot),p^0,u(\cdot))=0$ (the Hessian of the Hamiltonian is degenerate) is used to characterize singular controls. An extremal $(x(\cdot),p(\cdot),p^0,u(\cdot))$ is said \emph{totally singular} if this condition is satisfied.
The is especially the case when the control is affine (see Section \ref{hooc}).

\medskip
The PMP claims that if a trajectory is optimal, then it should be found among projections of extremals joining the initial set to the final target.
Nevertheless the projection of a given extremal is not necessarily optimal.
This motivates the next section on second-order optimality conditions.

\subsubsection{Second-Order Optimal Conditions}
The literature on first and/or second-order sufficient conditions with continuous control is rich (see, e.g., \cite{Dunn1995,MilyutinO1998,MaurerO2002,MaurerP1995,Zeidan1994}), which is less the case for discontinuous controls (see, e.g., \cite{OsmolovskiiL2002}).
We recall hereafter the Legendre type conditions with Poisson brackets to show that geometric optimal control allows a simple expression of the second-order necessary and sufficient conditions (see Theorem \ref{Gohcond}).

\paragraph{Legendre type conditions.}
For the optimal control problem ($\P_0$), we have the following second-order optimality conditions (see, e.g., \cite{Agrachev,Bolza,BonnardChyba}).
\begin{quote}
{\it
\begin{itemize}
\item If a trajectory $x(\cdot)$, associated to a control $u$, is optimal on $[0,t_f]$ in $L^\infty$ topology, then the Legendre condition holds along every extremal lift $(x(\cdot),p(\cdot),p^0,u(\cdot))$ of $x(\cdot)$, that is
$$
\frac{\partial^2 H}{\partial u^2} (x(\cdot),p(\cdot),p^0,u(\cdot)).(v,v) \leq 0, \quad \forall v\in \R^m.
$$
\item If the strong Legendre condition holds along the extremal $(x(\cdot),p(\cdot),p^0,u(\cdot))$, that is, there exists $\epsilon_0 >0$ such that
$$
\frac{\partial^2 H}{\partial u^2} (x(\cdot),p(\cdot),p^0,u(\cdot)).(v,v) \leq - \epsilon_0 \|v\|^2, \quad \forall v\in \R^m,
$$
then there exists $\epsilon_1 > 0$ such that $x(\cdot)$ is locally optimal in $L^\infty$ topology on $[0,\epsilon_1]$.
If the extremal is moreover normal, i.e., $p^0 \ne 0$, then $x(\cdot)$ is locally optimal in $C^0$ topology on $[0,\epsilon_1]$.
\end{itemize}
}
\end{quote}

The $C^0$ local optimality and $L^\infty$ local optimality are respectively called strong local optimality and weak local optimality\footnote{If the final time $t_f$ is fixed, then $\bar{x}(\cdot)$ is said to be locally optimal in $L^\infty$ topology (resp. in $C^0$ topology), if it is optimal in a neighborhood of $u$ in $L^\infty$ topology (resp. in a neighborhood of $\bar{x}(\cdot)$ $C^0$ topology). \\
\indent If the final time $t_f$ is not fixed, then a trajectory $\bar{x}(\cdot)$ is said to be locally optimal in $L^\infty$ topology if, for every neighborhood $V$ of $u$ in $L^\infty([0,t_f+\epsilon],U)$, for every real number $\eta$ so that $| \eta | \leq \epsilon$, for every control $v \in V$ satisfying $E(x_0,t_f+\eta,v) = E(x_0,t_f,u)$ there holds $C(t_f+\eta,v) \geq C(t_f,u)$.
Moreover, a trajectory $\bar{x}(\cdot)$ is said to be locally optimal in $C^0$ topology if, for every neighborhood $W$ of $\bar{x}(\cdot)$ in $M$, for every real number $\eta$ so that $| \eta | \leq \epsilon$, for every trajectory $x(\cdot)$, associated to a control $v \in V$ on $[0,t_f+\eta]$, contained in $W$, and satisfying $x(0) = \bar{x}(0) = x_0$, $x(t_f+\eta) = \bar{x}(t_f)$, there holds $C(t_f+\eta,v) \geq C(t_f,u)$.}.
The Legendre condition is a necessary optimality condition, whereas the strong Legendre condition is a sufficient optimality condition.
We say that we are in the \emph{regular case} whenever the strong Legendre condition holds along the extremal.
Under the strong Legendre condition, a standard implicit function argument allows expressing, at least locally, the control $u$ as a function of $x$ and $p$.

In the totally singular case, the strong Legendre condition is not satisfied and we have the following generalized condition \cite{Agrachev,Goh}.

\begin{thm} \label{Gohcond}
\emph{(Goh and Generalized Legendre condition)}
\begin{itemize}
\item If a trajectory $x(\cdot)$, associated to a piecewise smooth control $u$, and having a totally singular extremal lift $(x(\cdot),p(\cdot),p^0,u(\cdot))$, is optimal on $[0,t_f]$ in $L^\infty$ topology, then the Goh condition holds along the extremal, that is
$$
\left\{\frac{\partial H}{\partial u_i},\frac{\partial H}{\partial u_j} \right\} = 0,
$$
where $\left\{ \cdot,\cdot \right\}$ denotes the Poisson bracket on $T^\ast M$. Moreover, the generalized Legendre condition holds along every extremal lift $(x(\cdot),p(\cdot),p^0,u(\cdot))$ of $x(\cdot)$, that is
$$
\left\{ \left\{ H, \frac{\partial H}{\partial u}.v \right\},\frac{\partial H}{\partial u}.v \right\}
+ \left\{ \frac{\partial^2 H}{\partial u^2}.(\dot{u},v),\frac{\partial H}{\partial u}.v \right\}
\leq 0, \quad \forall v\in \R^m.
$$
\item If the Goh condition holds along the extremal lift $(x(\cdot),p(\cdot),p^0,u(\cdot))$,
if the strong Legendre condition holds along the extremal $(x(\cdot),p(\cdot),p^0,u(\cdot))$, that is, there exists $\epsilon_0 >0$ such that
$$
\left\{ \left\{ H, \frac{\partial H}{\partial u}.v \right\},\frac{\partial H}{\partial u}.v \right\}
+ \left\{ \frac{\partial^2 H}{\partial u^2}.(\dot{u},v),\frac{\partial H}{\partial u}.v \right\}
\leq - \epsilon_0 \|v\|^2, \quad \forall v\in \R^m,
$$
and if moreover the mapping $\frac{\partial f}{\partial u}(x_0,u(0)):\R^m \mapsto T_{x_0}M$ is one-to-one, then there exists $\epsilon_1 > 0$ such that $x(\cdot)$ is locally optimal in $L^\infty$ topology on $[0,\epsilon_1]$.
\end{itemize}
\end{thm}

\medskip

As we have seen, the Legendre (or generalized Legendre) condition is a necessary condition, while the strong (or strong generalized Legendre) condition is a sufficient condition.
However, these sufficient conditions are not easy to verify in practice. This leads to the next section where we explain how to use the so-called \emph{conjugate point} along the extremal to determine the time when the extremal is no longer optimal.

\paragraph{Conjugate points.}
We consider here the simplified problem ($\P_0$) with $M=\R^n$, $M_0=\left\{x_0\right\}$, $M_1=\left\{ x_1\right\}$, and $U=\R^m$.
Under the strict Legendre assumption assuming that the Hessian $ \frac{\partial^2 H}{\partial u^2} (x,p,p^0,u)$ is negative definite, the quadratic form $Q_{t_f}$ is negative definite if $t_f>0$ is small enough.

\begin{defi} \label{def_conj_p}
The first conjugate time is defined by the infimum of times $t>0$ such that $Q_t$ has a nontrivial kernel. We denote the first conjugate time along $x(\cdot)$ by $t_c$.
\end{defi}
The extremals are locally optimal (in $L^\infty$ topology) as long as we do not encounter any conjugate point.
Define the \emph{exponential mapping}
\begin{equation} \label{expmap}
\exp_{x_0}(t,p_0) = x(t,x_0,p_0),
\end{equation}
where the solution of \eqref{extremals} starting from $(x_0,p_0)$ at $t=0$ is denoted as $(x(t,x_0,p_0),p(t,x_0,p_0))$.
Then, we have the following result  (see, e.g., \cite{Agrachev,BonnardCT} for the proof and more precise results):
\begin{quote}
{\it
The time $t_c$ is a conjugate time along $x(\cdot)$ if and only if the mapping $\exp_{x_0}(t_c,\cdot)$ is not an immersion at $p_0$, i.e., the differential of the mapping $\exp_{x_0}(t_c,\cdot)$ is not injective.
}
\end{quote}
Essentially this result states that computing a first conjugate time $t_c$ reduces to finding the zero of some determinant along the extremal.  
In the smooth case (the control can be expressed as a smooth function of $x$ and $p$), the survey article \cite{BonnardCT} provides also some algorithms to compute first conjugate times.
In case of bang-bang control, a conjugate time theory has been developed (see \cite{SilvaT2010} for a brief survey of the approaches), but the computation of conjugate times remains difficult in practice (see, e.g., \cite{MaurerBKK2005}).

\medskip

When the singular controls are of order one (see Definition \ref{def_singorder}), the second-order optimality condition is sufficient for the analysis.
For higher order singular controls, higher order optimality conditions are needed which are recalled in the next section.

\subsubsection{Order of Singular Controls and Higher Order Optimality Conditions}
\label{hooc}
In this section we recall briefly the order of singular controls and the higher order optimality conditions. They will be used in Section \ref{example1} to analyze the example, which exhibits a singular control of order two.
It is worth noting that when the singular control is of \emph{order $1$} (also called \emph{minimal order} in \cite{BonnardChyba,Chitour}), these higher order optimality conditions are not required.

To illustrate how to use these conditions, we consider the minimal time control problem on $M$
\begin{equation} \label{pbtocp}
\left\{ \begin{split}
	& \min t_f , \\
	& \dot{x}(t) = f(x(t))+u_1(t) g_1(x(t))+u_2(t) g_2(x(t)),\quad u=(u_1,u_2) \\
	& \Vert u(t)\Vert^2 = u_1(t)^2+u_2(t)^2 \leq 1 , \\
	& x(0) = x_0,\ x(t_f) \in M_1 , \quad t_f\geq 0\ \textrm{free},
\end{split}\right.
\end{equation}
where $f$, $g_1$ and $g_2$ are smooth vector fields on $M$.
We assume that $M_1$ is accessible from $x_0$, and that there exists a constant $B_{t_f}$ such that for every admissible control $u$, the corresponding trajectory $x_u(t)$ satisfies $\| x_u(t) \| \leq B_{t_f}$ for all $t\in [0,t_f]$.
Then, according to classical results (see, e.g., \cite{Cesari,Trelatbook}), there exists at least one optimal solution $(x(\cdot),u(\cdot))$, defined on $[0,t_f]$.

\medskip

Let $h_0(x,p)=\langle p, f(x) \rangle $, $h_1(x,p)=\langle p, g_1(x) \rangle $, and $h_2(x,p)=\langle p, g_2(x) \rangle $.
According to the PMP (see Section \ref{subsubsec_langrange}), the Hamiltonian of the problem \eqref{pbtocp} is defined by  
$$
H(x,p,p^0,u) = h_0(x,p)+u_1 h_1(x,p) +u_2 h_2(x,p)+p^0
$$
where $p(\cdot)$ is the adjoint variable, and $p^0 \leq 0$ is a real number such that $(p(\cdot),p^0)\neq 0$.
Defining $\Phi(t)=(h_1(t),h_2(t))$, the maximization condition of the PMP yields
$$
u(t) = \frac{\Phi(t)}{\Vert\Phi(t)\Vert},
$$
almost everywhere on $[0,t_f]$, whenever $\Phi(t) \neq (0,0)$.

We call $\Phi$ (as well as its components) the \emph{switching function}.
We say that an arc (restriction of an extremal to a subinterval $I$) is \emph{regular} if $\Vert \Phi(t)\Vert \neq 0$ along $I$. Otherwise, the arc is said to be \emph{singular}.

\medskip

Following \cite{Gabasov}, we give here below a precise definition of the order of a singular control.
The use of Poisson (and Lie) brackets simplifies the formulation of the higher order optimality conditions. This is one of the reasons making geometric optimal control theory a valuable tool in practice.

\begin{defi} \label{def_singorder}
The singular control $u=(u_1,u_2)$ defined on a subinterval $I \subset [0,t_f]$ is said to be of \emph{order $q$} if
\begin{enumerate}
\item the first $(2q-1)$-th time derivatives of $h_i$, $i=1,2$, do not depend on $u$ and
$$
\frac{d^{k}}{dt^{k}}(h_i) = 0,\quad k=0,1,\cdots,2q-1,
$$
\item the $2q$-th time derivative of $h_i$, $i=1,2$, depends on $u$ linearly and
$$
\frac{\partial}{\partial u_i} \frac{d^{2q}}{dt^{2q}}(h_i) \neq 0, \quad
\det \left( \frac{\partial}{\partial u} \frac{d^{2q}}{dt^{2q}}\Phi \right) \neq 0,\quad
i=1,2,
$$
along $I$.
\end{enumerate}
The control $u$ is said to be of \emph{intrinsic order} $q$ if the vector fields satisfy also
$$
[g_i, \mathrm{ad}^k f.g_i] \equiv 0,\quad k=1,\cdots,2q-2,\quad i=1,2.
$$
\end{defi}

The condition of a nonzero determinant guarantees that the optimal control can be computed from the $2q$-th time derivative of the switching function. Note that this definition requires that the two components of the control have the same order.

\medskip

We next recall the Goh and generalized Legendre-Clebsch conditions (see \cite{Goh,Kelley,Krener}). It is worth noting that in \cite{Krener}, the following higher-order necessary conditions hold even when the components of the control $u$ have different orders.

\begin{lem} \label{necconds0}
\emph{(higher-order necessary conditions)}
We assume that a singular control $u=(u_1,u_2)$ defined on $I$ is of order $q$, that $u$ is optimal and not saturating, i.e., $\|u\|<1$. Then the Goh condition
$$
\frac{\partial}{\partial u_j} \frac{d^{k}}{dt^{k}}(h_i) = 0,\quad k=0,1,\cdots,2q-1,\quad i,j=1,2,\quad i \neq j,
$$
must be satisfied along $I$. Moreover, the matrix having as $(i,j)$-th component
$$
(-1)^q \frac{\partial}{\partial u_j} \frac{d^{2q}}{dt^{2q}}(h_i) ,\quad i,j=1,2,
$$
is symmetric and negative definite along $I$ (generalized Legendre-Clebsch Condition).
\end{lem}

In practice, it happens that the singular controls are often of intrinsic order $2$,
and that $ [ g_1,g_2 ] = 0$, $[ g_1,[f,g_2 ]] = 0$, and $[ g_2,[f,g_1 ]] = 0$.
The conditions given in the above definition yield $[ g_1,[f,g_1 ]] = 0$, $[ g_2,[f,g_2 ]] = 0$, $[ g_1,\mathrm{ad}^2 f.g_1 ] = 0$, $[ g_2,\mathrm{ad}^2 f.g_2] = 0$,
$\langle p,  [ g_1,\mathrm{ad}^3 f.g_1 ](x) \rangle \neq 0$, $\langle p,  [ g_2,\mathrm{ad}^3 f.g_2](x) \rangle \neq 0$, and
$$
\langle p,  [ g_1,\mathrm{ad}^3 f.g_1 ](x) \rangle \langle p,  [ g_2,\mathrm{ad}^3 f.g_2 ](x) \rangle
- \langle p,  [ g_2,\mathrm{ad}^3 f.g_1 ](x) \rangle \langle p,  [ g_1,\mathrm{ad}^3 f.g_2 ](x) \rangle \neq 0,
$$
We have thus the following higher-order necessary conditions, that will be used on the example in Section \ref{example1}.

\begin{cor} \label{necconds}
We assume that the optimal trajectory $x(\cdot)$ contains a singular arc, defined on the subinterval $I$ of $[0,t_f]$, associated with a non saturating control $u=(u_1,u_2)$ of intrinsic order $2$. If the vector fields satisfy $[g_1,g_2] = 0$, $[g_i,[f,g_j]]  = 0$, for $i,j=1,2$, then the Goh condition
$$
\langle p(t), [g_1,\mathrm{ad} f.g_2] (x(t)) \rangle = 0, \quad
\langle p(t), [g_1,\mathrm{ad}^2 f.g_2] (x(t)) \rangle = \langle p(t), [g_2,\mathrm{ad}^2 f.g_1] (x(t)) \rangle = 0,
$$
and the generalized Legendre-Clebsch condition (in short, GLCC )
$$
\langle p(t),[g_i,\mathrm{ad}^3f.g_i](x(t))\rangle \leq 0, \quad i=1,2,
$$
$$
\langle p(t),[g_1,\mathrm{ad}^3f.g_2](x(t))\rangle = \langle p(t),[g_2,\mathrm{ad}^3f.g_1](x(t))\rangle
$$
must be satisfied along $I$. Moreover, we say that the strengthened GLCC  is satisfied if we have a strict inequality above, that is, $\langle p(t),[g_i,\mathrm{ad}^3f.g_i](x(t))\rangle < 0$.
\end{cor}

\bigskip

In the next section, we recall the chattering phenomenon that may happen in the optimal control problem. This phenomenon is actually not rare as illustrated in \cite{Zelikin1} by many examples (in astronautics, robotics, economics, and etc.).
These examples are mostly single input systems.
The existence of chattering phenomenon for bi-input control affine systems is also proved in \cite{ztc2}.

\subsection{Chattering phenomenon}
\label{chatphenom}
\begin{figure}[h]
\centering
	\includegraphics[scale=1]{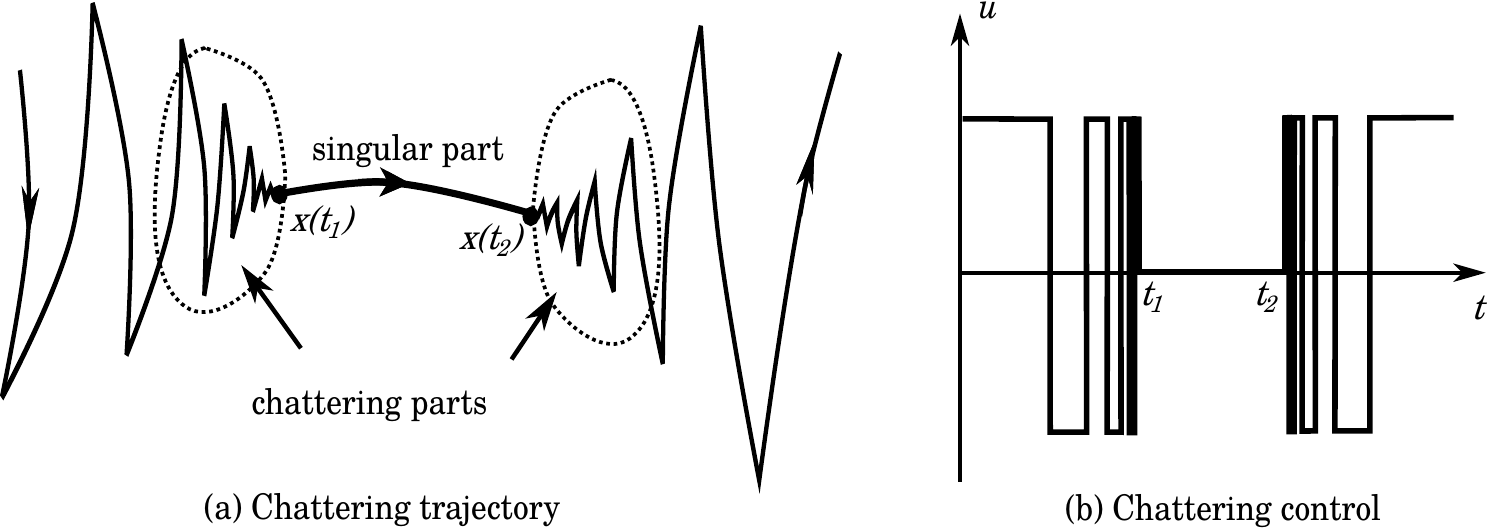}
	\caption{An illustration of chattering phenomenon.}
	\label{chattering}
\end{figure}

We call \textit{chattering phenomenon} (or Fuller's phenomenon) the situation when the optimal control switches an infinite number of times over a compact time interval. It is well known that, if the optimal trajectory involves a singular arc of higher order, then no connection with a bang arc is possible and the bang arcs asymptotically joining the singular arc must chatter. On Figure \ref{chattering}(b), the control is singular over $(t_1,t_2)$, and the control $u(t)$ with $t \in (t_1-\epsilon_1,t_1) \cup (t_2,t_2+\epsilon_2)$, $\epsilon_1>0$, $\epsilon_2>0$ is chattering. The corresponding optimal trajectory is called a chattering trajectory. On Figure \ref{chattering}(a), the chattering trajectory ``oscillates'' around the singular part and finally ``gets off" the singular trajectory with an infinite number of switchings.

The chattering phenomenon is illustrated by the Fuller's problem (see  \cite{Fuller,Marchal}), which is the optimal control problem
\begin{equation*}
\begin{split}
& \min \int_0^{t_f} x_1(t)^2 \, dt , \\
& \dot{x}_1(t)=x_2(t),\  \dot{x}_2(t)=u(t), \\
& |u(t)| \leq 1,\\
& x_1(0)=x_{10},\ x_2(0)=x_{20},\\
& x_1(t_f)=0,\ x_2(t_f)=0,\qquad t_f\ \textrm{free}.
\end{split}
\end{equation*}
We define $ \xi = \left( \frac{\sqrt{33}-1}{24} \right)^{1/2}$ as the unique positive root of the equation $\xi^4+\xi^2/12-1/18=0$, and we define the sets
\begin{equation*}
\begin{split}
\Gamma_{+}=\{ (x_1,x_2)\in \R^2 \mid\ & x_1= \xi x_2^2,\ x_2<0 \} , \\
&R_{+}=\{ (x_1,x_2)\in \R^2 \mid x_1 <  - \mathrm{sign} (x_2) \xi x_2^2 \} , \\
\Gamma_{-}=\{ (x_1,x_2)\in \R^2 \mid\ & x_1= - \xi x_2^2,\ x_2>0 \} , \\
&R_{-}=\{ (x_1,x_2)\in \R^2 \mid x_1 >  - \mathrm{sign} (x_2) \xi x_2^2 \} .
\end{split}
\end{equation*}
The optimal synthesis of the Fuller's problem yields the following feedback control (see \cite{Fuller,Schattler,Wonham}).
\begin{equation*}
	u^{\ast}=\begin{cases}
		\phantom{-}1 & \textrm{if}\ x \in R_{+} \bigcup \Gamma_{+} , \\
		-1& \textrm{if}\ x \in R_{-} \bigcup \Gamma_{-}  .
	  \end{cases}
\end{equation*}
The control switches from $u=1$ to $u=-1$ at points on $\Gamma_{-}$ and from $u=-1$ to $u=1$ at points on $\Gamma_{+}$. The corresponding trajectories crossing the switching curves $\Gamma_{\pm}$ transversally are chattering arcs with an infinite number of switchings that accumulate with a geometric progression at the final time $t_f>0$.

\begin{figure}[h]
\centering
	\includegraphics[scale=0.5]{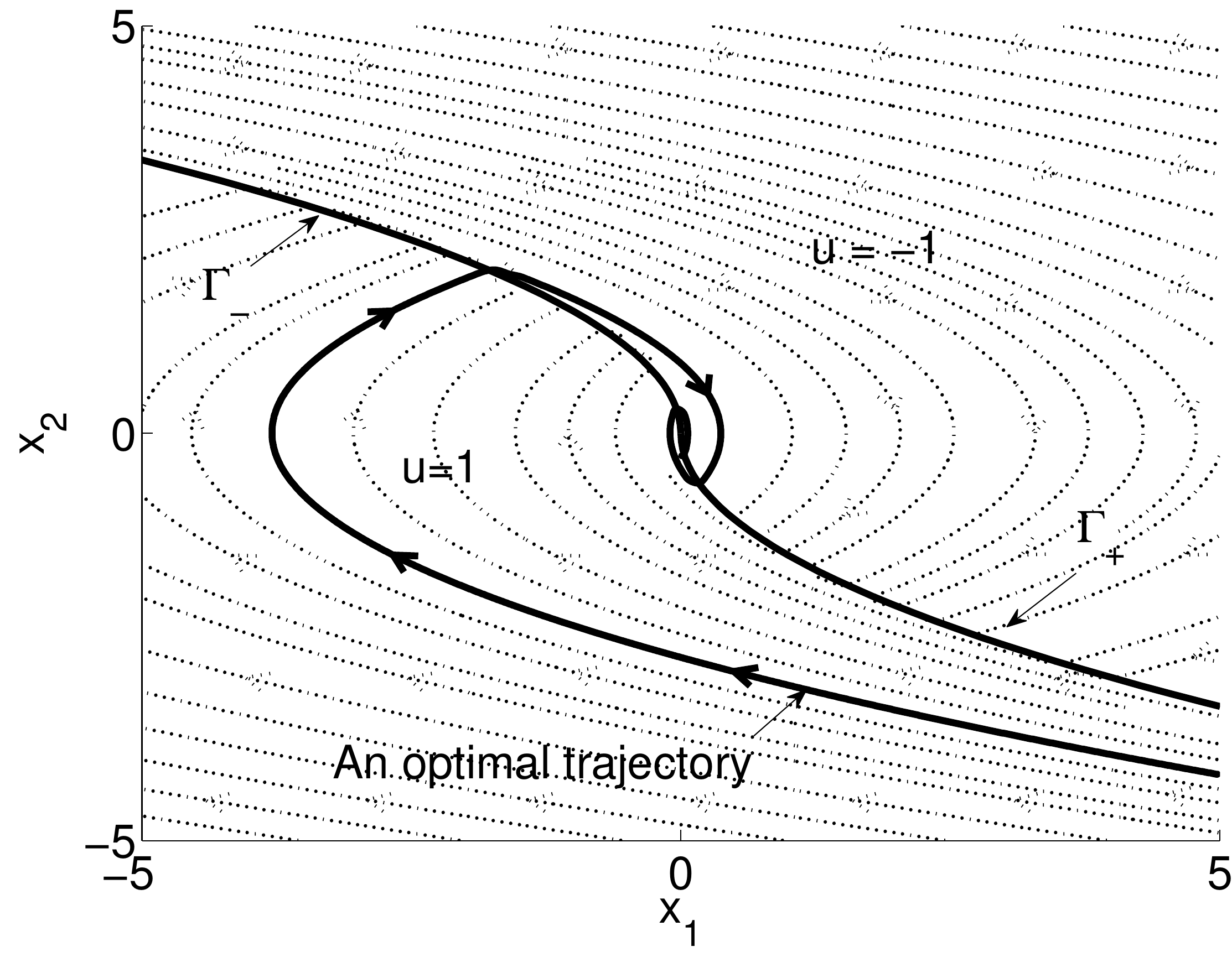}
	\caption{Optimal synthesis for the Fuller's problem.}
	\label{Fuller}
\end{figure}

The optimal synthesis for the Fuller's problem is drawn on Figure \ref{Fuller}.
The optimal control of the Fuller's problem, denoted $u^{\ast}$, contains a countable set of switchings of the form
\begin{equation*}
	u^{\ast}(t)=\begin{cases}
		\phantom{-}1 & \textrm{if}\ t \in [t_{2k},t_{2k+1}), \\
		-1& \textrm{if}\ t \in [t_{2k+1},t_{2k+2}] ,
	  \end{cases}
\end{equation*}
where $\left\{ t_{k} \right\}_{k \in \N}$ is a set of switching times that satisfies ${ (t_{i+2} - t_{i+1}) < (t_{i+1} - t_{i})}$, $i \in \N$ and converges to $t_f < +\infty$. This means that the chattering arcs contain an infinite number of switchings within a finite time interval $t_f>0$.

\section{Numerical Methods in Optimal Control}
\label{sec_NMOC}
Numerical approaches in optimal control are usually distinguished between direct and indirect methods.
Indirect methods consist in solving numerically the boundary value problem derived from the application of the PMP.
Direct methods consist in discretizing the state and the control, and solving the resulting nonlinear optimization problem.
The principles of both methods are recalled hereafter.

\subsection{Indirect Methods}
\label{subsec_indir}
In indirect approaches, the Pontryagin Maximum Principle (first-order necessary condition for optimality) is applied to the optimal control problem in order to express the control as a function of the state and the adjoint. This reduces the problem to a nonlinear system of $n$ equations with $n$ unknowns generally solved by Newton-like methods.
Indirect methods are also called shooting methods.
The principle of the simple shooting method and of the multiple shooting method are recalled. The problem considered in this section is ($\P_0$).

\paragraph{Simple shooting method.}
Using \eqref{Hamiltonian}, the optimal control can be expressed as a function of the state and the adjoint variable $(x(t),p(t))$.
Denoting $z(t)=(x(t),p(t))$, the extremal system \eqref{extremals} can be written under the form $\dot{z}(t)=F(z(t))$. The initial and final conditions \eqref{tercond_general}, the transversality conditions \eqref{condtrans_adjoint}, and the transversality condition on the Hamiltonian \eqref{condtrans_hamil} can be written under the form of $R(z(0),z(t_f),t_f)=0$. We thus get a two boundary value problem
$$
\dot{z}(t) = F(t,z(t)),\quad R(z(0),z(t_f),t_f)=0.
$$
Let $z(t,z_0)$ be the solution of the Cauchy problem 
$$\dot{z}(t) = F(t,z(t)),\qquad z(0)=z_0.$$
Then this two boundary value problem consists in finding a zero of the equation 
$$R(z_0,z(t_f,z_0),t_f)=0.$$
This problem can be solved by Newton-like methods or other iterative methods.

\paragraph{Multiple shooting method.}
\label{subsec_tirmultiple}
The drawback of the single shooting method is the sensitivity of the Cauchy problem to the initial condition $z_0$ . The multiple shooting aims at a better numerical stability by dividing the interval $[0,t_f]$ into $N$ subintervals $[t_i,t_{i+1}]$ and considering as unknowns the values of $z_i = (x(t_i),p(t_i))$ at the beginning of each subinterval.  
The application of the PMP to the optimal control problem yields a multi-point boundary value problem, which consists in finding $Z=(p(0),t_f,z_i)$, $i=1,\cdots,N-1$ such that the differential equation
\begin{equation*} 
\dot{z}_i (t)
= F(t,z(t))
= \begin{cases}
F_0(t,z(t)),\quad t_0 \leq t \leq t_1, \\
F_1(t,z(t)),\quad t_1 \leq t \leq t_2, \\
\cdots, \\
F_{N-1}(t,z(t)),\quad t_{N-1} \leq t \leq t_N, \\
\end{cases}
\end{equation*}
and the constraints
$$
x(0)\in M_0,
\quad x(t_f) \in M_1,
\quad p(0) \perp T_{x(0)}M_0,
$$
$$
p(t_f) - p^0 \frac{\partial g}{\partial x}(t_f,x(t_f))\perp T_{x(t_f)} M_1,
\quad H(t_f) = 0, \quad z(t_i^-)=z(t_i^+),\quad i=1,\cdots,N-1,
$$
are satisfied. 
The nodes of the multiple shooting method may involve the switching times (at which the switching function changes sign), and the junction times (entry, contact, or exit times) with boundary arcs.
In this case an a priori knowledge of the solution structure is required.

The multiple shooting method improves the numerical stability at the expense of a larger nonlinear system. An adequate node number must be chosen making a compromise between the system dimension and the convergence domain.

\subsection{Direct Methods}
\label{subsec_dir}
Direct methods are so called because they address directly the optimal control problem without using the first-order necessary conditions yielded by the PMP.
By discretizing both the state and the control, the problem reduces to a nonlinear optimization problem in finite dimension, also called NonLinear Programming problem (NLP).
The discretization may be carried out in many ways, depending on the problem features. As an example we may consider a subdivision $0=t_0<t_1<\cdot<t_N=t_f$ of the interval $[0,t_f]$. We discretize the controls such that they are piecewise constant on this subdivision with values in $U$. Meanwhile the differential equations may be discretized by an explicit Euler method : by setting $h_i=t_{i+1}-t_i$, we get $x_{i+1}=x_i+h_i f(t_i,x_i,u_i)$. The cost may be discretized by a quadrature procedure.
These discretizations reduces the optimal control problem $\P_0$ to a nonlinear programming problem of the form
\begin{equation*}
\begin{split}
\min \{ C(x_0,\cdots,x_N,u_0,\cdots,u_N) | & x_{i+1}=x_i+h_i f(t_i,x_i,u_i),\\
& u_i \in U, \, i=1,\cdots,N-1, \, x_0\in M_0,\, x_N \in M_1 \}.
\end{split}
\end{equation*}

From a more general point of view, a finite dimensional representation of the control and of the state has to be chosen such that the differential equation, the cost, and all constraints can be expressed in a discrete way.

The numerical resolution of a nonlinear programming problem is standard, by gradient methods, penalization, quasi-Newton, dual methods, etc. (see, e.g., \cite{BockPlitt,Gerdts,Kirches,StoerBulirsch}).
There exist many efficient optimization packages such as \texttt{IPOPT} (see \cite{IPOPT}), \texttt{MUSCOD-II} (see \cite{Diehl}), or the \texttt{Minpack} project (see \cite{More}) for many optimization routines.

Alternative variants of direct methods are the collocation methods, the spectral or pseudo-spectral methods, the probabilistic approaches, etc.

\medskip

Another approach to optimal control problems that can be considered as a direct method, consists in solving the Hamilton-Jacobi equation satisfied (in the viscosity sense) by the value function which is of the form 
$$\frac{\partial S}{\partial t}+H_r\left(x,\frac{\partial S}{\partial x}\right)=0.
$$
The value function is the optimal cost for the optimal control problem starting from a given point $(x,t)$ (see \cite{Sethian} for some numerical methods).

\subsection{Comparison Between Methods}
The main advantages and disadvantages of the direct and indirect methods are summarized in Table \ref{comparisons} (see also, e.g., \cite{Trelatbook,Trelatsurvey}).

\begin{table}[h]
\centering
\begin{tabular}{lll}
\hline
\hline
                  & Direct methods & Indirect methods  \\
\hline                  
 a priori knowledge of the solution structure & not required & required  \\
\hline
 sensible to the initial condition  & not sensible & very sensible \\
\hline
 handle the state constraints & easy & difficult \\
\hline
 convergence speed and accuracy & relatively slow and inaccurate & fast and accurate \\
\hline
 computational aspect & memory demanding & parallelizable\\
\hline
\hline
\end{tabular}
\caption{Pros and cons for direct and indirect methods}
\label{comparisons}
\end{table}

In practice no approach is intrinsically better than the other. The numerical method should be chosen depending on the problem features and on the known properties of the solution structure.
These properties are derived by a theoretical analysis using the geometric optimal control theory.
When a high accuracy is desired, as is generally the case for aerospace problems, indirect methods should be considered although they require more theoretical insight and may raise numerical difficulties.

Whatever the method chosen, there are many ways to adapt it to a specific problem (see \cite{Trelatsurvey}). Even with direct methods, a major issue lies in the initialization procedure.
In recent years, the numerical continuation has become a powerful tool to overcome this difficulty. The next section recalls some basic mathematical concepts of the continuation approaches, with a focus on the numerical implementations of these methods.

\section{Continuation Methods}
\label{sec_CM}
\subsection{Existence Results and Discrete Continuation}
The basic idea of continuation (also called homotopy) methods is to solve a difficult problem step by step starting from a simpler problem by parameter deformation.
The theory and practice of the continuation methods are well-spread (see, e.g., \cite{Allgower,Rheinboldt,Watson}).
Combined with the shooting problem derived from the PMP, a continuation method consists in deforming the problem into a simpler one (that can be easily solved) and then solving a series of shooting problems step by step to come back to the original problem.

One difficulty of homotopy methods lies in the choice of a sufficiently regular deformation that allows the convergence of the homotopy method. The starting problem should be easy to solve, and the path between this starting problem and the original problem should be easy to model.
Another difficulty is to numerically follow the path between the starting problem and the original problem. This path is parametrized by a parameter denoted $\lambda$.
When the homotopic parameter $\lambda$ is a real number and when the path is linear\footnote{meaning that in some coordinates, for $\lambda \in[0,1]$, the path consists in a convex combination of the simpler problem and of the original problem}  in $\lambda$, the homotopy method is rather called a continuation method.

The choice of the homotopic parameter may require considerable physical insight into the problem. This parameter may be defined either artificially according to some intuition, or naturally by choosing physical parameters of the system, or by a combination of both.

\medskip

Suppose that we have to solve a system of $N$ nonlinear equations in $N$ dimensional variable $Z$
$$
F(Z) = 0,
$$
where $F: \R^{N} \mapsto \R^{N}$ is a smooth map.
We define a deformation
$$
G:\R^{N} \times [0,1] \mapsto \R^{N},
$$
such that
$$
G(Z,0)=G_0(Z), \quad G(Z,1)=F(Z),
$$
where $G_0 : \R^{N} \mapsto \R^{N}$ is a smooth map having known zero points.

A \emph{zero path} is a curve $c(s) \in G^{-1}(0)$ where $s$ represents the arc length.
We would like to trace a zero path starting from a point $Z_0$ such that $G(Z_0,0)=0$ and ending at a point $Z_f$ such that $G(Z_f,1)=0$.

The first question to address is the existence of zero paths, since the feasibility of the continuation method lies on this assumption. The second question to address is how to numerically track such zero paths when they exist.

\medskip

\paragraph{Existence of zero paths}

The local existence of the zero paths is answered by the implicit function theorem. Some regularity assumptions are needed, as in the following statement (which is the contents of \cite[Theorem 2.1]{GarciaZangwill}).

\begin{thm} \label{thm_existencechemins}
\emph{(Existence of the zero paths)}
Let $\Omega$ be an open bounded subset of $\R^{N}$ and let the mapping $G:\Omega \times [0,1] \mapsto \R^{N}$ be continuously differentiable such that:
\begin{itemize}
\item Given any $(Z,\lambda) \in \left\{(Z,\lambda) \in \Omega \times [0,1] \ \mid\ G(Z,\lambda)=0 \right\}$, the Jacobian matrix
$$
G^\prime = \left( \frac{\partial G}{\partial Z_1},\cdots,\frac{\partial G}{\partial Z_{N}},\frac{\partial G}{\partial \lambda} \right),
$$
is of maximum rank $N$;
\item Given any $Z \in \left\{ Z\in\Omega \ \mid\  G(Z,0)=0 \right\} \cup \left\{ Z\in\Omega \ \mid\  G(Z,1)=0 \right\}$, the Jacobian matrix 
$$G^\prime = \left( \frac{\partial G}{\partial Z_1},\cdots,\frac{\partial G}{\partial Z_{N}}\right)$$
 is of maximum rank $N$;
\end{itemize}
Then $\left\{ (Z,\lambda) \in \Omega \times [0,1] \ \mid\ G(Z,\lambda)=0 \right\}$ consists of the paths that is either a loop in $\bar{\Omega} \times [0,1]$ or starts from a point of $\partial \bar{\Omega} \times [0,1]$ and ends at another point of $\partial \bar{\Omega} \times [0,1]$, where $\partial \bar{\Omega}$ denotes the boundary of $\bar{\Omega}$.
\end{thm}

\begin{figure}[h]
\centering
	\includegraphics[scale=0.45]{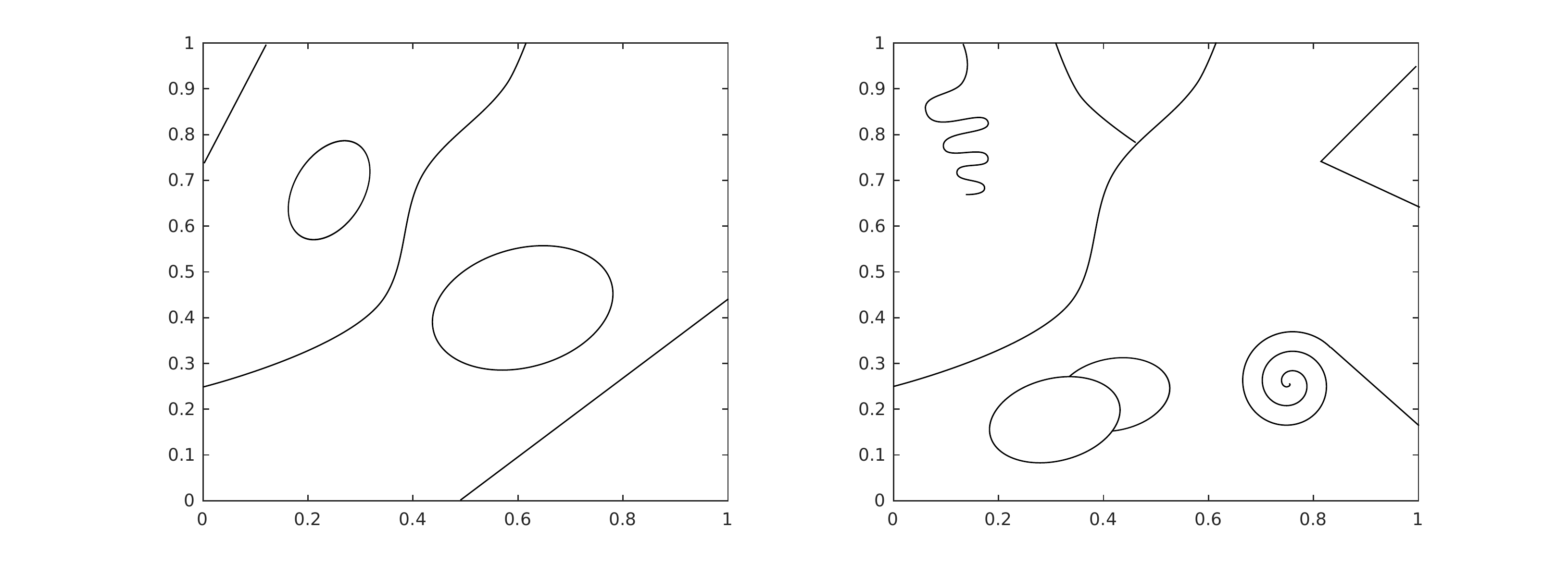}
	\caption{Possible zero paths (left) and impossible zero paths (right).}
	\label{zeropath}
\end{figure}

This means that the zero path is diffeomorphic to a circle or the real line.
The possible paths and impossible paths are shown in Figure \ref{zeropath} (borrowed from \cite{GarciaZangwill, Gergaud}).

\medskip

Now we provide basic arguments showing the feasibility of the continuation method (see Section 4.1 of \cite{Trelatsurvey} for more details).

Consider the simplified optimal control problem $\P_0$ with $M=\R^n$, $M_0=\left\{x_0\right\}$, $M_1=\left\{ x_1 \right\}$ and $U=\R^m$.
We assume that the real parameter $\lambda \in [0,1]$ is increasing monotonically from $0$ to $1$.
Under these assumptions, we are to solve a family of optimal control problems parameterized by $\lambda$, i.e.,
\begin{equation} \label{genpb_simplified}
\min_{E_{x_0,t_f,\lambda}(u_\lambda)=x_1} C_{t_f,\lambda}(u),
\end{equation}
where $E$ is the end-point mapping defined in Definition \ref{def_endpmap}.

We assume moreover that, along the continuation procedure:
\begin{itemize}
\item[(1)] there are no minimizing abnormal extremals; 
\item[(2)] there are no minimizing singular controls: by Definition \ref{def_sin_E}, the control $u$ is not singular means that the mapping $dE_{x_0,t_f,\lambda}(u)$ is surjective; 
\item[(3)] there are no conjugate points (by Definition \ref{def_conj_p} the quadratic form $Q_{t_f}$ is not degenerate). The absence of conjugate point can be numerically tested (see, e.g., \cite{BonnardCT}).
\end{itemize}
We will see that these assumptions are essential for the local feasibility of the continuation methods.

According to the Lagrange multipliers rule, especially the first-order condition \eqref{lagcond}, if $u_\lambda$ is optimal, then there exists $(\psi_\lambda,\psi^0_\lambda) \in \R^n \times \R \backslash \left\{(0,0)\right\}$ such that
$\psi_\lambda dE_{x_0,t_f,\lambda}(u_\lambda) + \psi^0_\lambda dC_{t_f,\lambda}(u)=0$.
Since we have assumed that there are no minimizing abnormal extremals in the problem and $(\psi_\lambda,\psi^0_\lambda)$ is defined up to a multiplicative scalar, we can set $\psi^0_\lambda=-1$.
Defining the Lagrangian by
$$L_{t_f,\lambda}(u,\psi) = \psi_\lambda E_{x_0,t_f,\lambda}(u) - C_{t_f,\lambda}(u),$$
we seek $(u_\lambda,\psi_\lambda)$ such that
$$
G(u,\psi,\lambda) =
\left( \begin{matrix}
\frac{\partial L_{t_f,\lambda}}{\partial u}(u,\psi)\\
E_{x_0,t_f,\lambda}(u) - x_1
\end{matrix}
\right)=0.
$$

Let $(u_{\bar{\lambda}},\psi_{\bar{\lambda}},\bar{\lambda})$ be a zero of $G$ and assume that $G$ is of class $C^1$. Then according to Theorem \ref{thm_existencechemins}, we require the Jacobian of $G$ with respect to $(u,\psi)$ at the point $(u_{\bar{\lambda}},\psi_{\bar{\lambda}},\bar{\lambda})$ to be invertible. More precisely, the Jacobian of $G$ is
\begin{equation} \label{mat_sensitive}
\left(
\begin{matrix}
Q_{t_f,\lambda} & dE_{x_0,t_f,\lambda}(u)^\ast \\
dE_{x_0,t_f,\lambda}(u) & 0
\end{matrix}
\right),
\end{equation}
where $Q_{t_f,\lambda}$ is the Hessian $\frac{\partial^2 L_{t_f,\lambda}}{\partial^2 u}(u,\psi,\psi^0)$ restricted to $\ker \frac{\partial L_{t_f,\lambda}}{\partial u}$, and $dE_{x_0,t_f,\lambda}(u)^\ast$ is the transpose of $dE_{x_0,t_f,\lambda}(u)$.

We observe that the matrix \eqref{mat_sensitive} is invertible if and only if the linear mapping $dE_{x_0,t_f,\lambda}(u)$ is surjective and the quadratic form $Q_{t_f,\lambda}$ is non-degenerate.
These properties correspond to the absence of any minimizing singular control and conjugate points, which are the assumptions done for the local feasibility of the continuation procedure.

The implicit function argument above is done on the control. In practice the continuation procedure is rather done on the exponential mapping (see \eqref{expmap}) and it consists in tracking a path of initial adjoint vectors $p_{0,\lambda}$.
Therefore we parameterize the exponential mapping by $\lambda$, and thus problem \eqref{genpb_simplified} is to solve
\begin{equation} \label{expmap}
\exp_{x_0,\lambda}(t_f,p_{0,\lambda}) = x_1.
\end{equation}
On the one hand, according to the PMP, the optimal control $u$ satisfies the extremal equations \eqref{Hamiltonian}, and thus $u_\lambda = u_\lambda(t,p_{0,\lambda})$ is a function of the initial adjoint $p_{0,\lambda}$.
On the other hand, the Lagrange multipliers are related to the adjoint vector by $p(t_f) = \psi$, and thus $\psi_\lambda = \psi_\lambda (p_{0,\lambda})$ is also a function of $p_{0,\lambda}$.
Therefore, the shooting function defined by $S(p_0,\lambda) = G(u(p_0),\psi(p_0),\lambda)$ has an invertible Jacobian if the matrix \eqref{mat_sensitive} is invertible.
We conclude then that the assumptions (1)-(3) mentioned above are sufficient to ensure the local feasibility.

Despite of local feasibility, the zero path may not be globally defined for any $\lambda \in [0,1]$. The path could cross some singularity or diverge to infinity before reaching $\lambda = 1$.

The first possibility can be eliminated by assuming (2) and (3) over all the domain $\Omega$ and for every $\lambda \in [0,1]$.
The second possibility is eliminated if the paths remain bounded  or equivalently by the properness of the exponential mapping (i.e., the initial adjoint vectors $p_{0,\lambda}$ that are computed along the continuation procedure remain bounded uniformly with respect to $\lambda$).
According to \cite{BonnardTrelat2,Trelat2}, if the exponential mapping is not proper, then there exists an abnormal minimizer.
By contraposition, if one assumes the absence of minimizing abnormal extremals, then the required boundedness follows.

For the simplified problem \eqref{genpb_simplified}, where the controls are unconstrained and the singular trajectories are the projections of abnormal extremals, if there are no minimizing singular trajectory nor conjugate points over $\Omega$, then the continuation procedure \eqref{expmap} is globally feasible on $[0,1]$.

\medskip

In more general homotopy strategies, the homotopic parameter $\lambda$ is not necessarily increasing monotonically from $0$ to $1$. There may be turning points (see, e.g., \cite{Watson}) and it is preferable to parametrize the zero path by the arc length $s$.
Let
$c(s) = (Z(s),\lambda(s))$
be the zero path such that $G(c(s)) = 0$.
Then, a turning point of order one is the point where $\lambda^\prime (\bar{s})= 0$, $\lambda^{\prime\prime} (\bar{s}) \neq 0$.
In \cite{CaillauDaoud}, the authors indicate that if $\lambda =\lambda(\bar{s})$ is a turning point of order one, then the corresponding final time $t_f$ is a conjugate time, and the corresponding point $E_{x_0,t_f,\lambda}(u(x_0,p_0,t_f,\lambda))$ is the corresponding conjugate point \footnote{There, the end-point mapping has been implemented with the exponential mapping $E_{x_0,t_f,\lambda}(u)=\exp_{x_0,\lambda}(t_f,p_0)$ with initial condition $(x(0),p(0))=(x_0,p_0)$.}.
By assuming the absence of conjugate points over $\Omega$ for all $\lambda \in [0,1]$, the possibility of turning points is discarded.

Unfortunately, assuming the absence of singularities is in general too strong, and weaker assumptions do not allow concluding to the feasibility of the continuation method.
In the literature, there are essentially two approaches to tackle this difficulty.
The first one is of local type. One detects the singularities or bifurcations along the zero path (see, e.g., \cite{Allgower}).
The second one is of global type, concerning the so-called globally convergent probability-one homotopy method.
We refer the readers to \cite{ChowMY,Watson} for more details concerning this method.

\paragraph{Numerical tracking the zero paths.}
There exists many numerical algorithms to track a zero path.
Among these algorithms, the simplest one is the so called \emph{discrete continuation} or embedding algorithm.
The continuation parameter denoted $\lambda$, is discretized by $0= \lambda^0 < \lambda^1 < \cdots < \lambda^{n_l} = 1$ and the sequence of problems $G(Z,\lambda^i)=0$, $i=1,\cdots,n_l$ is solved to end up with a zero point of $F(Z)$.
If the increment $\triangle \lambda = \lambda^{i+1} - \lambda^i$ is small enough, then the solution $Z^i$ associated to $\lambda^i$ such that $G(Z^i,\lambda^i)=0$ is generally close to the solution of $G(Z,\lambda^{i+1})=0$.
The discrete continuation algorithm is detailed in Algorithm \ref{algrithm_discon}.

\begin{algorithm}
\SetAlgoLined

\KwResult{The solution of the discrete continuation}
initialization $Z=Z_0$, $\lambda^{0}=0$, $\triangle\lambda \in (\triangle \lambda_{min},\triangle \lambda_{max})$\;
\While{$\lambda\leq 1$ and $\triangle \lambda_{min}\leq \triangle\lambda \leq \triangle\lambda_{max} $}
{
  $\triangle\lambda=\min(\triangle\lambda,1-\lambda)$\;
  $\tilde{\lambda}=\lambda+\triangle\lambda$\;
  Find the solution $\bar{Z}$ such that $G(\bar{Z},\tilde{\lambda})=0$\;
\eIf{successful}{
  $Z=\tilde{Z}$\;
  $\lambda=\bar{\lambda}$\;
  $\bigtriangleup\lambda=2\bigtriangleup\lambda$\;}{
  $\bigtriangleup\lambda=\bigtriangleup\lambda/2$\;}
  }
\eIf{successful}{The discrete continuation is successful\;}{
The discrete continuation has failed\;}
\caption{Discrete continuation algorithm}
\label{algrithm_discon}
\end{algorithm}

In some cases the parameter $\lambda$ may be ill suited to parameterize the zero path, and thus causes a slow progress or even a failure of the discrete continuation.
Two enhancements (predictor-corrector methods and piecewise-Linear methods) have been proposed in the literature.

\subsection{Predictor-Corrector (PC) Continuation}
\label{subsec_pccontinuation}
A natural parameter for the zero curve $(Z,\lambda)$ is the arc-length denoted $s$.

The zero curve parameterized by the arc length $s$ is denoted
$$
c(s)=(Z(s),\lambda(s)).
$$
Differentiating
$G(Z(s),\lambda(s))=0$ with respect to $s$, we obtain
\begin{equation} \label{pathfow}
J_G \,\, t(J_G) = 0, \quad \|t(J_G)\|=1, \quad c(Z(0),0)=(Z(0),0),
\end{equation}
where $J_G = \frac{\partial G(Z(s),\lambda(s))}{\partial (Z,\lambda)}$ is the Jacobian, and $t(J_G)=\frac{dc(s)}{ds}$ is the tangent vector of the zero path $c(s)$.

\begin{figure}[h]
\centering
	\includegraphics[scale=0.75]{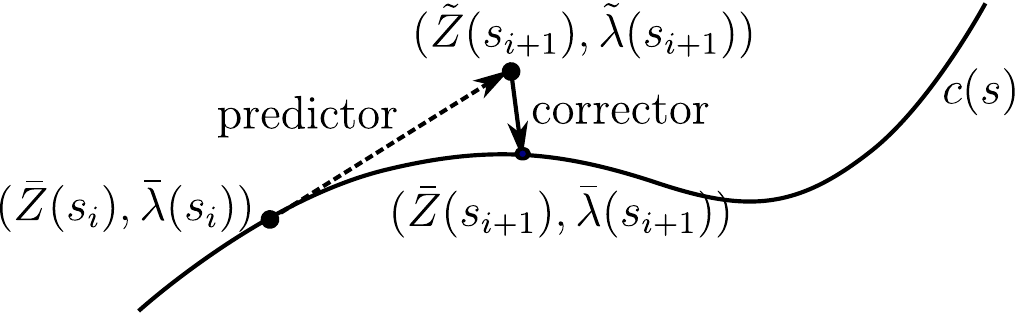}
	\caption{PC continuation.}
	\label{PCmethod}
\end{figure}

If we know a point of this curve $(\bar{Z}(s_i),\bar{\lambda}(s_i))$, and assuming that $c(s)$ is not a critical point (i.e., $t(J_G)$ is not null),we can predict a new zero point $(\tilde{Z}(s_{i+1}),\tilde{\lambda}(s_{i+1}))$ by
\begin{equation} \label{cont_pred}
(\tilde{Z}(s_{i+1}),\tilde{\lambda}(s_{i+1})) = (Z(s_i),\lambda(s_i)) + h_s \, t(J_G),
\end{equation}
where $h_s$ is the step size on $s$.
As shown in Figure \ref{PCmethod}, if the step size $h_s$ is sufficiently small, the prediction step yields a point $(\tilde{Z}(s_{i+1}),\tilde{\lambda}(s_{i+1}))$ close to a point $(\bar{Z}(s_{i+1}),\bar{\lambda}(s_{i+1}))$ on the curve, such that $G(c(s_{i+1}))=G(\bar{Z}(s_{i+1}),\bar{\lambda}(s_{i+1}))=0$. The correction step consists in coming back on the curve using a Newton-like method.
The PC continuation is described by Algorithm \ref{algrithm_pc}.
\begin{algorithm}
\SetAlgoLined
\KwResult{The solution of the PC continuation}
initialization $Z=Z_{0}$, $h_s >0$, $\lambda^{0}=0$, $\triangle\lambda\in(\triangle\lambda_{min},\triangle\lambda_{max})$\;
\While{$\lambda \leq 1$ and $\triangle\lambda_{min}\leq \triangle\lambda \leq \triangle\lambda_{max}$}
{
  (Predictor) Predict a point $(\tilde{Z},\tilde{\lambda})$ according to \eqref{cont_pred}\;
  (Corrector) Find the solution $(\bar{Z},\bar{\lambda})$ to $G(\tilde{Z},\tilde{\lambda})=0$\;
\eIf{successful}{
  $(Z,\lambda)=(\bar{Z},\bar{\lambda})$\;
  Increase the step length $h_s$\;}{
  Reduce the step length $h_s$\;}
}
\eIf{successful}{
   The continuation is successful\;}{
   The continuation has failed\;}
\caption{Prediction-Corrector continuation}
\label{algrithm_pc}
\end{algorithm}

When the optimal control problem is regular (in the sense of the Legendre condition are defined) and the homotopic parameter is a scalar, one can use the so called \emph{differential continuation} or \emph{differential pathfollowing}. This method consists in integrating accurately $t(J_G)$ satisfying \eqref{pathfow} (see details in \cite{CaillauCotsGergaud}). The correction step is replaced by the mere integration of an ordinary differential equation with the help of automatic differentiation (see, e.g., \cite{BischofCKM,CaillauN}).

\subsection{Piecewise-Linear (PL) Continuation}
The main advantage of the PL method is that it only needs the zero paths to be continuous (smoothness assumption of $G$ is not necessary).
For a detailed description of the PL methods, we refer the readers to \cite{Allgower,AllgowerGeorg,Gergaudthesis}.

Here we present the basic idea of the PL methods, which are also referred to as a \emph{simplicial methods}.
A PL continuation consists of following exactly a piecewise-linear curve $c_{\T}(s)$ that approximates the zero path $c(s) \in G^{-1}(0)$.

The approximation curve $c_\T(s)$ is a polygonal path relative to an underlying triangulation $\T$ of $\R^{N+1}$, which is a subdivision of $\R^{N+1}$ into $(N + 1)$-simplices.
\footnote{Let $v_1, \cdots, v_{j+1} \in \R^{N+1}$, $j \leq N+1$, be affinely independent points, i.e., $v_k - v_1$, $k=2,\cdots,j+1$ are linearly independent. A \emph{j-simplex} in $\R^{N+1}$ is defined by the convex hull of the set ${v_1, \cdots,v_{j+1}}$. The convex hull of any subset ${w_1,\cdots,w_{r+1}} \subset {v_1, \cdots,v_{j+1}}$ is an $r$-face.}

Then, for any map $G: \R^{N+1} \mapsto \R^N$, the \emph{piecewise linear approximation} $G_\T$ to $G$ relative to the triangulation $\T$ of $\R^{N+1}$ is the unique map defined by:
\begin{itemize}
\item[(1)] $G_\T(v) = G(v)$ for all vertices of $\T$;
\item[(2)] for any $N+1$-simplex $\sigma = [v_1,v_2,\cdots,v_{N+2}] \in \T$, the restriction $G_\T |_\sigma$ of $G_\T$ to $\sigma$ is an affine map.
\end{itemize}
Consequently a point $Z = \sum_{i=1}^{N+2} \alpha_i v_i$ (here $\alpha_i$ are barycentric coordinates that satisfy $\sum_{i=1}^{N+2} \alpha_i = 1$ and $\alpha_i \geq 0$) in a $N+1$-simplex satisfies
$$
G_\T (Z) = G\left(\sum_{i=1}^{N+2} \alpha_i v_i\right) = \sum_{i=1}^{N+2} \alpha_i G(v_i).
$$
The set $G_\T^{-1} (0)$ contains a polygonal path $c_\T : \R \mapsto \R^{N +1}$ which approximates the path $c$.
Tracking such a path is carried out via PL-steps similar to the steps used in linear programming methods such as the Simplex Method.
Figure \ref{PLmethod} portrays the basic idea of a PL method.
\begin{figure}[h]
\centering
	\includegraphics[scale=0.9]{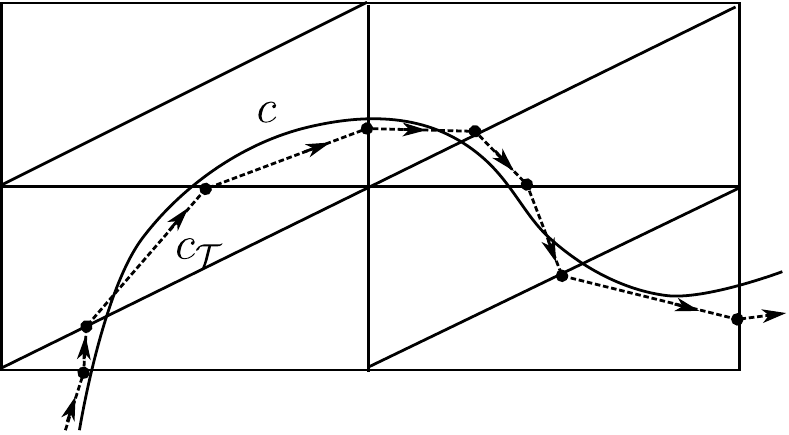}
	\caption{PL continuation.}
	\label{PLmethod}
\end{figure}

In aerospace applications, where the continuation procedure is in general differentiable, the PL methods are usually not as efficient as the PC methods or the differential continuation that we present in next sections.
Nevertheless when singularities exist in the zero path, the PL method is probably the most efficient one.

\section{Application to Attitude-Trajectory Optimal Control}
\label{fullexample}
In this section, the nonacademic attitude-trajectory optimal control problem for a launch vehicle (classical and airborne) is analyzed in detail.
Through this example, we illustrate how to analyze the (singular and regular) extremals of the problem with Lie and Poisson brackets, and how to elaborate numerical continuation procedures adapted to the solution structure. Indeed the theoretical analysis reveals the existence of a chattering phenomenon. Being aware of this feature is essential to devise an efficient numerical solution method.

\subsection{Geometric Analysis and Numerical Continuations for Optimal Attitude and Trajectory Control Problem ($\P_S$)}
\label{example1}
The problem is formulated in terms of dynamics, control, constraints and cost. The Pontryagin Maximum Principle and the geometric optimal control are then applied to analyze the extremals, revealing the existence of the chattering phenomenon.

\subsubsection{Formulation of ($\P_S$) and Difficulties}

\paragraph{Minimum time attitude-trajectory  control problem ($\P_S$).}
In this section, we formulate an attitude-trajectory minimum time control problem, denoted by ($\P_S$).

The trajectory of a launch vehicle is controlled by the thrust which can only have limited deflection angles with the vehicle longitudinal axis. Controlling the thrust direction requires controlling the vehicle attitude. When the attitude dynamics is slow, or when the orientation maneuver is large, this induces a coupling between the attitude motion and the trajectory, as explained in Section \ref{sec_appliaero}.

When this coupling is not negligible the dynamics and the state must account simultaneously for the trajectory variables (considering the launch vehicle as a mass point) and the attitude variables (e.g., the Euler angles or the quaternion associated to the body frame). 

The objective is then to determine the deflection angle law driving the vehicle from given initial conditions to the desired final attitude and velocity, taking into account the attitude-trajectory coupling.

The typical duration of such reorientation maneuvers is small compared to the overall launch trajectory. We assume therefore that 
the gravity acceleration is constant and we do not account for the position evolution. The aerodynamical forces (lift and drag) are supposed negligible in the first approach, and they will be introduced later in the system modelling.
The dynamics equations in an inertial frame $(O,x,y,z)$ are
\begin{equation}\label{sys_simple}
\begin{split}
 \dot{v}_x&= a \sin \theta \cos \psi + g_x , \\
 \dot{v}_y&= - a \sin \psi + g_y,\\
 \dot{v}_z&= a \cos \theta \cos \psi + g_z,\\
 \dot{\theta}&=(\omega_x \sin \phi + \omega_y \cos \phi)/ \cos \psi , \\
 \dot{\psi}&=\omega_x \cos \phi - \omega_y \sin \phi ,\\
 \dot{\phi}&= (\omega_x \sin \phi + \omega_y \cos \phi) \tan \psi ,\\
 \dot{\omega }_x&= - b u_2,\\
 \dot{\omega }_y&=  b u_1,
\end{split}
\end{equation}
where
($v_x$, $v_y$, $v_z$) represents the velocity,
($g_x$, $g_y$, $g_z$) represents the gravity acceleration,
$\theta$ (pitch), $\psi$ (yaw), $\phi$ (roll) are the Euler angles,
$a$ is the ratio of the thrust force to the mass,
and $b$ is the ratio of the thrust torque to the transverse inertia of the launcher ($a$ and $b$ are assumed constant).
$u=(u_1,u_2) \in \R^2$ is the control input of the system satisfying $|u|=u_1^2+u_2^2 \leq 1$.
See more details of the model and the problem formulation in \cite{ztc} or \cite{ztc3}.

Defining the state vector as $x=(v_x, v_y, v_z, \theta, \psi, \phi, \omega_x, \omega_y)$, we write the system \eqref{sys_simple} as the \emph{bi-input control-affine system}
\begin{equation} \label{sys_multi_smtc}
    \dot{x} = f(x) + u_1 g_1 (x)  + u_2 g_2 (x) ,
\end{equation}
where the controls $u_1$ and $u_2$ satisfy the constraint $u_1^2+u_2^2 \leq 1$,
and the vector fields $f$, $g_1$ and $g_2$ are defined by

\begin{multline}\label{fetgi_smtc}
f=  (a \sin \theta \cos \psi + g_x) \frac{\partial}{\partial v_x}
     + (- a \sin \psi + g_y) \frac{\partial}{\partial v_y}
     + (a \cos \theta \cos \psi + g_z) \frac{\partial}{\partial v_z}\\
    + (\omega_x \sin \phi + \omega_y \cos \phi)/ \cos \psi \frac{\partial}{\partial \theta}
     + (\omega_x \cos \phi - \omega_y \sin \phi)\frac{\partial}{\partial \psi} \\
    + \tan \psi (\omega_x \sin \phi + \omega_y \cos \phi) \frac{\partial}{\partial \phi} ,\quad
   g_1  = b \frac{\partial}{\partial \omega_y},\quad g_2 = - b \frac{\partial}{\partial \omega_x}.
\end{multline}

We define the target set (submanifold of $\R^8$)
\begin{multline} \label{target_smtc_PS}
M_{1S} =  \{  (v_x,v_y,v_z,\theta,\psi,\phi,\omega_x,\omega_y) \in \R^8 \ \mid\  v_z \sin \psi_f + v_y \cos \theta_f \cos \psi_f =0,\\
  v_z \sin \ theta_f - v_x \cos \theta_f \ =0, \quad \theta =\theta_f,\quad \psi=\psi_f,\quad \phi=\phi_f,\\
  \omega_x=\omega_{x_f},\quad \omega_y=\omega_{y_f} \} .
\end{multline}

The first two conditions in \eqref{target_smtc_PS} define a final velocity direction parallel to the longitudinal axis of the launcher, or in other terms a zero angle of attack.

The problem ($\P_S$) consists in steering the bi-input control-affine system \eqref{sys_multi_smtc} from $x(0)=x_0 = ({v_{x_0}},{v_{y_0}},{v_{z_0}},\theta_0,\psi_0,\phi_0,{\omega_{x_0}},{\omega_{y_0}}) \in \R^8$ to the final target $M_{1S}$ in minimum time $t_f$, with controls satisfying the constraint $u_1^2+u_2^2 \leq 1$.
The fixed initial condition is $x(0)=x_0$ and the final condition of problem $\P_S$ is
\begin{equation} \label{fc_smtc}
x(t_f) \in M_{1S},
\end{equation}
The initial and final conditions are also called terminal conditions.

\paragraph{Difficulties.}
The problem ($\P_S$) is difficult to solve directly due to the coupling of the attitude and the trajectory.
The system is of dimension $8$ and its dynamics contains both slow (trajectory) and fast (attitude) components.
This observation is be particularly important in order to design an appropriate solution method.
The idea is to define a simplified starting problem and then to apply continuation techniques.
However the essential difficulty of this problem is the chattering phenomenon making the control switch an infinite number of times over a compact time interval.
Such a phenomenon typically occurs when trying to connect bang arcs with higher-order singular arcs (see, e.g., \cite{Fuller,Marchal,Zelikin1,Zelikin2}, or Section \ref{chatphenom}).

In a preliminary step, we limited ourselves to the planar problem, which is a single-input control affine system. This planar problem is close to real flight conditions of a launcher ascent phase.
We have used the results of M.I. Zelikin and V.F. Borisov \cite{Zelikin1,Zelikin2} to understand the chattering phenomenon and to prove the local optimality of the chattering extremals. We refer the readers to \cite{ztc2} for details.

In a second step using the Pontryagin Maximum Principle and the geometric optimal control theory (see \cite{Agrachev,Schattler,Trelatsurvey}), we have established an existence result of the chattering phenomenon for a class of bi-input control affine systems and we have applied the result to the problem ($\P_s$).
More precisely, based on Goh and generalized Legendre-Clebsch conditions, we have proved that there exist optimal chattering arcs when connecting the regular arcs with a singular arc of order two.

\subsubsection{Geometric Analysis for ($\P_S$)}
\paragraph{Singular Arcs and Necessary Conditions for Optimality}\label{sec32}
The first step to analyze the problem is to apply the PMP (see Theorem \ref{thm_pmp}).
Let us consider the system \eqref{sys_multi_smtc}, with the vector fields $f$, $g_1$ and $g_2$ defined by \eqref{fetgi_smtc}.
According to the PMP, there must exist an absolutely continuous mapping $p(\cdot)=(p_{v_x}(\cdot), p_{v_y}(\cdot), p_{v_z}(\cdot), p_{\theta}(\cdot), p_{\psi}(\cdot), p_{\phi}(\cdot), p_{\omega_x}(\cdot), p_{\omega_y}(\cdot))$ defined on $[0,t_f]$, such that $p(t)\in T^*_{x(t)}M$ (cotangent space) for every $t\in[0,t_f]$, and a real number $p^0 \leq 0$, with $(p(\cdot),p^0)\neq 0$, such that
$\dot{x}(t) = \frac{\partial H}{\partial p}(x(t),p(t),p^0,u(t))$
and almost everywhere on $[0,t_f]$
\begin{equation*} 
\begin{split}
& \dot{p}_{v_x}= 0,\quad \dot{p}_{v_y}= 0,\quad \dot{p}_{v_z}= 0,\\
& \dot{p}_{\theta}= -a \cos \psi (p_{v_x} \cos \theta -p_{v_z} \sin \theta),\\
& \dot{p}_{\psi}= a \sin \psi \sin \theta p_{v_x} + a \cos \psi p_{v_y} + a \cos \theta \sin \psi p_{v_z} \\
&\qquad -\sin \psi (\omega_x \sin \phi + \omega_y \cos \phi)/\cos^2 \psi p_{\theta}
  -(\omega_x \sin \phi + \omega_y \cos \phi)/ \cos^2 \psi p_{\phi}, \\
& \dot{p}_{\phi} = -(\omega_x \cos \phi -\omega_y \sin \phi)/\cos \psi p_{\theta}
+(\omega_x \sin \phi + \omega_y \cos \phi) p_{\psi} \\
&\qquad\qquad -\tan \psi (\omega_x \cos \phi -\omega_y \sin \phi) p_{\phi}, \\
& \dot{p}_{\omega_x}= - \sin \phi / \cos \psi p_{\theta}- \cos \phi p_{\psi}
 -\sin \psi \sin \phi / \cos \psi p_{\phi}, \\
& \dot{p}_{\omega_y}=  - \cos \phi / \cos \psi p_{\theta}+ \sin \phi p_{\psi}
 -\sin \psi \cos \phi / \cos \psi p_{\phi} ,
\end{split}
\end{equation*}
The Hamiltonian of the optimal control problem ($\P_S$) is defined by
$$
H(x,p,p^0,u) = h_0(x,p)+u_1 h_1(x,p) +u_2 h_2(x,p)+p^0,
$$
with $h_0(x,p)=\langle p, f(x) \rangle $, $h_1(x,p)=\langle p, g_1(x) \rangle $, and $h_2(x,p)=\langle p, g_2(x) \rangle $. 
With a slight abuse of notation as before, we will denote $h_i(t)=h_i(x(t),p(t))$, $i=0,1,2$.

The maximization condition of the PMP yields, almost everywhere on $[0,t_f]$,
\begin{equation*} 
(u_1(t),u_2(t)) = \frac{(h_1(t),h_2(t)) }{\sqrt{h_1(t)^2+h_2(t)^2}} = \frac{\Phi(t)}{\Vert\Phi(t)\Vert},
\end{equation*}
whenever $\Phi(t)=(h_1(t),h_2(t)) = ( b p_{\omega_y}(t),- b p_{\omega_x}(t))\neq (0,0)$. The function $\Phi$ is of class $C^1$ and is called (as well as its components) the switching function. 
The switching manifold $\Gamma$ is the submanifold of $\R^{16}$ of codimension two defined by
$\Gamma = \left\{ z=(x,p)\in\R^{16} \mid p_{\omega_x} = p_{\omega_y} = 0 \right\}$.

The transversality condition $p(t_f) \perp T_{x(t_f)} M_1$ yields
\begin{equation}\label{trc}
    p_{v_x}(t_f) \sin \theta_f \cos \psi_f - p_{v_y}(t_f) \sin \psi_f + p_{v_z}(t_f) \cos \theta_f \cos \psi_f =0 .
\end{equation}
where $T_{x(t_f)}M_1$ is the tangent space to $M_1$ at the point $x(t_f)$. The final time $t_f$ being free and the system being autonomous, we have also
$h_0(x(t),p(t))+\Vert\Phi(t)\Vert+p^0=0,\: \forall t\in[0,t_f]$.

We say that an arc (restriction of an extremal to a subinterval $I$) is \emph{regular} if $\Vert \Phi(t)\Vert \neq 0$ along $I$. Otherwise, the arc is said to be \emph{singular}.
An arc that is a concatenation of an infinite number of regular arcs is said to be \emph{chattering}. The chattering arc is associated with a \emph{chattering control} that switches an infinite number of times, over a compact time interval.
A junction between a regular arc and a singular arc is said to be a \emph{singular junction}.

We next compute the singular control, since it is important to understand and explain the occurrence of chattering.
The usual method for to computing singular controls is to differentiate repeatedly the switching function until the control explicitly appears.
Note that here we need to use the notion of Lie bracket and Poisson bracket (see Section \ref{Lie}).

Assuming that $\Vert \Phi(t)\Vert = 0$ for every $t\in I$, i.e., $h_1(t)=h_2(t)=0$, and differentiating with respect to $t$, we get, using the Poisson brackets, 
\begin{align*}
\dot{h}_1 = \left\{ h_0,h_1 \right\} + u_2 \left\{ h_2,h_1 \right\}  = 0 ,\\
\dot{h}_2 = \left\{ h_0,h_2 \right\} + u_1 \left\{ h_1,h_2 \right\}  = 0,
\end{align*}
along $I$.
If the singular arc is optimal and the associated singular control is not saturating, then the Goh condition (see \cite{Goh}, see also Theorem \ref{Gohcond}) $ \left\{ h_1,h_2 \right\}  = \langle p, [g_1,g_2](x)\rangle= 0$ must be satisfied along $I$. Therefore we get that 
\begin{align*}
\dot{h}_1 = \left\{ h_0,h_1 \right\}  = \langle p, [f,g_1](x)\rangle=0,\\
\dot{h}_2 = \left\{ h_0,h_2 \right\}  = \langle p, [f,g_2](x)\rangle=0,
\end{align*}
along $I$.

Since the vector fields $g_1$ and $g_2$ commute, i.e., $[g_1,g_2] = 0$, we get by differentiating again that
\begin{align*}
\ddot{h}_1 = \left\{ h_0,\left\{h_0,h_1 \right\}\right\} + u_1 \left\{ h_1,\left\{h_0,h_1 \right\}\right\} + u_2 \left\{ h_2,\left\{h_0,h_1 \right\}\right\} = 0, \\
\ddot{h}_2 = \left\{ h_0,\left\{h_0,h_2 \right\}\right\}+ u_1 \left\{ h_1,\left\{h_0,h_2 \right\}\right\} + u_2 \left\{ h_2,\left\{h_0,h_2 \right\}\right\} = 0.
\end{align*}
Assuming that
\begin{equation*}
\det \Delta_1 = \det
    \begin{pmatrix}
       \left\{ h_1,\left\{h_0,h_1 \right\}\right\}    & \left\{ h_2,\left\{h_0,h_1 \right\}\right\}   \\
       \left\{ h_1,\left\{h_0,h_2 \right\}\right\}    & \left\{ h_2,\left\{h_0,h_2 \right\}\right\}
    \end{pmatrix}
\neq 0
\end{equation*}
along $I$, we obtain that
\begin{equation*} 
\begin{split}
u_1 &= \big( -\left\{ h_0,\left\{h_0,h_1 \right\}\right\}  \left\{ h_2,\left\{h_0,h_2 \right\}\right\} + \left\{ h_0,\left\{h_0,h_2 \right\}\right\} \left\{ h_2,\left\{h_0,h_1 \right\}\right\} \big) / \det \Delta_1, \\
u_2 &= \big(   \left\{ h_0,\left\{h_0,h_1 \right\}\right\} \left\{ h_1,\left\{h_0,h_2 \right\}\right\} - \left\{ h_0,\left\{h_0,h_2 \right\}\right\} \left\{ h_1,\left\{h_0,h_1 \right\}\right\} \big) / \det \Delta_1,
\end{split}
\end{equation*}
so that the control $u=(u_1,u_2)$ is said of order $1$.
$u_1$ and $u_2$ must moreover satisfy the constraint $u_1^2+u_2^2 \leq 1$.

However, in problem ($\P_S$), we have $[g_1,[f,g_2]]=0$, $[g_2,[f,g_1]]=0$, and $\left\{ h_1,\left\{h_0,h_1 \right\}\right\} =\left\{ h_2,\left\{h_0,h_2 \right\}\right\}=0$ along $I$, which indicates that the singular control is of higher order.
According to the Goh condition (see \cite{Goh, Krener}, see also Definition \ref{def_singorder} and Theorem \ref{Gohcond}), we must have $\left\{ h_i,\left\{h_0,h_j \right\}\right\} = 0$, $i,j=1,2$, $i\neq j$, and we can go on differentiating. It follows from $[ g_1,[f,g_1 ]] = 0$ and $[ g_2,[f,g_2 ]] = 0$ that
$$
[g_i, \mathrm{ad}^2 f.g_i]] = [g_i, [f, \mathrm{ad}\, f.g_i]] = -[f,[\mathrm{ad}\, f.g_i,g_i]] - [\mathrm{ad}\, f.g_i, [g_i,f]] = 0, \quad i=1,2 ,
$$
and we get
\begin{equation*} 
\begin{split}
& h_1^{(3)} = \left\{ h_0,\mathrm{ad}^2 h_0. h_1 \right\} + u_2 \left\{ h_2,\mathrm{ad}^2 h_0. h_1 \right\}  = 0,\\
& h_2^{(3)} = \left\{ h_0,\mathrm{ad}^2 h_0. h_2 \right\} + u_1 \left\{ h_1,\mathrm{ad}^2 h_0. h_2 \right\}  = 0.
\end{split}
\end{equation*}
Using $[g_1,g_2] = 0$ and $[ g_i,[f,g_i ]] = 0$, $i=1,2$, it follows that $[g_k,[g_i,[f,g_j ]]=0$, $i,j,k=1,2$ and
$$
\frac{d}{dt} \left\{ h_2,\left\{h_0,h_1 \right\}\right\}=\frac{d}{dt} \left\{ h_1,\left\{h_0,h_2 \right\}\right\}
=\left\{ h_0, \left\{ h_1,\left\{h_0,h_2 \right\}\right\} \right\}=0.
$$
This is a new constraint along the singular arc. The time derivative of this constraint is equal to zero and therefore does not induce any additional constraint.

The higher-order necessary conditions for optimality (see Definition \ref{def_singorder}) state that an optimal singular control can only appear explicitly within an even derivative. Therefore we must have 
$$\left\{ h_2,\mathrm{ad}^2 h_0. h_1 \right\}=\left\{ h_1,\mathrm{ad}^2 h_0. h_2 \right\} = 0$$
 along $I$. Accordingly, $h_i^{(3)}=0$, $i=1,2$, gives three additional constraints along the singular arc:
$$\left\{ h_0,\mathrm{ad}^2 h_0. h_1 \right\}=\left\{ h_0,\mathrm{ad}^2 h_0. h_2 \right\}=\left\{ h_2,\mathrm{ad}^2 h_0. h_1 \right\}  = \left\{ h_1,\mathrm{ad}^2 h_0. h_2 \right\}=0.
$$
By differentiating the first two constraints with respect to $t$, we get
\begin{align*}
h_1^{(4)}=\mathrm{ad}^4h_0.h_1 + u_1 \left\{ h_1,\mathrm{ad}^3h_0.h_1\right\} + u_2 \left\{ \mathrm{ad}^2h_0.h_1,\mathrm{ad}\, h_0.h_2\right\} = 0, \\
h_2^{(4)}=\mathrm{ad}^4h_0.h_2 + u_1 \left\{ \mathrm{ad}^2h_0.h_2,\mathrm{ad}\, h_0.h_1\right\} + u_2 \left\{ h_2,\mathrm{ad}^3h_0.h_2\right\}= 0 .
\end{align*}
Assuming that $\left\{h_i,\mathrm{ad}^3h_0.h_i\right\}<0$ for $i=1,2$ (generalized Legendre-Clebsch condition, see Corollary \ref{necconds}) and since 
$$\left\{ \mathrm{ad}^2h_0.h_1,\mathrm{ad}\, h_0.h_2\right\} =\left\{ \mathrm{ad}^2h_0.h_2,\mathrm{ad}\, h_0.h_1\right\} = \mathrm{ad}^4h_0.h_1=\mathrm{ad}^4h_0.h_2=0$$
along $I$ for problem ($\P_S$), the singular control is 
\begin{equation*}
u_1  =0, \quad
u_2  =0.
\end{equation*}
The singular control $u=(u_1,u_2)$ is then said of \emph{intrinsic order two} (see the precise definition in Definition \ref{def_singorder}).

Let us assume that $(x(\cdot),p(\cdot),p^0,u(\cdot))$ is a singular arc of ($\P_S$) along the subinterval $I$, which is locally optimal in $C^0$ topology.
Then we have $u=(u_1,u_2)=(0,0)$ along $I$, and $u$ is a singular control of intrinsic order two.
Moreover, we can establish (see the proof in \cite{ztc}) that this singular extremal must be normal, i.e., $p^0 \neq 0$, and according to  Lemma \ref{necconds}, the Generalized Legendre-Clebsch Condition (GLCC) along $I$ takes the form
\begin{equation}\label{lille1532}
a + g_x\sin\theta\cos\psi-g_y\sin\psi+g_z\cos\theta\cos\psi \geq 0,
\end{equation}

We define next the singular surface $S$, which is filled by singular extremals of ($\P_S$), by
\begin{multline} \label{ssurface}
S = \Big\{(x,p)\ \mid\ \omega_x=\omega_y=0, \quad p_{\theta}=p_{\psi}=p_{\phi}=p_{\omega_x}=p_{\omega_y}=0,\quad p_{v_x} = \tan \theta p_{v_z}, \\
 p_{v_z}= \frac{-p^0 \cos\theta\cos\psi}{a+g_x\sin \theta \cos \psi-g_y\sin \psi+g_z\cos \theta \cos \psi}, \quad p_{v_y}=-\tan \psi / \cos \theta p_{v_z} \Big\}.
\end{multline}
We will see later that the solutions of the problem of order zero (defined in the following Section) lie on this singular surface $S$.

Finally, the possibility of chattering in problem ($\P_S$) is demonstrated in \cite{ztc}. A chattering arc appears when trying to connect a regular arc with an optimal singular arc. More precisely, let $u$ be an optimal control, solution of ($\P_S$), and assume that $u$ is singular on the sub-interval $(t_1,t_2)\subset[0,t_f]$ and is regular elsewhere. If $t_1>0$ (resp., if $t_2<t_f$) then, for every $\varepsilon>0$, the control $u$ switches an infinite number of times over the time interval $[t_1-\varepsilon,t_1]$ (resp., on $[t_2,t_2+\varepsilon]$).
The condition \eqref{lille1532} was required in the proof.

The knowledge of chattering occurrence is essential for solving the problem ($\P_S$) in practice. Chattering raises indeed numerical issues that may prevent any convergence, especially when using an indirect approach (shooting). The occurrence of the chattering phenomenon in ($\P_S$) explains the failure of the indirect methods for certain terminal conditions (see also the recent paper \cite{Caponigro}).

\subsubsection{Indirect Method and Numerical Continuation Procedure for ($\P_S$)}
The principle of the continuation procedure is to start from the known solution of a simpler problem (called hereafter the \emph{problem of order zero}) in order to initialize an indirect method for the more complicated problem ($\P_S$).
This simple low-dimensional problem will then be embedded in higher dimension, and appropriate continuations will be applied to come back to the initial problem.

The problem of order zero defined below considers only the trajectory dynamics which is much slower than the attitude dynamics. Assuming an instantaneous attitude motion simplifies greatly the problem and provides an analytical solution.
It is worth noting that the solution of the problem of order zero is contained in the singular surface $S$ filled by the singular solutions for ($\P_S$), defined by \eqref{ssurface}.

\paragraph{Auxiliary Problems.}
We define the \emph{problem of order zero}, denoted by ($\P_0$) as the ``subproblem'' of problem ($\P_S$) reduced to the trajectory dynamics. The control for this problem is directly the vehicle attitude, and the attitude dynamics is not simulated.

Denoting the vehicle longitudinal axis as $\vec{e}$ and considering it as the control vector (instead of the attitude angles $\theta$, $\psi$), we formulate the problem as follows:
\begin{align*}
& \dot{\vec{v}} = a \vec{e}+\vec{g}, \\
& \vec{v}(0)=\vec{v}_0,\quad \vec{v}(t_f) // \vec{w},\\
& \Vert \vec{w} \Vert =1,  \\
& \min t_f ,
\end{align*}
where $\vec{w}$ is a given vector that refers to the desired target velocity direction, and $\vec{g}$ is the gravitational acceleration vector.
The solution of this problem is straightforward and gives :
the optimal solution of ($\P_0$) is given by
\begin{equation*}
 \vec{e}^{\ast}=\frac{1}{a}\left(\frac{k \vec{w}-\vec{v}_0}{t_f}-\vec{g}\right), \quad
 t_f=\frac{-a_2 + \sqrt{a_2^2-4a_1a_3}}{2 a_1},\quad
 \vec{p}_v=\frac{-p^0}{a + \langle \vec{e}^{\ast},\vec{g} \rangle } \vec{e}^{\ast}.
\end{equation*}
with
$$k=\langle \vec{v}_0,\vec{w} \rangle + \langle \vec{g},\vec{w} \rangle t_f,$$
$$a_1 = a^2-\Vert  \langle \vec{g},\vec{w} \rangle \vec{w}-\vec{g} \Vert ^2,$$
$$a_2 = 2 ( \langle \vec{v}_0,\vec{w} \rangle \langle \vec{g},\vec{w} \rangle - \langle \vec{v}_0,\vec{g} \rangle),$$
and
$$a_3 = - \Vert  \langle \vec{v}_0,\vec{w} \rangle \vec{w}-\vec{v}_0 \Vert ^2.$$
We refer the readers to \cite{ztc} for the detailed calculation.

The Euler angles $\theta^{\ast} \in (-\pi,\pi)$ and $\psi^{\ast} \in (-\pi/2,\pi/2)$ are retrieved from the components of the vector $\vec{e}^\ast$ since
$
\vec{e}^\ast = (\sin \theta^\ast \cos \psi^\ast, -\sin \psi^\ast, \cos \theta^\ast \sin \psi^\ast)^{\top}
$.

We can check that these optimal angles $\theta = \theta^\ast$, $\psi=\psi^\ast$ and $\phi=\phi^\ast$ (whatever the value of $\phi^\ast$) satisfy the equations \eqref{ssurface}, so that the solution of ($\P_0$) is contained in the singular surface $S$. 
The optimal solution of ($\P_0$) actually corresponds to a singular solution of ($\P_S$) with the terminal conditions given by
\begin{equation} \label{inicond_ocpztomtcp}
\begin{split}
v_x(0) = {v_{x_0}}, \quad v_y(0) = {v_{y_0}}, \quad v_z(0) = {v_{z_0}}, \\
\theta(0)=\theta^\ast, \quad \psi(0)=\psi^\ast,\quad, \phi(0)=\phi^\ast,\quad \omega_x(0) = 0,\quad \omega_y(0)=0,
\end{split}
\end{equation}
\begin{equation} \label{fincond_ocpztomtcp}
v_z(t_f)\sin \psi_f + v_y(t_f) \cos \theta_f \cos \psi_f=0,\quad v_z(t_f)\sin \theta_f-v_x(t_f)\cos \theta_f=0,
\end{equation}
\begin{equation} \label{fincond_ocpztomtcp2}
 \theta(t_f)=\theta^\ast, \quad \psi(t_f)=\psi^\ast,\quad, \phi(t_f)=\phi^\ast,\quad \omega_x(t_f) = 0,\quad \omega_y(t_f)=0.
 \end{equation}

A natural continuation strategy consists in changing continuously  these terminal conditions \eqref{inicond_ocpztomtcp}-\eqref{fincond_ocpztomtcp2} to come back to the terminal conditions \eqref{fc_smtc} of ($\P_S$).

Unfortunately the chattering phenomenon may prevent the convergence of the shooting method. When the terminal conditions are in the neighborhood of the singular surface $S$, the optimal extremals are likely to contain a singular arc and thus chattering arcs causing the failure of the shooting method. In order to overcome the numerical issues we define a regularized problem with a modified cost functional.

\medskip

The regularized problem ($\P_R$) consists in minimizing the cost functional
\begin{equation} \label{cost_opcr}
    C_K = t_f + K \int_0^{t_f} (u_1^2 + u_2^2) \, dt ,
\end{equation}
for the bi-input control-affine system \eqref{sys_multi_smtc}, under the control constraints $-1\leq u_i\leq 1$, $i=1,2$, and with the terminal conditions \eqref{fc_smtc}. The constant $K > 0$ is arbitrary.
We have replaced the constraint $u_1^2+u_2^2\leq 1$ (i.e., $u$ takes its values in the unit Euclidean disk) with the constraint that $u$ takes its values in the unit Euclidean square. Note that we use the Euclidean square (and not the disk) because we observed that our numerical simulations worked better in this case.
This regularized optimal control problem with the cost \eqref{cost_opcr} has continuous extremal controls and it is therefore well suited to a continuation procedure.

The Hamiltonian of problem ($\P_R$) is
\begin{equation}\label{Hamil_n}
H_K=\langle p,f(x) \rangle + u_1\langle p,g_1(x) \rangle + u_2 \langle p,g_2(x) \rangle +p^0(1+K u_1^2 + K u_2^2),
\end{equation}
and according to the PMP, the optimal controls are
\begin{equation} \label{OCPn_u1}
\begin{split}
u_1(t) &= \mathrm{sat}\left(-1,- \frac{\bar{b} p_{\omega_y}(t) }{ 2 K p^0},1\right),\\\
u_2(t) &= \mathrm{sat}\left(-1, \frac{\bar{b} p_{\omega_x}(t) }{ 2 K p^0},1\right),
\end{split}
\end{equation}
where the saturation operator $\mathrm{sat}$ is defined by 
$$
\mathrm{sat}(-1,f(t),1) =
\begin{cases}
-1 &\textrm{if}\quad f(t)\leq -1, \\
1&\textrm{if}\quad f(t)\geq 1,\\
f(t) & \textrm{if}\quad -1\leq f(t)\leq 1.
\end{cases}
$$
An important advantage of considering problem ($\P_R$) is that when we embed the solutions of ($\P_0$) into the ($\P_R$), they are not singular, whereas the solution of ($\P_0$) is a \emph{singular} trajectory of the full problem ($\P_S$) and thus passing directly from ($\P_0$) to ($\P_S$) causes essential difficulties due to chattering.
More precisely,
an extremal of ($\P_0$) can be embedded into ($\P_R$), by setting
$$
u(t) = (0,0), \quad \theta(t)=\theta^\ast, \quad \psi(t)=\psi^\ast, \quad \phi(t)=\phi^\ast, \quad \omega_x(t) = 0, \quad \omega_y(t)=0,
$$
$$
p_{\theta}(t)=0, \quad p_{\psi}(t)= 0, \quad p_{\phi}(t)=0, \quad p_{\omega x}(t) = 0, \quad p_{\omega y}(t)=0,
$$
where $\theta^\ast$ and $\psi^\ast$ are given by solving problem $\P_0$,
with the natural terminal conditions given by \eqref{inicond_ocpztomtcp} and \eqref{fincond_ocpztomtcp}-\eqref{fincond_ocpztomtcp2}. This solution is not a singular extremal for ($\P_R$).
The extremal equations for ($\P_R$), are the same than for ($\P_S$), as well as the transversality conditions.

\paragraph{Numerical Continuation Procedure.}
The objective is to find the optimal solution of ($\P_S$), starting from the explicit solution of $\P_0$. We proceed as follows:
\begin{itemize}
\item First, we embed the solution of ($\P_0$) into ($\P_R$).
For convenience, we still denote ($\P_0$) the problem ($\P_0$) formulated in higher dimension.
\item Then, we pass from ($\P_0$) to ($\P_S$) by means of a numerical continuation procedure, involving three continuation parameters. The first two parameters $\lambda_1$ and $\lambda_2$ are used to pass continuously from the optimal solution of ($\P_0$) to the optimal solution of the regularized problem ($\P_R$) with prescribed terminal attitude conditions, for some fixed $K>0$. The third parameter $\lambda_3$ is then used to pass to the optimal solution of ($\P_S$) (see Figure \ref{homotopie}).
\end{itemize}

\begin{figure}[h]
\centering
	\includegraphics[scale=1.0]{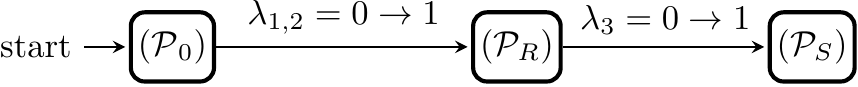}
	\caption{Continuation procedure for ($\P_S$).}
	\label{homotopie}
\end{figure}

In a first step, we use the continuation parameter $\lambda_1$ to act on the initial conditions, according to
$$
\theta(0) = \theta^{\ast} (1-\lambda_1) + \theta_0 \lambda_1, \quad \psi(0) = \psi^{\ast} (1-\lambda_1) + \psi_0 \lambda_1, \quad \phi(0) = \phi^{\ast} (1-\lambda_1) + \phi_0 \lambda_1,
$$
$$
\omega_x(0) = \omega_x^{\ast} (1-\lambda_1) + {\omega_{x_0}} \lambda_1, \quad \omega_y(0) = \omega_y^{\ast} (1-\lambda_1) + {\omega_{y_0}} \lambda_1,
$$
where $\omega_x^{\ast} = \omega_y^{\ast} =0$, $\phi^{\ast} = 0$, and $\theta^{\ast}$, $\psi^{\ast}$ are given by the explicit solution of the problem ($\P_0$).

Using the transversality condition \eqref{trc} and the extremal equations, the shooting function $ S_{\lambda_1} $ for the $\lambda_1$-continuation is of dimension 8 and defined by
\begin{equation*}
\begin{split}
S_{\lambda_1} = \big(
        &  p_{\omega_x }(t_f) , \,\,
            p_{\omega_y} (t_f),\,\,
            p_{\theta} (t_f),\,\,
            p_{\psi} (t_f),\,\,
            p_{\phi} (t_f),\,\,
           H_K(t_f), \\
        &  v_z(t_f) \sin \psi_f + v_y(t_f) \cos \theta_f \cos \psi_f,\,\,
         v_z(t_f) \sin \theta_f - v_x(t_f) \cos \theta_f \big),
\end{split}
\end{equation*}
where $H_K(t_f)$ with $p^0=-1$ is calculated from \eqref{Hamil_n} and $u_1$ and $u_2$ are given by \eqref{OCPn_u1}.
Recall that we have proved that a singular extremal of problem ($\P_S$) must be normal, and since we are starting to solve the problem from a singular extremal, we can assume that $p^0 = -1$.

Note again that there is no concern using $S_{\lambda_1}$ as shooting function for ($\P_R$). This would not be the case for ($\P_S$) : if $S_{\lambda_1}=0$, then together with $\omega_x(t_f)=0$ and $\omega_y(t_f)=0$, the final point $(x(t_f),p(t_f))$ of the extremal would lie on the singular surface $S$ defined by \eqref{ssurface} and this would cause the failure of the shooting method. On the opposite, for problem ($\P_R$), even when $x(t_f) \in S$, the shooting problem is smooth and it can still be solved.

The solution of ($\P_0$) is a solution of ($\P_R$) for $\lambda_1=0$, corresponding to the terminal conditions \eqref{inicond_ocpztomtcp}-\eqref{fincond_ocpztomtcp}  (the other states at $t_f$ being free). By continuation, we vary $\lambda_1$ from $0$ to $1$, yielding the solution of ($\P_R$), for $\lambda_1=1$. The final state of the corresponding extremal gives some unconstrained Euler angles denoted by
$\theta_e = \theta(t_f)$, $\psi_e = \psi(t_f)$, $\phi_e = \phi(t_f)$, $\omega_{xe} = \omega_x(t_f)$ and $\omega_{ye} = \omega_y(t_f)$.

\medskip

In a second step, we use the continuation parameter $\lambda_2$ to act on the final conditions, in order to make them pass from the values $\theta_e$, $\psi_e$, $\phi_e$, $\omega_{xe}$ and $\omega_{ye}$, to the desired target values $\theta_f$, $\psi_f$, $\phi_f$, $\omega_{xf}$ and $\omega_{yf}$. The shooting function $ S_{\lambda_2} $ for the $\lambda_2$-continuation is still of dimension $8$ and defined by
\begin{equation*}
 S_{\lambda_2} = \begin{pmatrix}
        \omega_x(t_f) - (1-\lambda_2) \omega_{xe} -\lambda_2 \omega_{x_f} \\
        \omega_y(t_f) - (1-\lambda_2) \omega_{ye} -\lambda_2 \omega_{y_f}\\
        \theta(t_f)-(1-\lambda_2) \theta_e - \lambda_2 \theta_f \\
        \psi(t_f)-(1-\lambda_2) \psi_e - \lambda_2 \psi_f \\
        \phi(t_f)-(1-\lambda_2) \phi_e - \lambda_2 \phi_f \\
        v_z(t_f) \sin \psi_f + v_y(t_f) \cos \theta_f \cos \psi_f \\
        v_z(t_f) \sin \theta_f - v_x(t_f) \cos \theta_f \\
        H_K(t_f)
        \end{pmatrix}.
\end{equation*}
Solving this problem by varying $\lambda_2$ from $0$ to $1$, we obtain the solution of ($\P_R$), with the terminal condition \eqref{fc_smtc}.

\medskip

Finally, in order to compute the solution of ($\P_S$), we use the continuation parameter $\lambda_3$ to pass from ($\P_R$) to ($\P_S$). We introduce the parameter $\lambda_3$ into the cost functional \eqref{cost_opcr} and the Hamiltonian $H_{K}$ as follows:
$$
C_K = t_f + K \int_0^{t_f} (u_1^2 + u_2^2)(1-\lambda_3) \, dt ,
$$
$$
H(t_f,\lambda_3)=\langle p,f \rangle+\langle p,g_1 \rangle u_1 +\langle p,g_2 \rangle u_2 +p^0 + p^0K( u_1^2 +u_2^2)(1-\lambda_3).
$$
According to the PMP, the extremal controls of this problem are given by $u_i=\mathrm{sat}(-1,u_{ie},1)$, $i=1,2$, where
\begin{align*}
u_{1e} &=   \frac{\bar{b} p_{\omega_y}}{-2 p^0 K (1-\lambda_3) + \bar{b}  \lambda_3 \sqrt{p_{\omega_x}^2+p_{\omega_y}^2} } ,\\
u_{2e} &= \frac{ -\bar{b} p_{\omega_x}}{-2 p^0 K (1-\lambda_3) + \bar{b}  \lambda_3 \sqrt{p_{\omega_x}^2+p_{\omega_y}^2} } .
\end{align*}
The shooting function $S_{\lambda_3}$ is defined similarly to $S_{\lambda_2}$, replacing $H_K(t_f)$ with $H_K(t_f,\lambda_3)$. The solution of ($\P_S$) is then obtained by varying $\lambda_3$ continuously from $0$ to $1$.

\medskip

This last continuation procedure fails in case of chattering, and thus it cannot be successful for any arbitrary terminal conditions. In particular, if chattering occurs then the $\lambda_3$-continuation is expected to fail for some value $\lambda_3 = \lambda_3^\ast<1$. In such a case this value of $\lambda_3$ corresponds to a sub-optimal solution of ($\P_S$), which is practically valuable since it satisfies the terminal conditions with a reduced final time (also not minimal), with a continuous control.
The numerical experiments show that this continuation procedure is very efficient. In most cases, optimal solutions with prescribed terminal conditions can be obtained within a few seconds (without parallel calculations).

\subsubsection{Direct Method}
In this section we envision a direct approach for solving ($\P_S$), with a piecewise constant control over a given time discretization.
The solutions obtained with such a method are sub-optimal, especially when the control is chattering (the number of switches being limited by the time step).

Since the initialization of a direct method may also raise some difficulties, we propose the following strategy. The idea is to start from the problem ($\P_S$) with relaxed terminal requirements, in order to get a first solution, and then to reintroduce step by step the final conditions \eqref{fc_smtc} of ($\P_S$). We implement this direct approach with the software \texttt{BOCOP} and its batch optimization option (see \cite{BonnansMartinon}).

\begin{itemize}
\item Step 1: we solve ($\P_S$) with the initial condition $x(0)=x_0$ and the final conditions
$$
\omega_{y}(t_f)=0,\quad \theta (t_f) = \theta_f ,\quad v_{z} (t_f) \sin \theta_f - v_{x} (t_f) \cos \theta_f =0.
$$
These final conditions are those of the planar version of ($\P_S$) (see \cite{ztc2} for details). This problem is easily solved by a direct method without any initialization care (a constant initial guess for the discretized variables suffices to ensure convergence).

\item Then, in Steps 2, 3, 4 and 5, we add successively (and step by step) the final conditions
$$v_{z} (t_f) \sin \psi_f + v_{y} (t_f) \cos \theta_f \cos \psi_f =0,$$
$$\psi (t_f)= \psi_f,\quad \phi (t_f)=\phi_f,\quad \omega_{x}(t_f) = \omega_{xf},$$
and for each new step we use the solution of the previous one as an initial guess.
\end{itemize}
At the end of this process, we have obtained the solution of ($\P_S$).

\subsubsection{Comparison of the Indirect and Direct Approaches}
\label{example2}
So far, in order to compute numerically the solutions of ($\P_S$), we have implemented two approaches. The indirect approach, combining shooting and numerical continuation, is time-efficient when the solution does not contain any singular arcs.

Depending on the terminal conditions, the optimal solution of ($\P_S$) may involve a singular arc of order two, and the connection with regular arcs generates chattering.
The occurrence of chattering causes the failure of the indirect approach. For such cases, we have proposed two alternatives.
The first alternative is based on an indirect approach involving three continuations. The last continuation starting from a regularized problem with smooth controls aims at coming back to the original problem that may be chattering. When chattering appears the continuation fails, but the last successful step provides a valuable smooth solution meeting the terminal conditions.

The second alternative is based on a direct approach, and it yields as well a sub-optimal solution having a finite number of switches. The number of switches is limited by the discretization step.
In any case, the direct strategy is much more time consuming than the indirect approach and the resulting control may exhibit many numerical oscillations as can be observed on Figure \ref{directBocop}. This kind of solutions is practically undesirable.

\begin{figure}[h]
\centering
\includegraphics[width=0.85\textwidth]{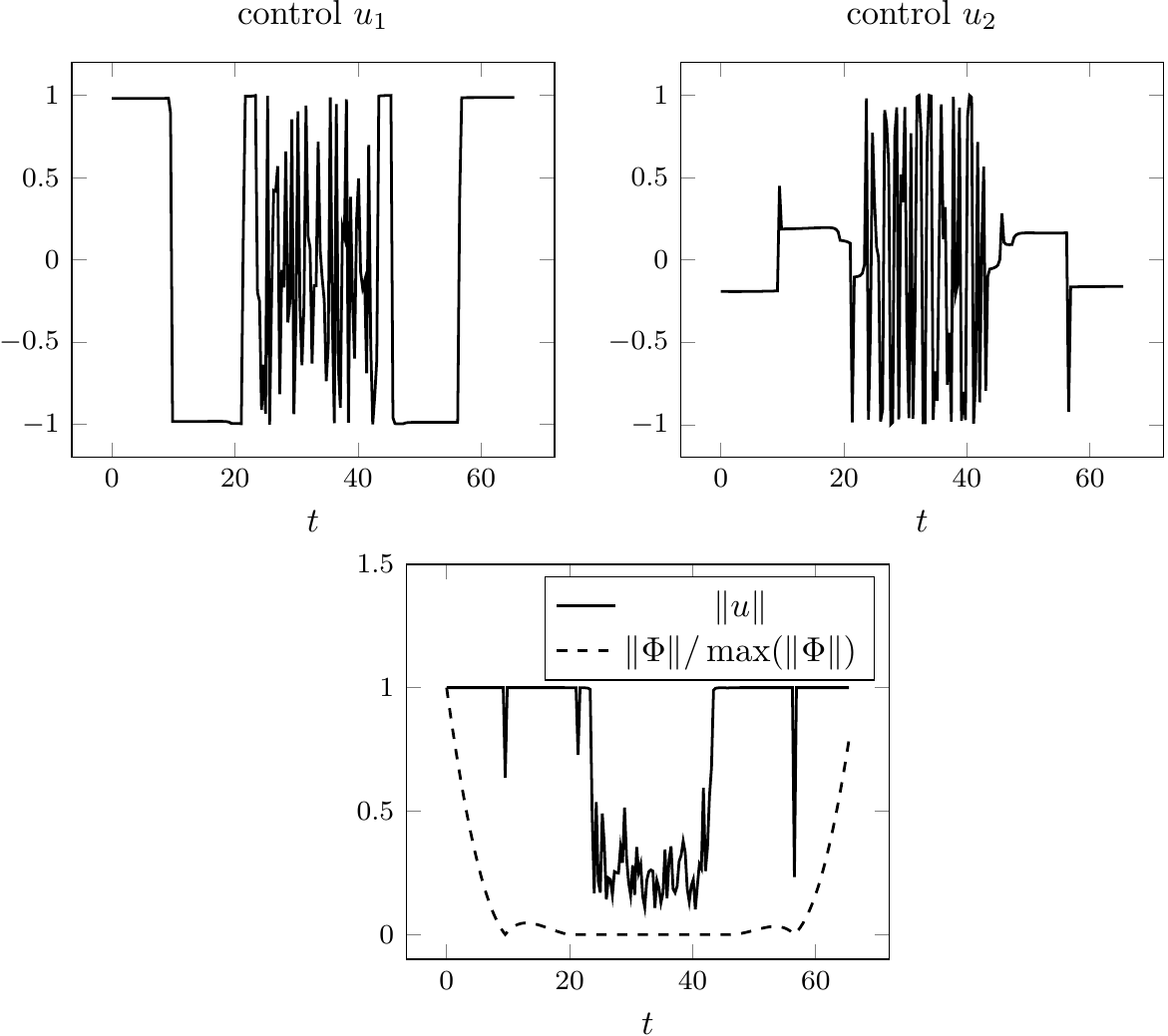}
\caption{Control $u(t)$ for ($\P_S$) obtained by BOCOP.}
\label{directBocop}
\end{figure}

Note that with both approaches, no a priori knowledge of the solution structure is required (in particular, the number of switches is unknown).

\medskip

As a conclusion about this example ($\P_S$), we can emphasize that the theoretical analysis has revealed the existence of singular solutions with possible chattering. This led us to introduce a regularized problem in order to overcome this essential difficulty.
On the other hand a continuation procedure is devised considering the dynamics slow-fast rates. This procedure is initiated by the problem of order zero reduced to the trajectory dynamics.

In the next section, we extend this approach to a more complicated problem (optimal pull-up maneuvers of airborne launch vehicles), in order to further illustrate the potential of continuation methods in aerospace applications.

\subsection{Extension to Optimal Pull-up Maneuver Problem ($\P_A$)}
\label{sec_ext_pb}
Since the first successful flight of Pegasus vehicle in April 1990, the airborne launch vehicles have always been a potentially interesting technique for small and medium-sized space transportation systems. The mobility and deployment of the airborne launch vehicles provide increased performance and reduced velocity requirements due to non-zero initial velocity and altitude.
Airborne launch vehicles consist of a carrier aircraft and a rocket-powered launch vehicle. The launch vehicle is released almost horizontally from the carrier aircraft and its engine is ignited a few seconds later once the carrier aircraft has moved away.
The flight begins with a pull-up maneuver \cite{sarigul3,sarigul2} targeting the optimal flight path angle for the subsequent ascent at zero angle of attack.
The kinematics conditions for the Pegasus vehicle are recalled here after \cite{Berend,Clegern,Mosier,Roble}.
The release takes place horizontally at an altitude of $12.65\,km$. The first stage is ignited at an altitude of $12.54 \,km$ and a velocity of $236.8\, m/s$ ($0.8$ Mach). The pull-up maneuver targets a flight path angle of $13.8^\circ$ at the end of the first stage flight. The load factor is limited to $2.5\, g$ and the dynamic pressure is limited to $47.6\, kPa$.

The pull-up maneuver consists in an attitude maneuver such that the flight path angle increases up to its targeted value, while satisfying the state constraints on the load factor and the dynamic pressure.
In this section, we address the minimum time-energy pull-up maneuver problem for airborne launch vehicles with a focus on the numerical solution method.

The model of the control system is more complex than \eqref{sys_simple} due to the aerodynamics forces that depend on the flight conditions (atmospheric density depending on the altitude, vehicle angle of attack):
\begin{equation}\label{sys_full_com}
\begin{split}
& \dot{r}_x = v_x, \quad \dot{r}_y = v_y,\quad \dot{r}_z = v_z,\\
& \dot{v}_x= a \sin \theta \cos \psi + g_x +(D_x+L_x)/m, \\
& \dot{v}_y= - a \sin \psi + g_y +(D_y+L_y)/m,\\
& \dot{v}_z= a \cos \theta \cos \psi + g_z +(D_z+L_z)/m,\\
& \dot{\theta}=(\omega_x \sin \phi + \omega_y \cos \phi)/ \cos \psi , \\
&  \dot{\psi}=\omega_x \cos \phi - \omega_y \sin \phi ,\\
& \dot{\phi}= (\omega_x \sin \phi + \omega_y \cos \phi) \tan \psi ,\\
& \dot{\omega }_x= - b u_2,\qquad \dot{\omega }_y=  b u_1.
\end{split}
\end{equation}
where ($r_x$, $r_y$, $r_z$) is the position,
$m$ is the mass,
($L_x$, $L_y$, $L_z$) is the lift force,
and ($D_x$, $D_y$, $D_z$) is the drag force.

Defining the state variable
$x=(r_x,r_y,r_z,v_x, v_y, v_z, \theta, \psi, \phi, \omega_x, \omega_y)$, we write the system \eqref{sys_full_com} as a bi-input control-affine system
\begin{equation} \label{sys_full}
    \dot{x} = \hat{f}(x) + u_1 \hat{g}_1 (x)  + u_2 \hat{g}_2 (x) ,
\end{equation}
where the controls $u_1$ and $u_2$ satisfy the constraint $u_1^2+u_2^2 \leq 1$,
and the smooth vector fields $\hat{f}$, $\hat{g}_1$ and $\hat{g}_2$ are defined by
\begin{equation*}
\begin{split}
\hat{f}= &v_x \frac{\partial}{\partial r_x} + v_y \frac{\partial}{\partial r_y} + v_z \frac{\partial}{\partial r_z}
   + (a \sin \theta \cos \psi + g_x+(D_x+L_x)/m) \frac{\partial}{\partial v_x}\\
   &+ (- a \sin \psi + g_y+(D_y+L_y)/m) \frac{\partial}{\partial v_y}
   + (a \cos \theta \cos \psi + g_z+(D_z+L_z)/m) \frac{\partial}{\partial v_z}\\
   &+ (\omega_x \sin \phi + \omega_y \cos \phi)/ \cos \psi \frac{\partial}{\partial \theta}
   + (\omega_x \cos \phi - \omega_y \sin \phi)\frac{\partial}{\partial \psi}\\
   &+ \tan \psi (\omega_x \sin \phi + \omega_y \cos \phi) \frac{\partial}{\partial \phi} ,\\
\hat{g}_1 = & b \frac{\partial}{\partial \omega_y},\quad
\hat{g}_2 =  - b \frac{\partial}{\partial \omega_x}.
\end{split}
\end{equation*}

The initial state is fixed
$
x_0 = (r_{x0},r_{y0},r_{z0},{v_{x_0}},{v_{y_0}},{v_{z_0}},\theta_0,\psi_0,\phi_0,{\omega_{x_0}},{\omega_{y_0}}) \in \R^{11}
$,
and the target set is defined by (submanifold of $\R^{11}$)
\begin{multline*} 
M_1 =  \Big\{  (r_x,r_y,r_z,v_x,v_y,v_z,\theta,\psi,\phi,\omega_x,\omega_y) \in \R^{11} \ \mid\  v_z \sin \psi_f + v_y \cos \theta_f \cos \psi_f =0,\\
 v_z \sin \theta_f - v_x \cos \theta_f =0, \quad \theta =\theta_f,\quad \psi=\psi_f,\quad \phi=\phi_f,\\
 \omega_x=\omega_{x_f},\quad \omega_y=\omega_{y_f} \Big\} .
\end{multline*}
The optimal pull-up maneuver problem ($\P_A$) consists in steering the bi-input control-affine system \eqref{sys_full} from
\begin{equation} \label{ic_cmtc}
x(0)=x_0
\end{equation}
to a point belonging to the final target $M_1$, i.e.,
\begin{equation} \label{fc_cmtc}
\begin{split}
x(t_f) \in M_1,
\end{split}
\end{equation}
while minimizing the cost functional
\begin{equation} \label{cost_cmtc}
C(t_f,u,K_p)=t_f + K \int_0^{t_f} \|u\|^2 dt,
\end{equation}
with controls satisfying the constraint $u_1^2+u_2^2 \leq 1$, and with the state satisfying constraints on the lateral load factor and the dynamic pressure due to aerodynamic forces
$$
\bar{n} = \frac{\rho |v|^2 S C_N}{2 m g_0} \leq \bar{n}_{max},\quad
\bar{q}= \frac{1}{2} \rho |v|^2 \leq \bar{q}_{max},
$$
where $\rho$ is the air density, $S$ is the reference surface of the launcher, $C_N$ is the lift coefficient approximated by $C_N = C_{N0} + C_{N \alpha} \alpha$ with given constants $C_{N0}$ and $C_{N \alpha}$. $\alpha$ is the angle of attack given by
$$
\alpha = (v_x \sin \theta \cos \psi - v_y \sin \psi + v_z \cos \theta \cos \phi)/v ,
$$
and $|v|$ is the module of the velocity $|v|=\sqrt{v_x^2+v_y^2+v_z^2}$.
Compared to ($\P_S$), a significant additional difficulty comes from the state constraints.

\paragraph{Hard constraint formulation.}
Recall that a state constraint $c(x) \leq 0$ is of order $m$ if
$\hat{g}_i.c=\hat{g}_i \hat{f}.c=\cdots=\hat{g}_i \hat{f}^{m-2}.c=0$ and $g_i f^m.c \ne 0$, $i=1,2$. Here we use the notation of Lie derivatives, see Section \ref{Lie}.
A \emph{boundary arc} is an arc (not reduced to a point) satisfying the system $$
c(x(t))= c^{(1)}(x(t))=\cdots = c^{(m-1)}(x(t))=0,
$$
and the control along the boundary arc is a feedback control obtained by solving
$$
c^{(m)} = \hat{f}^m.c + u_1 \, \hat{g}_1 \hat{f}^{(m-1)}.c + u_2 \, \hat{g}_2 \hat{f}^{(m-1)}.c =0.
$$
After calculations, we find that the constraint on the load factor $\bar{n}$ is of order $2$ and the constraint on the dynamic pressure $\bar{q}$ is of order $3$.

According to the maximum principle with state constraints (see, e.g., \cite{Hartl}), there exists a nontrivial triple of Lagrange multipliers $(p,p^0,\eta)$, with $p^0 \leq 0$, $p \in BV(0,t_f)^{11}$ and $\eta=(\eta_1,\eta_2) \in BV(0,t_f)^2$, where $BV(0,t_f)$ is the set of functions of bounded variation over $[0, t_f]$, such that almost everywhere on $[0,t_f]$
\begin{equation*}
\begin{split}
\dot{x} &= \frac{\partial H(x,p,u,p^0,\eta)}{\partial p}, \\
dp &= -\frac{\partial H(x,p,u,p^0,\eta)}{\partial x} dt - \sum_{i=1}^2 \frac{\partial c_i(x)}{\partial x} d\eta_i,
\end{split}
\end{equation*}
where the Hamiltonian of the problem is
$$
H(x,p,u,p^0,\eta)= \langle p,\hat{f}(x)+ u_1 \hat{g}_1(x) + u_2 \hat{g}_2(x) \rangle + \sum_{i=1}^2 \eta_i c_i(x)+p^0 (1+K\Vert u\Vert^2),
$$
and we have the maximization condition
$$u(t) \in {\rm argmax}_{w} H(x(t),p(t),w,p^0,\eta(t))$$
for almost every $t$.
In addition, we have
$d\eta_i \geq 0$ and $\int_0^{t_f} c_i(x) \, d\eta_i = 0$ for $i=1,2$.

Along a boundary arc, we must have $h_i = \langle p,\hat{g}_i(x) \rangle = 0$, $i=1,2$. Assuming that only the first constraint (which is of order $2$) is active along this boundary arc, and
differentiating twice the switching functions $h_i$, $i=1,2$, we have
$d^2 h_i = \langle p, {\rm ad^2}\hat{f}.\hat{g}_i (x) \rangle dt^2 -  d\eta_1 \cdot ({\rm ad}\hat{f}.\hat{g}_i).c_1 dt$.
Moreover, at an entry point occurring at $t=\tau$, we have
$dh_i(\tau^+)=dh_i(\tau^-)-  d\eta_1 \cdot ({\rm ad}\hat{f}.\hat{g}_i).c_1 =0$, which yields $d\eta_1$. A similar result is obtained at an exit point.

The main drawback of this formulation is that the adjoint vector $p$ is no longer absolutely continuous. A jump $d\eta$ may occur at the entry or at the exit point of a boundary arc, which complexifies significantly the numerical solution.

An alternative approach to address the dynamic pressure state constraint, used in \cite{Corvin,DukemanCalise}, is to design a feedback law that reduces the commanded throttle based on an error signal. According to \cite{DukemanCalise}, this approach works well when the trajectory does not violate too much the maximal dynamic pressure constraint, but it may cause instability if the constraint is violated significantly. In any case the derived solutions are suboptimal.

Another alternative is the \emph{penalty function method} (also called \emph{soft constraint method} ).
The soft constraint consists in introducing a penalty function to discard solutions entering the constrained region \cite{DenhamArthur,MarkopoulosCalise,Trelatbook}.
For the problem ($\P_A$), this soft constraint method is well suited in view of a continuation procedure starting from an unconstrained solution. This initial solution generally violates significantly the state constraint. The continuation procedure aims at reducing progressively the infeasibility.

\paragraph{Soft constraint formulation.}
The problem ($\P_A$) is recast as an unconstrained optimal control problem by adding a penalty function to the cost functional defined by \eqref{cost_cmtc}. The penalized cost is
\begin{equation*} 
C(t_f,u,K_p)=t_f + K \int_0^{t_f} \|u\|^2 dt + K_p \int_0^{t_f} P(x(t))dt,
\end{equation*}
where the penalty function $P(\cdot)$ for the state constraints is defined by
$$
P(x) = (\max(0,\bar{n} - \bar{n}_{max}))^2 + (\max(0,\bar{q} - \bar{q}_{max}))^2.
$$
The constraint violation is managed by tuning the parameter $K_p$.
For convenience we still denote this unconstrained problem by ($\P_A$) and we apply the PMP.

\paragraph{Application of the PMP.}
\label{par_pmp}
The Hamiltonian is now given by
$$
H(x,p,p^0,u) = \langle p, \hat{f}(x)\rangle+u_1 \langle p, \hat{g}_1(x)\rangle+u_2 \langle p, \hat{g}_2(x) \rangle + p^0 (1+ K \|u\|^2 + K_p P(x)).
$$
The adjoint equation is
\begin{equation} \label{adjointsys_cmtc}
\dot{p}(t) = -\frac{\partial H}{\partial x}(x(t),p(t),p^0,u(t)) ,
\end{equation}
where we have set $p =(p_{r_x},p_{r_y},p_{r_z},p_{v_x},p_{v_y},p_{v_z},p_{\theta},p_{\psi},p_{\phi},p_{\omega_x},p_{\omega_y})$.
Let $h = (h_1,h_2)$ be the switching function and let
\begin{align*}
& h_1 (t)= \langle p(t), \hat{g}_1(x(t))\rangle=b p_{\omega_y}(t), \\
& h_2 (t)= \langle p(t), \hat{g}_2(x(t))\rangle=-b p_{\omega_x}(t).
\end{align*}
The maximization condition of the PMP gives
\begin{equation} \label{OCPn_u1_cmtc}
u = \begin{cases}
(h_1,h_2)  / (2 K) & \textrm{if}\ \|h\| \leq 2 K,\\
(h_1,h_2)/\|h\| & \textrm{if}\  \|h\| > 2 K.
\end{cases}
\end{equation}
The transversality condition $p(t_f) \perp T_{x(t_f)} M_1$, where $T_{x(t_f)}M_1$ is the tangent space to $M_1$ at the point $x(t_f)$, yields the additional conditions
$$
p_{v y}(t_f) \sin \psi_f = p_{v x}(t_f) \sin \theta_f \cos \psi_f+ p_{v z}(t_f) \cos \theta_f \cos \psi_f
$$
and 
$$p_{r_x}(t_f) = p_{r_y}(t_f) = p_{r_z}(t_f)=0.$$
The final time $t_f$ being free and the system being autonomous, we have in addition that 
$$H(x(t),p(t),p^0,u(t))=0,$$ 
almost everywhere on $[0,t_f]$. As previously we can assume $p^0=-1$.

The optimal control given by \eqref{OCPn_u1_cmtc} is regular unless $K=0$ and $\|h(t)\| =0$, in which case it becomes \emph{singular}.
As before the term $K \int_0^{t_f} \|u(t)\|^2 dt$ in the cost functional \eqref{cost_cmtc} is used to avoid chattering \cite{Marchal,Fuller,Robbins,Zelikin1,Zelikin2}, and the exact minimum time solution can be approached by decreasing step by step the value of $K\geq 0$ until the shooting method possibly fails due to chattering.

\paragraph{Solution algorithm and comparison with ($\P_S$)}
We aim at extending the continuation strategy developed for ($\P_S$) in order to address ($\P_A$). Comparing ($\P_A$) with ($\P_S$), we see that in ($\P_A$):
\begin{itemize}
\item[(a)] the position of the launcher is added to the state vector;
\item[(b)] the gravity acceleration $\vec{g}$ depends on the position and the aerodynamic forces (lift force $\vec{L}$ and the drag force $\vec{D}$) are considered;
\item[(c)] the cost functional is penalized by the state constraints violation;
\end{itemize}
 

\medskip

Regarding the point (a), we need embedding the solution of ($\P_0$) into a larger dimension problem with the adjoint variable of the position $\vec{p}_r = (p_{rx},p_{ry},p_{rz})^\top$ being zero.
More precisely, consider the following problem, denoted by ($\P_0^H$), in which the position and the velocity are considered
\begin{align*}
& \dot{\vec{r}} = {\vec{v}}, \quad \dot{\vec{v}} = a \vec{e} +\vec{g}_0, \\
& {\vec{r}}(0)={\vec{r}}_0, \quad \vec{v}(0)=\vec{v}_0,\quad \vec{v}(t_f) // \vec{w},\\
& \Vert \vec{w} \Vert =1,\\
& \min t_f.
\end{align*}

The solution of ($\P_0^H$) is retrieved from the solution of ($\P_0$) completed by the new state components,
\begin{equation*}
 t_f=\frac{-a_2 + \sqrt{a_2^2-4a_1a_3}}{2 a_1},\qquad \vec{p}_r = \vec{0},\qquad
 \vec{p}_v=\frac{-p^0}{a + \langle \vec{e}^{\ast},\vec{g} \rangle } \vec{e}^{\ast}\, ,
\end{equation*}
and the optimal control is
$$
\vec{e}=\vec{e}^\ast
=\frac{1}{a}\left(\frac{k \vec{w}-\vec{v}_0}{t_f}-\vec{g}_0\right),
$$
with
$$k=\langle \vec{v}_0,\vec{w} \rangle + \langle \vec{g}_0,\vec{w} \rangle t_f,$$
$$a_1 = a^2-\Vert  \langle \vec{g}_0,\vec{w} \rangle \vec{w}-\vec{g}_0 \Vert ^2,$$
$$a_2 = 2 ( \langle \vec{v}_0,\vec{w} \rangle \langle \vec{g}_0,\vec{w} \rangle - \langle \vec{v}_0,\vec{g}_0 \rangle),$$
and
$$a_3 = - \Vert  \langle \vec{v}_0,\vec{w} \rangle \vec{w}-\vec{v}_0 \Vert ^2.$$

We use this solution as the initialization of the continuation procedure for solving ($\P_A$).

\medskip

The point (b) can be addressed with a new continuation parameter $\lambda_4$ introducing simultaneously the variable gravity acceleration, the aerodynamic forces and the atmospheric density $\rho$ (exponential model) as follows:
\begin{equation*}
\begin{split}
& \dot{v}_x= a \sin \theta \cos \psi + g_{0x} (1-\lambda_4) + \lambda_4 g_x  + \lambda_4 \frac{D_x+L_x}{m}, \\
& \dot{v}_y= - a \sin \psi + g_{0y} (1-\lambda_4) + \lambda_4 g_y + \lambda_4 \frac{D_y+L_y}{m},\\
& \dot{v}_z= a \cos \theta \cos \psi + g_{0z} (1-\lambda_4) + \lambda_4 g_z + \lambda_4 \frac{D_z+L_z}{m},
\end{split}
\end{equation*}
and
\begin{equation*}
\begin{split}
\rho(t) = \rho_0 \big( (1-\lambda_4)  \exp(-(\sqrt{(R_E+r_x(0))^2+r_y(0)^2+r_z(0)^2}-R_E)/h_s) \\
+ \lambda_4 \exp(-(\sqrt{(R_E+r_x)^2+r_y^2+r_z^2}-R_E)/h_s) \big),
\end{split}
\end{equation*}
where $R_E = 6378137\ m$ is the radius of the Earth, $h_s=7143\ m$, $\rho_0=1.225\ kg/m^3$, and $g_x$, $g_y$, $g_z$ are given by
\begin{equation*} 
(g_x,g_y,g_z)^\top =  - \frac{ g_0 \sqrt{(R_E+r_x(0))^2+r_y(0)^2+r_z(0)^2} }{ \sqrt{(R_E+r_x)^2+r_y^2+r_z^2}}  (\cos l_2, \sin l_1 \sin l_2, \cos l_1 \sin l_2)^\top,
\end{equation*}
with 
$$g_0=\sqrt{g_{x0}^2 +g_{y0}^2+ g_{z0}},$$
and
$$\tan l_1 = r_y/r_x,\qquad \tan l_2 = \sqrt{r_y^2+r_z^2}/(r_x+R_E).$$
The parameter $\lambda_4$ acts only on the dynamics. Applying the PMP, $\lambda_4$ appears explicitly in the adjoint equations, but not in the shooting function.

\medskip

Finally, regarding the point (c), the penalty parameter $K_p$ in the cost functional \eqref{cost_opcr} has to be large enough in order to produce a feasible solution.
Unfortunately, too large values of $K_p$ may generate ill conditioning and raise numerical difficulties. In order to obtain an adequate value for $K_p$, a simple strategy \cite{FrangosSnyman,Snyman} consists in starting with a quite small value of $K_p=K_{p0}$ and solving a series of problems with increasing $K_p$. The process is stopped as soon as $\| c(x(t)) \| < \epsilon_c$, for every $t\in [0,t_f]$, for some given tolerance $\epsilon_c>0$.

\bigskip

For convenience, we define the \emph{exo-atmospheric pull-up maneuver problem} ($\P_A^{exo}$) as ($\P_A$) without state constraints and without aerodynamic forces and the \emph{unconstrained pull-up maneuver problem} ($\P_A^{unc}$) as ($\P_A$) without state constraints.

We proceed as follows:
\begin{itemize}
\item First, we embed the solution of ($\P_0 $), into the larger dimension problem ($\P_A$). This problem is denoted ($\P_0^H$).
\item Then, we pass from ($\P_0^H$), to ($\P_A$) by using a numerical continuation procedure, involving four continuation parameters:
two parameters $\lambda_1$ and $\lambda_2$ introduce the terminal conditions \eqref{ic_cmtc}-\eqref{fc_cmtc} into ($\P_A^{exo}$); $\lambda_4$ introduces the variable gravity acceleration and the aerodynamic forces in ($\P_A^{unc}$); $\lambda_5$ introduces the soft constraints in ($\P_A$).
\end{itemize}

The overall continuation procedure is depicted on Figure \ref{homotopie_cmtc}. The final step of the procedure is to increase $\lambda_3$ (or equivalently decrease $K$) in order to minimize the maneuver duration.
\begin{figure}[ht]
\centering
\includegraphics[scale=1.1]{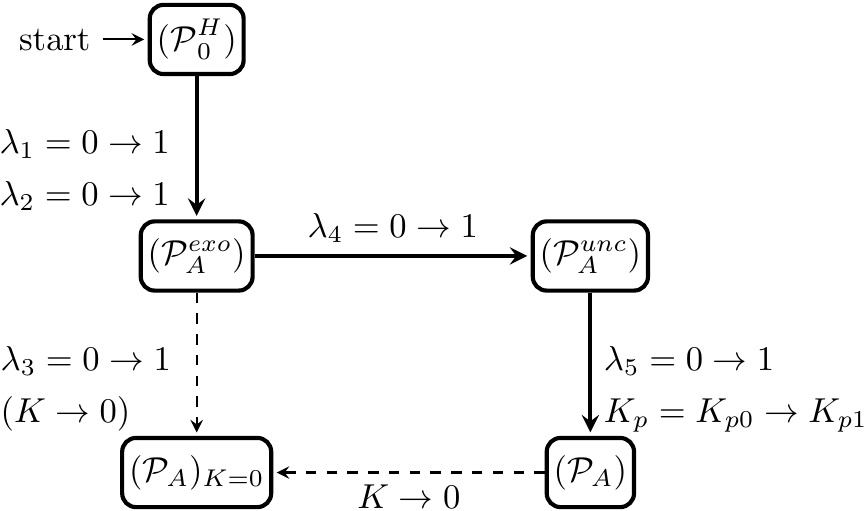}           
\caption{Continuation procedure for solving ($\P_A$).}
\label{homotopie_cmtc}
\end{figure}

More precisely, we have to solve the following problem with continuation parameters $\lambda_i$, $i=1,2,4,5,3$

$$
\min t_f + (1-\lambda_3) \int_0^{t_f} \|u\|^2 dt + \lambda_5 K_p \int_0^{t_f} P(x(t)) dt,
$$
subject to the dynamics
\begin{align*}
& \dot{r}_x = v_x, \qquad \dot{r}_y = v_y, \qquad \dot{r}_z = v_z, \\
& \dot{v}_x= a \sin \theta \cos \psi + g_{0x} (1-\lambda_4) + \lambda_4 g_x  + \lambda_4 \frac{D_x+L_x}{m}, \\
& \dot{v}_y= - a \sin \psi + g_{0y} (1-\lambda_4) + \lambda_4 g_y + \lambda_4 \frac{D_y+L_y}{m}, \\
& \dot{v}_z= a \cos \theta \cos \psi + g_{0z} (1-\lambda_4) + \lambda_4 g_z + \lambda_4 \frac{D_z+L_z}{m}, \\
& \dot{\theta}=(\omega_x \sin \phi + \omega_y \cos \phi)/ \cos \psi , \\
& \dot{\psi}=\omega_x \cos \phi - \omega_y \sin \phi , \\
& \dot{\phi}= (\omega_x \sin \phi + \omega_y \cos \phi) \tan \psi , \\
& \dot{\omega }_x= -\bar{b} u_2, \\
& \dot{\omega }_y=  \bar{b} u_1 ,
\end{align*}
and with initial conditions
\begin{align*}
& \vec{r}(0) = (r_{x0},r_{y0},r_{z0})^\top,\\
& \vec{v}(0) = (v_{x0},v_{y0},v_{z0})^\top, \\
& \theta(0) = \theta^\ast (1-\lambda_1) + \theta_0 \lambda_1,\\
& \psi(0) = \psi^\ast (1-\lambda_1) + \psi_0 \lambda_1,\\
& \phi(0) = \phi_0 \lambda_1, \\
& \omega_x(0)= \omega_{x0} \lambda_1,\\
& \omega_y(0)= \omega_{y0} \lambda_1,
\end{align*}
and final conditions
\begin{align*}
& \vec{r}(t_f) \, \textrm{free},\\
& \vec{v}(t_f) \perp \hat{z}_b, \\
& \theta(t_f) = \theta_e (1-\lambda_2) + \theta_f \lambda_2,\\
& \psi(t_f) = \psi_e (1-\lambda_2) + \psi_f \lambda_2,\\
& \phi(t_f) = \phi_e (1-\lambda_2) + \phi_f \lambda_2,\\
& \omega_x(t_f)= \omega_{xe} (1-\lambda_2) + \omega_{xf} \lambda_2,\\
& \omega_y(t_f)= \omega_{ye} (1-\lambda_2) + \omega_{xf} \lambda_2.
\end{align*}
The attitude angles $\theta_e$, $\psi_e$, $\phi_e$, $\omega_{xe}$, and $\omega_{ye}$ are those obtained at the end of the first continuation on $\lambda_1$. $\theta^\ast$, $\psi^\ast$ are the explicit solutions of ($\P_0^H$).

\medskip

These successive continuations are implemented using the PC continuation combined with the multiple shooting method. Some additional enhancements regarding the inertial frame choice and the Euler angle singularities help improving the overall robustness of the solution process.

\paragraph{Multiple shooting.}
The unknowns of this shooting problem are $p(0) \in \R^{11}$, $t_f \in \R$, and $z_i=(x_i,p_i) \in \R^{22}$, $i=1,\cdots,N-1$, where $z_i$ are the node points of the multiple shooting method (see Section \ref{subsec_indir}).
We set $Z=(p(0),t_f,z_i)$, and let $E=(\theta,\psi,\phi)$, $\omega = (\omega_x,\omega_y)$, $p_r=(p_{rx},p_{ry},p_{rz})$, $p_E=(p_{\theta},p_{\psi},p_{\phi})$, and $p_{\omega}=(p_{\omega x},p_{\omega y})$. Then, the shooting function with the continuation parameter $\lambda_1$ is given by
\begin{equation*} 
G_{\lambda_1} = \begin{pmatrix}
         v_z(t_f) \sin \psi_f + v_y(t_f) \cos \theta_f \cos \psi_f \\
         v_z(t_f) \sin \theta_f - v_x(t_f) \cos \theta_f  \\
        p_{v y}(t_f) \sin \psi_f- (p_{v x}(t_f) \sin \theta_f \cos \psi_f+ p_{v z}(t_f) \cos \theta_f \cos \psi_f)\\
        p_r(t_f),\\
        p_{\omega}(t_f)\\
        p_E (t_f)\\
        H(t_f)\\
        \left\{z_i(t_i^-)=z_i(t_i)^+,i=1,\cdots,N-1 \right\} 
        \end{pmatrix},
\end{equation*}
where the Hamiltonian is given by
$$
H=\langle p, f(x)\rangle+u_1 \langle p, g_1(x)\rangle+u_2 \langle p, g_2(x) \rangle + p^0 (1+ (1-\lambda_3) K \|u\|^2 + \lambda_5 K_p P(x)).
$$

The shooting function with the continuation parameter $\lambda_2$ is 
\begin{equation*} 
 G_{\lambda_2} = \begin{pmatrix}
  v_z(t_f) \sin \psi_f + v_y(t_f) \cos \theta_f \cos \psi_f,\\
   v_z(t_f) \sin \theta_f - v_x(t_f) \cos \theta_f  \\
  p_{v y}(t_f) \sin \psi_f- (p_{v x}(t_f) \sin \theta_f \cos \psi_f+ p_{v z}(t_f) \cos \theta_f \cos \psi_f) \\
  {E}(t_f)-(1-\lambda_2) {E}_e - \lambda_2 {E}_f \\
 \omega(t_f) - (1-\lambda_2) \omega_e -\lambda_2 \omega_f \\
  p_{r }(t_f),\, H(t_f)\\
  \left\{ z_i(t_i^-)=z_i(t_i)^+,i=1,\cdots,N-1 \right\}
  \end{pmatrix},
\end{equation*}
and the shooting functions $G_{\lambda_4} $ and $G_{\lambda_5} $ are identical to $G_{\lambda_2} $.

\paragraph{PC continuation.}
The predictor-corrector continuation requires the calculation of the Jacobian matrix $J_G$ (see Section \ref{subsec_pccontinuation}) which is computationally expenssive. In order to speed up the process, an approximation is used based on the assumption of no conjugate point. According to \cite{CaillauDaoud}, the first turning point of $\lambda(\bar{s})$ (where $\frac{d  \lambda}{ds}(\bar{s})=0$ and $\frac{d^2 \lambda}{ds^2}(\bar{s}) \neq 0 $) corresponds to a conjugate point (the first point where extremals lose local optimality).
If we assume the absence of the conjugate point, there is no turning point for $\lambda(s)$, and $\lambda$ increases monotonically along the zero path.
Knowing three zeros $(Z_{i-2},\lambda_{i-2})$, $(Z_{i-1},\lambda_{i-1})$ and $(Z_{i},\lambda_{i})$, and let
$s_1 = \| (Z_{i-1},\lambda_{i-1}) - (Z_{i-2},\lambda_{i-2})\|$,
$s_2 = \| (Z_i,\lambda_i) - (Z_{i-2},\lambda_{i-2})\|$,
$s_3 = \| (Z_i,\lambda_i) - (Z_{i-1},\lambda_{i-1})\|$,
we can approximate the tangent vector $t(J_G)$ by
\begin{equation} \label{tJG}
t(J_G) = \frac{(Z_i,\lambda_i)-(Z_{i-1},\lambda_{i-1})}{s_2-s_1} \frac{|s_2-s_1|}{|s_3|}.
\end{equation}
When the step length $h_s$ is small enough, this approximation yields a predicted point \eqref{cont_pred} very close to the true zero.

\paragraph{Change of Frame.}
\label{precond}
Changing the inertial reference frame can improve the problem conditioning and enhance the numerical solution process. The new frame $S_R^\prime$ is defined from the initial frame $S_R$ by two successive rotations of angles $(\beta_1,\beta_2)$. 
The problem ($\P_A$) becomes numerically easier to solve when the new reference frame $S_R^\prime$ is adapted to the terminal conditions. However we do not know a priori which reference frame is the best suited. We propose to choose a reference frame associated to $(\beta_1,\beta_2)$ such that $\psi^\prime_f=-\psi^\prime_0$ and $| \psi^\prime_f |+| \psi^\prime_0 |$ being minimal (the subscribe $\prime$ here means the new variable in $S_R^\prime$).
This choice centers the terminal values on the yaw angle on zero. Thus we can hope that the solution remains far from the Euler angle singularities occurring when $\psi \to \pi/2+ k\pi$.

This frame rotation defines a nonlinear state transformation, which acts as a preconditionner.  We observe from numerical experiments that it actually enhances the robustness of the algorithm.
The reader is referred to \cite{ztc2} for more details of the change of frame.

\paragraph{Singularities of Euler Angles.} \label{SingularEuler}
The above frame change is not sufficient to avoid Euler angle singularities in all cases. Smoothing the vector fields at these singular configurations is another enhancement improving the overall robustness. The state and costate equations are smoothened as follows.
Assuming first that $\dot{\theta}$ is bounded, we have $\omega_x \sin \phi + \omega_y \cos \phi \to 0$ when $\psi \to \pi/2+ k\pi$. Since
$$\dot{\theta} \dot{\phi} = \lim_{\psi \to \pi/2+ k\pi} (\omega_x \sin \phi + \omega_y \cos \phi)^2 \sin \psi \to 0$$ 
and $\dot{\theta}/\dot{\phi} \to 1$ as $\psi \to \pi/2+ k\pi$, we can smoothen the state equations by $\dot{\theta} = \dot{\phi} = 0$ when $\psi \to \pi/2+ k\pi$.
Assuming then that
$- \frac{p_{\theta}+p_{\phi} \sin \psi}{\cos \psi} \to A < \infty$ as $\psi \to \pi/2+ k\pi$,
and taking the first-order derivatives of the numerator and denominator 
$$
A = \lim_{\psi \to \pi/2+ k\pi}   - \frac{p_{\theta}+p_{\phi} \sin \psi}{\cos \psi}
=\lim_{\psi \to \pi/2+ k\pi}   \frac{\dot{p}_{\theta}+\dot{p}_{\phi} \sin \psi+p_{\phi} \cos \psi \dot{\psi}}{\sin \psi \dot{\psi}} = -A
$$
we obtain $A = 0$. We can smoothen the costate equations by $\dot{p}_{\theta} = 0$, $\dot{p}_{\phi} = 0$, $\dot{p}_{\psi} = a \sin \theta p_{v x} + a \cos \theta p_{v z}$, $\dot{p}_{\omega_x} = -p_{\psi} \cos \phi$, $\dot{p}_{\omega_y} = p_{\psi} \sin \phi$.
Summing up, at points $\psi \to \pi/2+k\pi$, the attitude equations in system \eqref{sys_full} and \eqref{adjointsys_cmtc} become
\begin{equation}\label{sys_atti_limit-sys_adj_limit}
\begin{split}
& \dot{\theta}=0 , \\
& \dot{\psi}=\omega_x \cos \phi - \omega_y \sin \phi ,\\
& \dot{\phi}= 0 ,\\
& \dot{\omega }_x= -\bar{b} u_2,\\
& \dot{\omega }_y=  \bar{b} u_1, \\
& \dot{p}_\theta=0 , \\
& \dot{p}_\psi= a \sin \theta p_{v_x} + a \cos \theta p_{v_z},\quad \dot{p}_\phi= 0 ,\\
& \dot{p}_{\omega_x}= -p_{\psi} \cos \phi,\quad \dot{p}_{\omega_y}= p_{\psi} \sin \phi .
\end{split}
\end{equation}
These equations \eqref{sys_atti_limit-sys_adj_limit} are used close to the singularities.

\paragraph{Algorithm}
We describe the whole numerical strategy of solving ($\P_A$) in the following algorithm.
\begin{algorithm}[h]
\SetAlgoLined
\KwResult{The solution of the problem ($P_A$)}
$\cdot$ Change of frame: compute $(\beta_1,\beta_2)$ and the new initial condition $x(0)=x_0^\prime$\;
$\cdot$ Solve ($P_0^H$), to get a solution $Z_0$\;
$\cdot$ Initialize the multiple shooting method with $Z_0$ and $\lambda_i = 0$, $i=1,\cdots 5$\;
\For{ $i=1,2,4,5,(3)$}{
\While{$\lambda_i \leq 1$ and $\triangle\lambda_{min}^i\leq \triangle\lambda_i \leq \triangle\lambda_{max}^i$}
{
  (Predictor) Predict a point $(\tilde{Z},\tilde{\lambda}_i)$ according to \eqref{cont_pred} and \eqref{tJG}\;
  (Corrector) Find the solution $(\bar{Z},\bar{\lambda}_i)$ of $G_{\lambda i}(\tilde{Z},\tilde{\lambda}_i)=0$\;
\eIf{successful}{
  $(Z,\lambda_i)=(\bar{Z},\bar{\lambda}_i)$\;}{
  Reduce the step-length $h_s$\;}
}
\eIf{successful}{
   The $\lambda_i$-continuation is successful\;}{
   The $\lambda_i$-continuation has failed\;}
\caption{Prediction-Corrector continuation for ($\P_A$)}
}
\label{alg_pc_cmtc}
\end{algorithm}

\subsection{Numerical Results of Solving ($\P_A$)}
The algorithm \ref{alg_pc_cmtc} is first applied to a pull-up maneuver of an airborne launch vehicle just after its release from the carrier. We present some statistical results showing robustness of our algorithm. A second example considers the three-dimensional reorientation maneuver of a launch vehicle upper stage after a stage separation.

\subsubsection{Pull-Up Maneuvers of an airborne launch vehicle (AVL)}
We consider a pull-up maneuver of an airborne launch vehicle close to the Pegasus configuration : $a = 15.8$, $b = 0.2$, $S=14\, m^2$, $C_{x0} = 0.06$, $C_{x \alpha}=0$, $C_{z0}= 0$, and $C_{z \alpha} = 4.7$. Let $\bar{n}_{max}=2.2 g$ and $\bar{q}_{max}=47 \,kPa$.
The initial conditions \eqref{ic_cmtc} correspond to the engine ignition just after the release.
$$
r_{x 0}=11.9\,km, \quad r_{y 0}=r_{z 0}=0, \quad v_{0}=235\, m/s, \quad \theta_{v 0}=-10^\circ,
$$
$$
\psi_{v0}=0^\circ, \quad \theta_0=-10^\circ, \quad \psi_0=\phi_0=0, \quad \omega_{x 0}=\omega_{y 0}=-1^\circ/s,
$$
The final conditions \eqref{fc_cmtc} correspond to the beginning of the atmospheric ascent flight at zero angle of attack.
$$
\theta_f= 42^\circ, \quad \psi_f=10^\circ,\quad \phi_f=0, \quad \omega_{x f}=\omega_{y f}=0.
$$
Such pull-up maneuvers are generally planar ($\psi_f=0^\circ$). Here we set $\psi_f=10^\circ$ in order to show that the algorithm can also deal efficiently with non-planar pull-up maneuvers.

\begin{figure} [H]
\centering
\includegraphics[width=\textwidth]{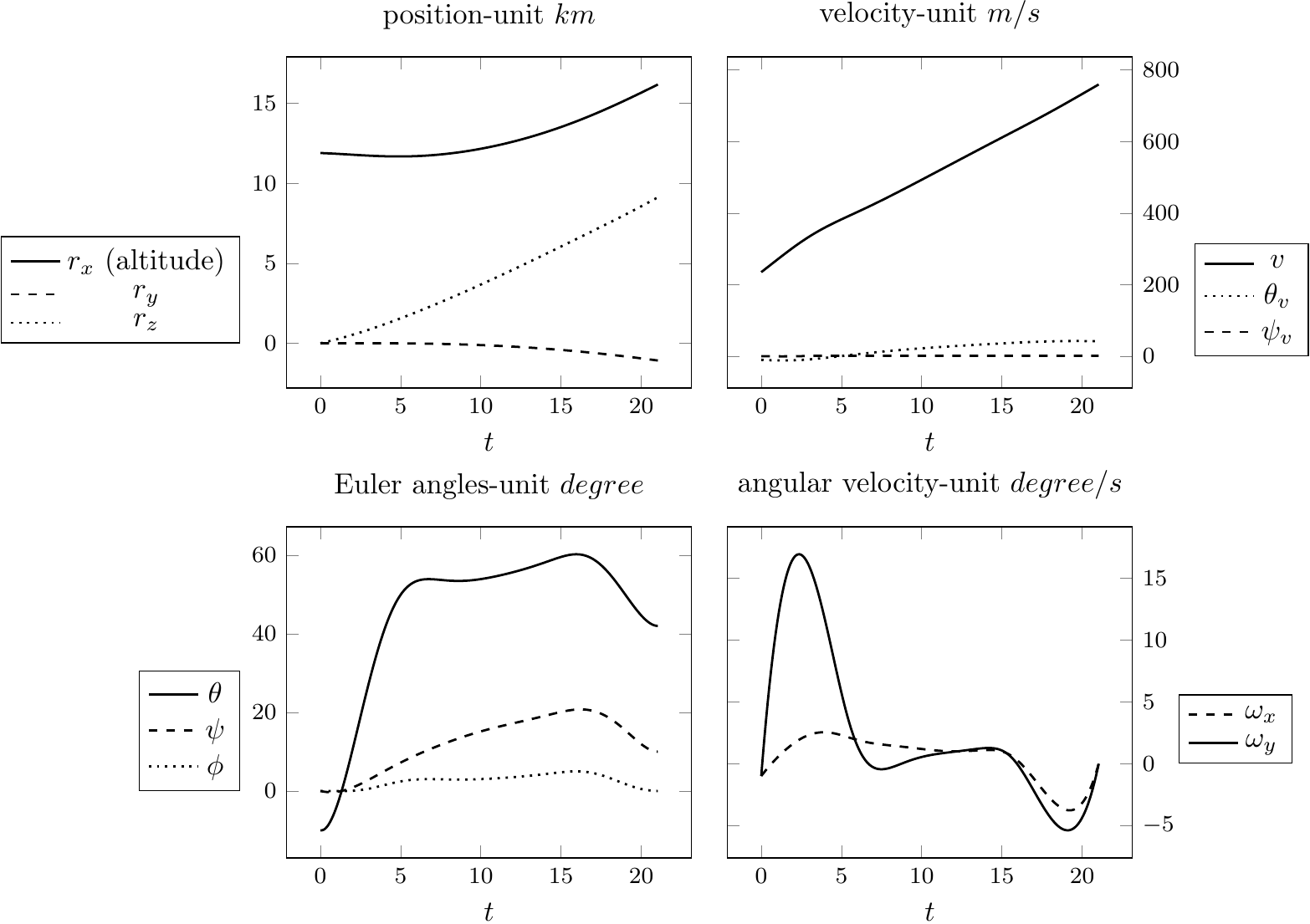}   
\caption{Time history of state variable $x(t)$ during the pull-up maneuver.}
\label{state_air}
\end{figure}

\begin{figure} [H]
\centering
\includegraphics[scale=0.7]{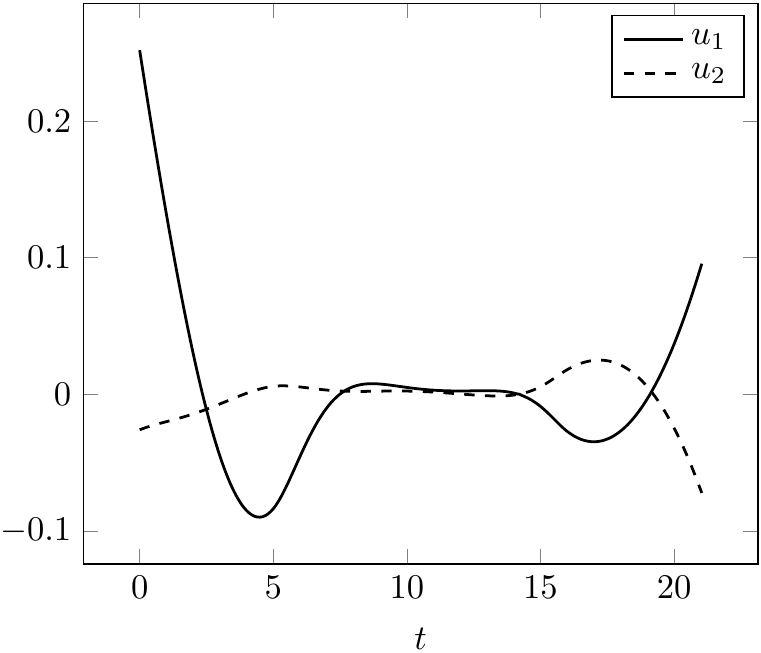}   
\caption{Time history of control variable $u(t)$ during the pull-up maneuver.}  
\label{control_air}
\end{figure}

\begin{figure} [H]
\centering
\includegraphics[width=0.8\textwidth]{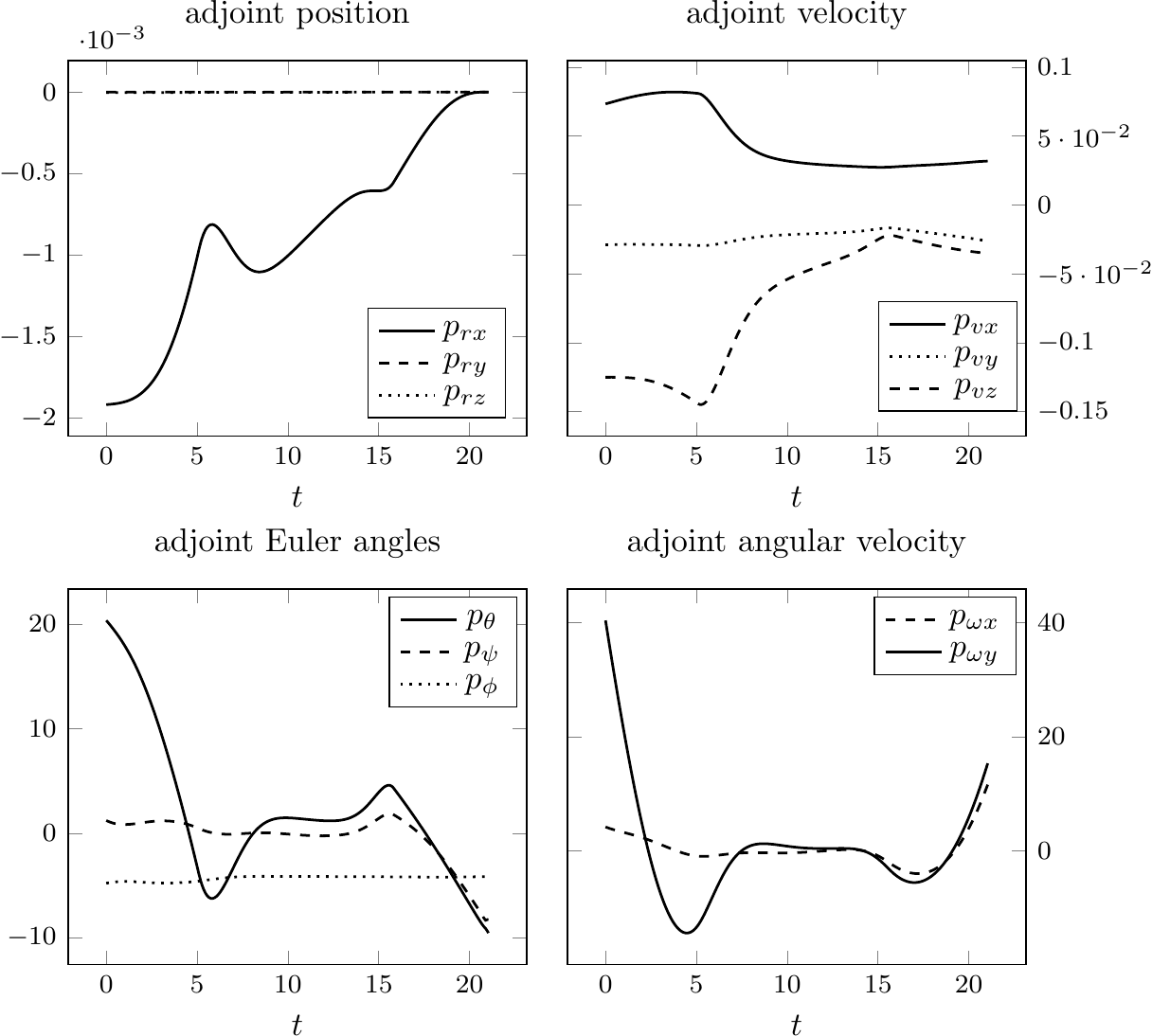}   
\caption{Time history of the adjoint variable $p(t)$ during the pull-up maneuver.}  
\label{adjoint_air}
\end{figure}

\begin{figure} [H]
\centering
\includegraphics[width=0.7\textwidth]{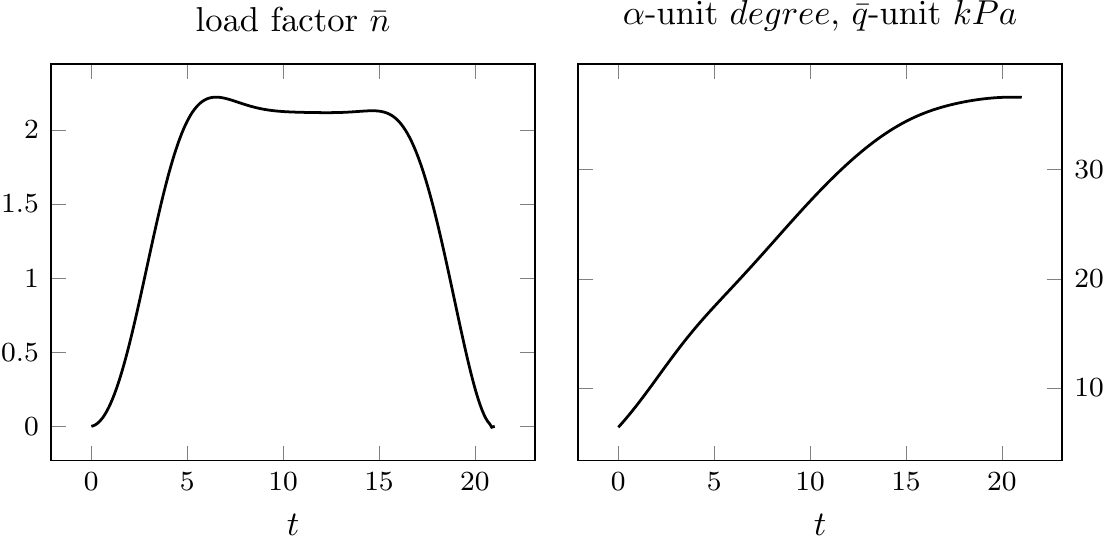}    
\caption{Time history of the constraints $c(x(t))$ during the pull-up maneuver.}  
\label{constraint_air}
\end{figure}

The multiple shooting method is applied with three node points.
The components of the state variable $x$ and the control $u$ are plotted on Figures \ref{state_air} and \ref{control_air}, the components of the adjoint variable $p$ are plotted on Figure \ref{adjoint_air}, the time histories of the load factor $\bar{n}$ and of the dynamic pressure $\bar{q}$ are plotted on Figure \ref{constraint_air}. The position components are given in the geographic local frame with the vertical along the first axis (denoted x, not to be confused with the state vector). The control vector first component $u_1$ lies mainly in the trajectory plane and it acts mainly on the pitch angle.


We observe on Figure \ref{constraint_air} a boundary arc on the load factor constraint near the maximal level $\bar{n}_{max}=2.2 g$.
This corresponds on Figure \ref{adjoint_air} to the switching function $h(t)=b(p_{\omega_y},-p_{\omega_x})$ being close to zero.
Comparing Figs. \ref{control_air} and \ref{adjoint_air}, we see that the control follows the form of the switching function. On the other hand, the state constraint of the dynamic pressure is never active.

We observe also on Figure \ref{adjoint_air} a steeper variation of $p_{\theta}(t)$ at $t=5.86\,s$. The penalty function $P(x)$ starts being positive at this date and adds terms in the adjoint differential equation.

Running this example requires $24.6 \, s$ to compute the optimal solution, with CPU: Intel(R) Core(TM) i5-2500 CPU 3.30GHz; Memory: 3.8 Gio; Compiler: gcc version 4.8.4 (Ubuntu 14.04 LTS).
The number of nodes for the multiple shooting has been set to 3 from experiments. Passing to four node increases the computing time to $31.2\,s$ without obvious robustness benefit.

We next present some statistical results obtained with the same computer settings.

\paragraph{Statistical results.}
($\P_A$) is solved for various terminal conditions.
The initial and final conditions are swept in the range given in Table \ref{table_statistic}.
The last cell of the table indicates that the initial angle of attack is bounded to 10 degrees in order to exclude unrealistic cases.
\begin{table} [ht]
\centering
\caption{Parameter ranges.}
\begin{tabular}{|c|c|c|c|c|}
\hline
$v_0$ & $\theta_{v0}$ & $\psi_{v0}$ & $\theta_{0}$  & $\psi_{0}$  \\
\hline
fixed $ 0.8 \, Mach$ & $[-10,0]^\circ$ & fixed $0^\circ$ & $[-10,10]^\circ$ & fixed $0^\circ$\\
\hline
 $\theta_{f}$ &  $\psi_f$ & $\omega_{x 0}$ & $\omega_{y 0}$ & $\theta_0-\theta_{v0}$\\
\hline
 $[20,80]^\circ$ & $[-10,10]^\circ$ & $[-2,2]^\circ/s$ & $[-2,2]^\circ/s$ & $[0,10]^\circ$\\
\hline
\end{tabular}
\label{table_statistic}
\end{table}
For each variable, we choose a discretization step and we solve all possible combinations resulting from this discretization (factorial experiment). The total number of cases is $1701$. All cases are run with the penalty parameter varying from $K_{p 0}=0.1$ to $K_{p 1}=100$ during the third continuation. For each continuation stage the number of simulations is limited to $200$.

The $1701$ cases are run for different settings of the number of nodes ($N=0$ or $N=2$) and of the regularization parameter ($K=800$ or $K=1000$). 

The statistical results are reported in Table \ref{table_statistic_res1}-\ref{table_statistic_res3}.

\begin{table} [H]
\centering
\caption{Statistical results ($N=2$ and $K=8 \times 10^2$).}
\begin{tabular}{|c|c|c|}
\hline
  & planar & non-planar \\
\hline
 Number of cases & 567 & 1134  \\
\hline
  Rate of success (\%) & 89.07 & 80.04  \\
\hline
 Number of failure cases & & \\
  \hspace{0.5cm}- In $\lambda_1$-continuation & 0 & 14 \\
  \hspace{0.5cm}- In $\lambda_2$-continuation & 21 & 172  \\
  \hspace{0.5cm}- In $\lambda_4$-continuation & 41 & 26  \\
  \hspace{0.5cm}- In $\lambda_5$+$K_p$-continuation & 0  & 10  \\      
\hline
  Average execution time (s) &  &  \\
  - Total & 26.94 & 44.05  \\
  \hspace{0.5cm}- In $\lambda_1$-continuation & 0.49 & 0.48  \\
  \hspace{0.5cm}- In $\lambda_2$-continuation & 2.07 & 2.37  \\
  \hspace{0.5cm}- In $\lambda_4$-continuation & 2.54 & 2.99  \\
  \hspace{0.5cm}- In $\lambda_5$+$K_p$-continuation & 23.16 & 37.35  \\
\hline
\end{tabular}
\label{table_statistic_res1}
\end{table}

\begin{table} [H]
\centering
\caption{Statistical results ($N=2$ and $K=1 \times 10^3$).}
\begin{tabular}{|c|c|c|}
\hline
  & planar & non-planar \\
\hline
 Number of cases & 567 & 1134  \\    
\hline
  Rate of success (\%) & 85.89 & 86.94  \\
\hline
 Number of failure cases & & \\
  \hspace{0.5cm}- In $\lambda_1$-continuation & 0 & 4 \\
  \hspace{0.5cm}- In $\lambda_2$-continuation & 36 & 120  \\
  \hspace{0.5cm}- In $\lambda_4$-continuation & 44 & 16  \\
  \hspace{0.5cm}- In $\lambda_5$+$K_p$-continuation & 0  & 8  \\   
\hline
  Average execution time (s) &  &  \\
  - Total & 26.55 & 47.96  \\
  \hspace{0.5cm}- In $\lambda_1$-continuation & 0.49 & 0.51  \\
  \hspace{0.5cm}- In $\lambda_2$-continuation & 2.12 & 2.40  \\
  \hspace{0.5cm}- In $\lambda_4$-continuation & 2.71 & 2.74  \\
  \hspace{0.5cm}- In $\lambda_5$+$K_p$-continuation & 22.60 & 42.28  \\
\hline  
\end{tabular}
\label{table_statistic_res2}
\end{table}

Tables \ref{table_statistic_res1}-\ref{table_statistic_res2} show the results with a multiple shooting using $2$ nodes, with different values of the regularization parameter $K$.
The algorithm appears fairly robust with respect to the terminal conditions.
The choice of the regularization parameter $K$ affects the resolution results:
(i) the rate of success increases (resp. decreases) in the non-planar case (resp. planar case) when $K$ increases from $K=800$ to $K=1000$;
(ii) in term of the execution time, we see that in both cases, it is faster to get a result in planar case than in non-planar case, and most time is devoted to deal with the state constraints during the last continuation.

\begin{figure} [H]
\centering
\includegraphics[width=0.7\textwidth]{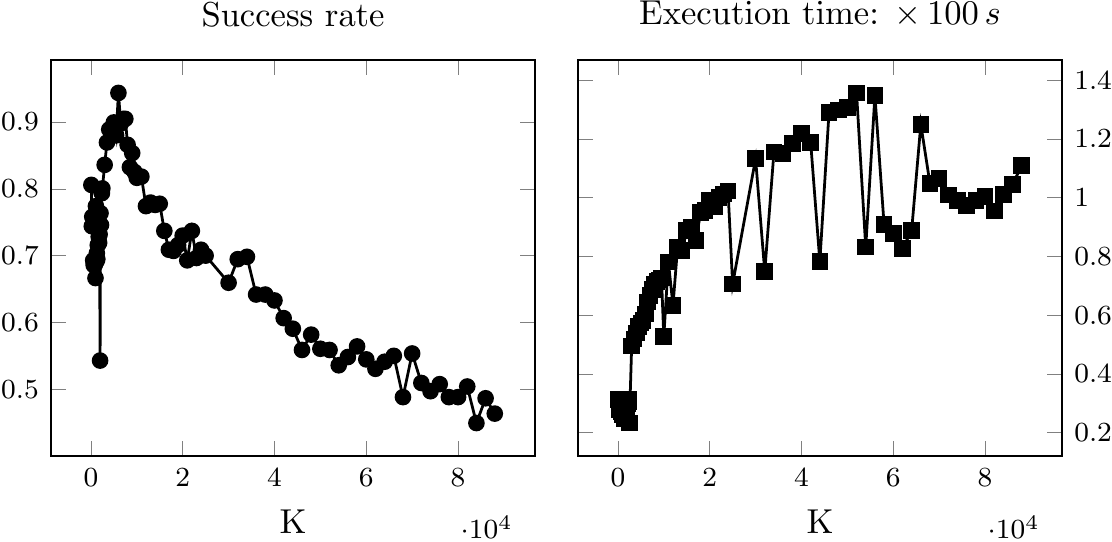}
\caption{Rate of success with respect to $K$ for solving ($\P_A$)}\label{ParaK}
\end{figure}

This suggests that for each specific problem (defined by the launcher configuration and the terminal conditions) a systematical experiment should be processed to find out the best $K$ value.
For example, we have tested the planar cases with different values of $K$. The success rate and the execution time are plotted with respect to $K$ in Figure \ref{ParaK}.

We see that the value of $K$ should neither be too large nor too small.
From Tables \ref{table_statistic_res1}-\ref{table_statistic_res3}, we observe also that the $\lambda_2$-continuation causes most failures in the non-planar case. The success rate could be possibly improved by adapting the $K$ value.


Table \ref{table_statistic_res1} and Table \ref{table_statistic_res3} compare the multiple and the single shooting method ($N=0$). The multiple shooting method ($N=2$) clearly improves the robustness of the algorithm, without significant increase of the execution time.

\begin{table} [H]
\centering
\caption{Statistical results ($N=0$ and $K=8 \times 10^2$).}
\begin{tabular}{|c|c|c|}
\hline
  & planar & non-planar \\
\hline
 Number of cases & 567 & 1134  \\
\hline
  Rate of success (\%) & 83.95 & 74.96  \\
\hline
 Number of failure cases & & \\
   \hspace{0.5cm}- In $\lambda_1$-continuation & 4 & 10 \\
  \hspace{0.5cm}- In $\lambda_2$-continuation & 29 & 210  \\
  \hspace{0.5cm}- In $\lambda_4$-continuation & 21 & 24  \\
  \hspace{0.5cm}- In $\lambda_5$+$K_p$-continuation & 37  & 40  \\   
\hline
  Average execution time (s) &  &  \\
  - Total & 28.93 & 33.36  \\
  \hspace{0.5cm}- In $\lambda_1$-continuation & 0.47 & 0.57  \\
  \hspace{0.5cm}- In $\lambda_2$-continuation & 1.17 & 1.71  \\
  \hspace{0.5cm}- In $\lambda_4$-continuation & 10.80 & 10.51  \\
  \hspace{0.5cm}- In $\lambda_5$+$K_p$-continuation & 18.17 & 21.56  \\
\hline    
\end{tabular}
\label{table_statistic_res3}
\end{table}

\begin{figure}[H]
\centering
\includegraphics[width=0.7\textwidth]{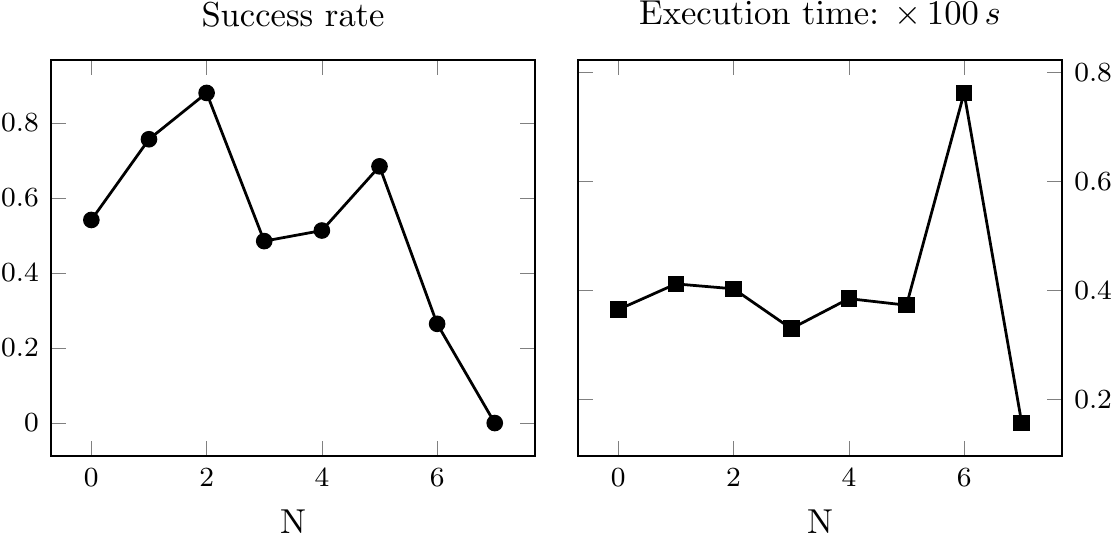}
\caption{Rate of success with respect to $N$ by solving ($\P_A$)}\label{ParaN}
\end{figure}

Figure \ref{ParaN} plots the success rate and the execution time depending of the number of nodes. The test case is the planar maneuver with the regularization parameter $K$ set to $5.5 \times 10^3$. 
The rate of success does not increase monotonically with respect to the number of node points, and the execution time does not change significantly for $N$ less than 6. When $N \geq 6$, the success rate decreases quickly and equals to zero when $N=7$.
When the number of unknowns for the shooting method becomes too large, the domain of convergence of a Newton-type method reduces which finally leads to lower rate of success.

\subsubsection{Reorientation Maneuver of a launch vehicle}
Along multi-burn ascent trajectories, the control (Euler angles) exhibit jumps at the stage separations (see for example \cite[Figure 3]{LuForbes}).
In this case, a reorientation maneuver is necessary to follow the optimal thrust direction. For this reason, we apply the above algorithm as well to the maneuver problem of the upper stages of the launch vehicles.

Opposite to the airborne launch vehicle's pull-up maneuvers, these reorientation maneuvers are in general three-dimensional and of lower magnitude.
They occur at high altitudes (typically higher than 50 km since a sufficiently low dynamic pressure is required to ensure the separation safety) and high velocity (since the first stage has already delivered a large velocity increment).

The maneuver occurs in vacuum so that no state constraints apply. Finding the minimum time maneuver corresponds to solving the problem ($\P_S$).

In the example, we set the system parameters in \eqref{sys_full} to $a = 20$, $b = 0.2$, which approximate an Ariane-like launcher.
The initial conditions \eqref{ic_cmtc} are
$$
r_{x 0}=100\,km, \quad r_{y 0}=r_{z 0}=0, \quad v_{0} = 5000\, m/s, \quad \theta_{v 0}=30^\circ,
$$
$$
\psi_{v0}=0^\circ, \quad \theta_0=40^\circ, \quad \psi_0=\phi_0=0, \quad \omega_{x 0}=\omega_{y 0}=0,
$$
and the final conditions \eqref{fc_cmtc} are
$$\theta_f= 60^\circ, \quad \psi_f=10^\circ, \quad \phi_f=0, \quad \omega_{x f}=\omega_{y f}=0.$$

\begin{figure} [H]
\centering
\includegraphics[width=0.9\textwidth]{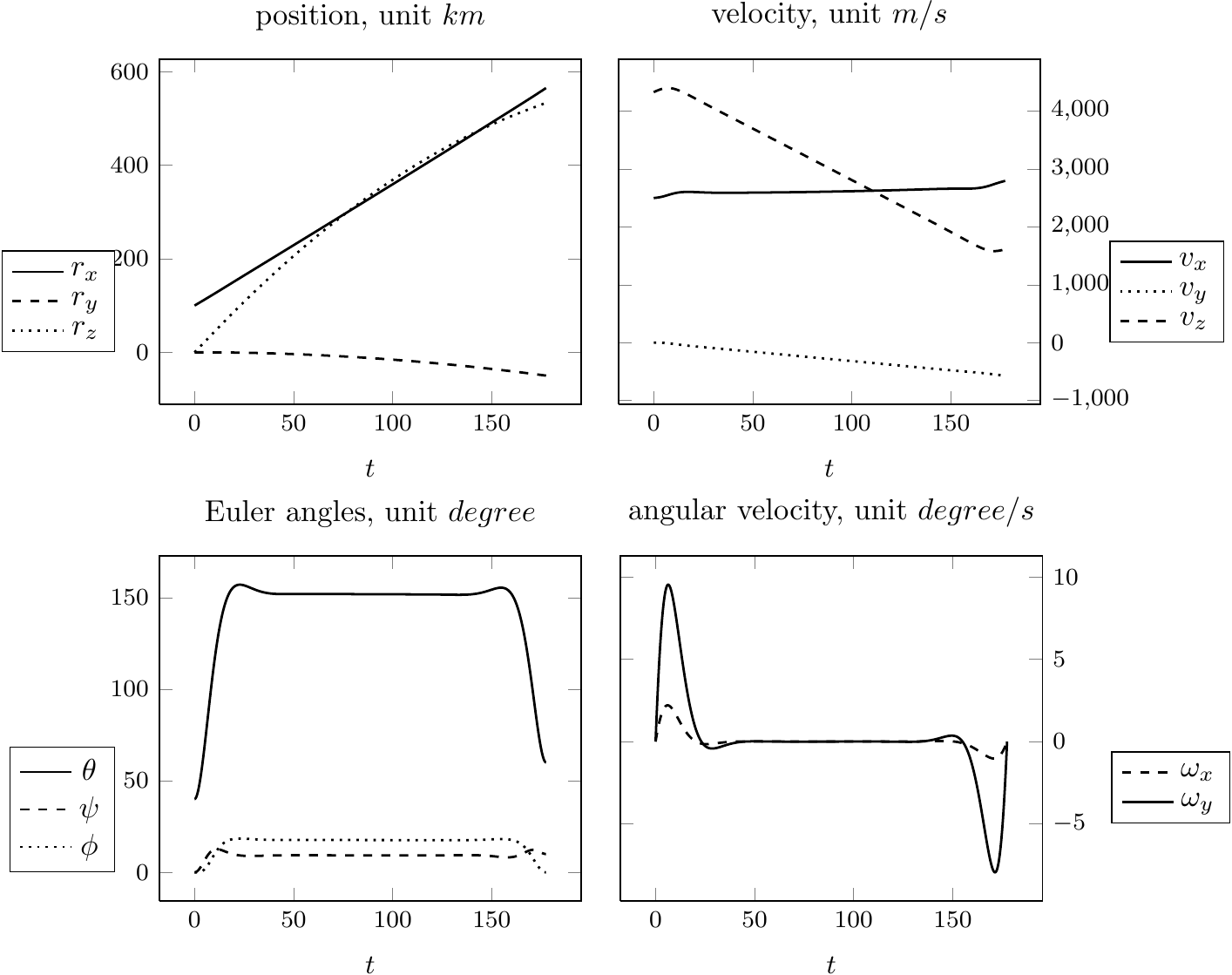}   
\caption{Time history of state variable $x(t)$ for a reorientation maneuver.}
\label{state_mul}
\end{figure}
\begin{figure} [H]
\centering
 \includegraphics[scale=1.0]{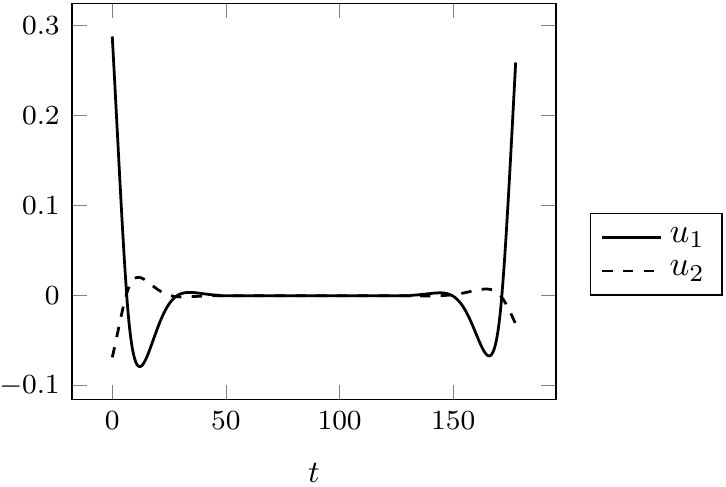}  
\caption{Time history of control variable $u(t)$ for a reorientation maneuver.}  
\label{control_mul}
\end{figure}

The multiple shooting method is applied with four node points. On Figures \ref{state_mul} and  \ref{control_mul}, we report the components of state and control variables.
We observe that, when $t \in [32,145]\, s$, the control is quasi null, and the attitude angles take the solution values of the zero order problem ($\P_0^H$):  $\theta=151.5^\circ \approx \theta^\ast=151.57^\circ$, $\psi=8.6^\circ \approx \psi^\ast=8.85^\circ$.
The regularization term $K \int_0^{t_f} \|u\|^2 dt$ in the cost functional yields a continuous control plotted on Figure \ref{control_mul} and avoids chattering
For this application case the regularization parameter moves from $(1-\lambda_3)K=8\times 10^4$ ($\lambda = 0$) until $(1-\lambda_3)K = 240$ ($\lambda_3 = 0.997$) and the computing time is about $110 \, s$.

The maneuver duration $t_f$ is about $175\ s$ due to the large direction change required on the velocity. During a real flight the velocity direction change is much smaller and the maneuver takes at most a few seconds.
Our purpose when presenting this ``unrealistic'' case is rather to show that the proposed algorithm is robust in a large range of system configurations and terminal conditions.

\bigskip

\section{Applications to Trajectory Optimization}
\label{sec_ATO}
The previous section was devoted to an ascent trajectory application. The example dealt with the pull-up maneuver of an airborne launch vehicle just after its release from the carrier.
This section gives a brief overview of optimal geometric control and continuation techniques applied to other mission categories, namely orbital transfer and atmospheric reentry.

\subsection{Orbital Transfer Problems}
The orbital transfer problem consists in steering the engine from an initial orbit to a final one while minimizing either the duration or the consumption. This problem has been widely studied in the literature, and the solution algorithms involve direct methods as well as indirect methods. The reader is referred to \cite{Betts1998} and \cite{CHT2012} for a list of methods and references.

Our aim is here to recall how geometric optimal control theory and numerical continuation methods can help solving such problems. The dynamics is modelled by the controlled Kepler equations
\begin{align*}
\ddot{r}(t) &=  - r(t) \frac{\mu}{\|r(t)\|^3} + \frac{T(t)}{m(t)}, \\
\dot{m}(t)& = -\beta \|T(t)\|,\\
\| T(t)\| &\leq T_{max},
\end{align*}
where $r(\cdot)$ is the position of the spacecraft, $\mu$ is the gravitational constant of the planet, $T(\cdot) \leq T_{max}$ is the bounded thrust, and $m(\cdot)$ is the mass with $\beta$ a constant depending on the specific impulse of the engine.

Controllability properties ensuring the feasibility of the problem have been studied in \cite{BonnardCT,BFT}, based on the analysis of Lie algebra generated by the vector fields of the system.

The minimum time low thrust transfer is addressed for example in \cite{CGN}. It is observed that the domain of convergence of the Newton-type method in the shooting problem becomes smaller when the maximal thrust decreases.
Therefore, a natural continuation process consists in starting with larger values of the maximal thrust and then decreasing step by step the maximal thrust. In \cite{CGN}, the authors started with the maximal thrust $T_{max}=60 \, N$ and achieved the continuation up to $T_{max} = 0.14 \, N$.

The minimum fuel consumption orbit transfer problem has also been widely studied. With the cost functional $\int_0^{t_f} \| T(t)\| dt$, the problem is more difficult than minimizing the time, since the optimal control is discontinuous. In \cite{Gergaud,Gergaud2}, the authors propose a continuation on the cost functional, starting from the minimum-energy problem. The cost functional is defined by 
$$\int_0^{t_f} \left( (1-\lambda)\| T(t)\|^2 + \lambda \| T(t)\| \right) dt,
$$where $\lambda \in [0,1]$ is the homotopy parameter. When $\lambda=0$ (minimum-energy), the control derived from the PMP is continuous and the shooting problem is thus easier to solve. The authors prove the existence of a zero path from $\lambda=0$ and to $\lambda=1$. This continuation approach is later applied in \cite{CCC2016} for studying the $L^1$-minimization of trajectory optimization problem. 
Such minimum-fuel low-thrust transfers are very important for deep space explorations, since all the propellant must be carried on board by the satellite.
Similar continuation procedures have also been applied to the well-known Goddard's problem, and to its three-dimensional variants (\cite{BMT2008,BLM2009}). The possible singular arcs (along which the norm of the thrust is neither zero nor maximal) have thus been analyzed and numerically computed.

Continuation procedures are also valuable for high-thrust orbital transfer problems. In \cite{CHT2012}, the authors proposed a continuation procedure starting from a flat Earth model with constant gravity. The variable gravity and the Earth curvature are introduced step by step by homotopy parameters. The theoretically analysis of the flat Earth model shows that the solution structure consists in a single boost followed by a coast arc. This helps solving the starting problem in a direct way, before coming back by continuation to the real round Earth problem. The round Earth solution exhibits a different solution structure (boost – coast – boost) which appears progressively along the continuation process.

\subsection{Atmospheric Reentry Problem}
An atmospheric reentry typically begins at an altitude of $120 km$ and ends with a landing phase. The final landing phase until the touchdown is generally studied apart and it is highly dependent on the mission specifications (ground or sea landing, manned or unmanned flight, etc). 
The so-called atmospheric leg aims at reducing the vehicle energy before the final landing phase. No fuel is used and the braking has to be fully achieved by aerodynamics while satisfying the state constraints, in particular on the thermal flux.
The final conditions specify a target position at a low altitude, typically less than $20\ km$.

The vehicle is considered as a glider submitted to the gravity and the aerodynamical forces, the control $u$ being the bank angle and possibly the angle of attack.
The optimal control problem consists thus in steering the space shuttle from given entry conditions to targeted final conditions while minimizing the total heat and satisfying state constraints on the thermal flux, the normal acceleration, and the dynamic pressure.
We refer the readers to \cite{BTreentry,T03} for a formulation of this problem.

The control $u$ acts on the lift force orientation, changing simultaneously the descent rate and the heading angle.

A practical guidance strategy consists in following the constraint boundaries, successively : thermal flux, normal acceleration, and dynamic pressure.
This strategy does not care about the cost functional and it is therefore not optimal.
Applying the Pontryagin Maximum Principle with state constraints is not promising due to a narrow domain of convergence of the shooting method. Finding a correct guess for the initial adjoint vector proves quite difficult.
Therefore direct methods are generally preferred for this atmospheric reentry problem (see, e.g., \cite{Betts1998,Betts2010,Pesch}).

Here we recall two alternative approaches to address the problem by indirect methods.

The first approach is to analyze the control system using geometric control theory.
For example, in \cite{BFLT,BFT2005,T03}, a careful analysis of the control system provides a precise description of the optimal trajectory. The resulting problem reduction makes it tractable by the shooting method.
More precisely, the control system is rewritten as a single-input control-affine system in dimension three under some reasonable assumptions. Local optimal syntheses are derived from extending existing results in geometric optimal control theory.
Based on perturbation arguments, this local nature of the optimal trajectory is then used to provide an approximation of the optimal trajectory for the full problem in dimension six, and finally simple approximation methods are developed to solve the problem.

A second approach is to use the continuation method.
For example, in \cite{Hermant}, the problem is solved by a shooting method, and a continuation is applied on the maximal value of the thermal flux.
It is shown in \cite{BH2008,Hermant2010} that under some appropriate assumptions, the change in the structure of the trajectory is regular, i.e., when a constraint becomes active along the continuation, only one boundary arc appears.
Nevertheless it is possible that an infinite number of boundary arcs appear (see, e.g., \cite{Robbins}). This phenomenon is possible when the constraint is of order three at least.
By using a properly modified continuation procedure, the reentry problem was solved in \cite{Hermant} and the results of \cite{BFT2005} were retrieved.

\section{Conclusion}
The aim of this article was to show how to apply techniques of geometric optimal control and numerical continuation to aerospace problems. After an overview of space transportation missions,  some classical techniques of optimal control have been recalled, including Pontryagin Maximum Principle, first and second-order optimality conditions, and conjugate time theory. Techniques of geometric optimal control have then been recalled, such as higher-order optimality conditions and singular controls.

A quite difficult problem has illustrated in detail how to design an efficient solution method with the help of geometric optimal control tools and continuation methods. Other applications in space trajectory optimization have also been recalled.

Though geometric optimal control and numerical continuation provide a nice way to design efficient approaches for many aerospace applications, the answer to ``how to select a reasonably simple problem for the continuation procedure'' for general optimal control problems remains open.
A deep understanding of the system dynamics is necessary to devise a simple problem that is ``physically'' sufficiently close to the original problem, while being numerically suited to initiate a continuation procedure.

In practice, many problems remain difficult due to the complexity of real-life models. In general, a compromise should be found between the complexity of the model under consideration and the choice of an adapted numerical method.

As illustrated by the example of airborne launch vehicles, many state and/or control constraints should also be considered in a real-life problem, and such constraints makes the problem much more difficult. For the airborne launch problem a penalization method combined with the previous geometric analysis proves satisfying.
But this approach has to be customized to the specific problem under consideration. A challenging task is then to combine an adapted numerical approach with a thorough geometric analysis in order to get more information on the optimal synthesis.
We refer the readers to \cite{Trelatsurvey} for a summary of open challenges in aerospace applications.

\paragraph{Acknowledgment.}
The second author acknowledges the support by FA9550-14-1-0214 of the EOARD-AFOSR.


\end{document}